\pdfoutput=1
\RequirePackage{silence}
\documentclass[a4paper,11pt]{amsart}
\usepackage[hmarginratio={1:1},vmarginratio={1:1},lmargin=60.0pt,tmargin=60.0pt]{geometry}

\allowdisplaybreaks

\usepackage[utf8]{inputenc}

\usepackage{latexsym,exscale,mathtools,textcomp}
\usepackage{amssymb,amsmath,amsthm,amsfonts,mathrsfs,bbm,enumitem,stmaryrd,dsfont}
\usepackage[table]{xcolor}
\usepackage{graphicx}
\usepackage{mathtools}
\usepackage{ytableau}
\usepackage{relsize}
\usepackage{xfrac}
\usepackage{caption}
\usepackage{float}
\usepackage{xparse}
\usepackage{graphbox}
\usepackage{verbatim}

\usepackage{scalefnt}

\usepackage{dynkin-diagrams}

\usepackage{array}
\newcolumntype{C}{>{$}c<{$}}

\definecolor{mygray}{gray}{0.6}
\definecolor{mygraydark}{gray}{0.4}
\definecolor{mygraylight}{gray}{0.85}
\definecolor{spinach}{RGB}{46,139,87}
\definecolor{tomato}{RGB}{255,99,71}
\definecolor{orchid}{RGB}{143,40,194}
\definecolor{neon}{RGB}{77,77,255}
\definecolor{pumpkin}{RGB}{224,180,80}
\definecolor{citron}{RGB}{190,180,90}

\definecolor{lava}{RGB}{207,16,32}
\definecolor{cream}{RGB}{255,253,208}
\definecolor{verdigris}{RGB}{67,179,174}
\definecolor{Black}{RGB}{0,0,0}
\definecolor{mydarkblue}{RGB}{10,10,170}
\definecolor{darkspinach}{RGB}{20,70,20}
\definecolor{darktomato}{RGB}{155,40,30}
\definecolor{darkorchid}{RGB}{50,10,100}
\definecolor{darklava}{RGB}{150,8,16}

\usepackage{todonotes}

\usepackage{enumitem}
\setlist[enumerate]{itemsep=0.15cm,label=\emph{\upshape(\alph*)}}
\setlist[enumerate,2]{itemsep=0.15cm,label=\emph{\upshape(\roman*)}}
\setlist[enumerate,3]{itemsep=0.15cm,label=\emph{\upshape(\Alph*)}}

\let\emph\relax
\DeclareTextFontCommand{\emph}{\bfseries\em}

\newcommand{\placeholder}{{}_{-}}

\renewcommand{\dots}{\text{...}}
\newcommand{\mystrut}{\rule[-0.2\baselineskip]{0pt}{1.00\baselineskip}}

\let\<=\langle
\let\>=\rangle

\newcommand{\acts}{\centerdot}

\renewcommand{\dots}{\text{...}}

\DeclarePairedDelimiterX{\set}[1]{\{}{\}}{\setargs{#1}}
\NewDocumentCommand{\setargs}{>{\SplitArgument{1}{|}}m}{\setargsaux#1}
\NewDocumentCommand{\setargsaux}{mm}
{\IfNoValueTF{#2}{#1} {#1\,\delimsize|\,\mathopen{}#2}}

\newcommand{\etc}{\text{etc.}}

\newcommand{\C}{\mathbb{C}}
\newcommand{\R}{\mathbb{R}}
\newcommand{\N}{\mathbb{Z}_{\geq 0}}

\newcommand{\Z}{\mathbb{Z}}
\newcommand{\K}{\mathbb{K}}

\renewcommand{\dim}[1][\K]{\mathrm{dim}_{#1}}

\newcommand{\End}{\mathrm{End}}
\newcommand{\Hom}{\mathrm{Hom}}

\newcommand{\HH}{\mathrm{H}}

\newcommand{\id}{id}

\newcommand{\algebra}[1][A]{\mathscr{#1}}

\newcommand{\cellbasis}[1][\mathscr{A}]{B_{\mathscr{A}}}

\newcommand{\sand}[1][\lambda]{\mathscr{H}_{#1}}

\newcommand{\rad}[1][\jcell]{\mathrm{Rad}}

\newcommand{\lcell}{\mathcal{L}}
\newcommand{\rcell}{\mathcal{R}}
\newcommand{\jcell}{\mathcal{J}}
\newcommand{\hcell}{\mathcal{H}}

\newcommand{\jideal}[1][\lambda]{\algebra^{>_{lr}\lambda}}

\newcommand{\oneb}{\mathbbm{1}_{b}}
\newcommand{\onet}{\mathbbm{1}_{t}}

\newcommand{\monoid}[1][M]{{#1}}

\newcommand{\onemon}{{1}}

\newcommand{\sym}[1][n]{{S}_{#1}}
\newcommand{\psym}[1][n]{{pS}_{#1}}

\newcommand{\cob}[1][n]{\mathbf{Cob}_{#1}}
\newcommand{\pacat}[1][n]{\mathbf{Pa}_{#1}}
\newcommand{\robrcat}[1][n]{\mathbf{RoBr}_{#1}}
\newcommand{\rocat}[1][n]{\mathbf{Ro}_{#1}}
\newcommand{\ppacat}[1][n]{\mathbf{pPa}_{#1}}

\newcommand{\brcat}[1][n]{\mathbf{Br}_{#1}}
\newcommand{\obrcat}[1][n]{\mathbf{oBr}_{#1}}
\newcommand{\tlcat}[1][n]{\mathbf{TL}_{#1}}
\newcommand{\mocat}[1][n]{\mathbf{Mo}_{#1}}
\newcommand{\procat}[1][n]{\mathbf{pRo}_{#1}}

\newcommand{\pamon}[1][n]{{Pa}_{#1}}
\newcommand{\robrmon}[1][n]{{RoBr}_{#1}}
\newcommand{\romon}[1][n]{{Ro}_{#1}}
\newcommand{\ppamon}[1][n]{{pPa}_{#1}}

\newcommand{\brmon}[1][n]{{Br}_{#1}}
\newcommand{\obrmon}[1][n]{{oBr}_{#1}}
\newcommand{\tlmon}[1][n]{{TL}_{#1}}
\newcommand{\momon}[1][n]{{Mo}_{#1}}
\newcommand{\promon}[1][n]{{pRo}_{#1}}

\usepackage[all]{xy}
\usepackage{tikz}
\usetikzlibrary{cd}
\usetikzlibrary{decorations}
\usetikzlibrary{decorations.markings}
\usetikzlibrary{decorations.pathreplacing}
\usetikzlibrary{decorations.pathmorphing}
\usetikzlibrary{arrows.meta,shapes,positioning,matrix,calc}
\usetikzlibrary{shapes.callouts}
\usetikzlibrary{tqft}
\tikzset{
anchorbase/.style={baseline={([yshift=#1]current bounding box.center)}},
anchorbase/.default={-0.5ex},
tinynodes/.style={font=\tiny,text height=0.25ex,text depth=0.05ex},
smallnodes/.style={font=\scriptsize,text height=0.75ex,text depth=0.15ex},
mor/.style={line width=0.75,color=black,fill=cream},
mor2/.style={line width=0.75,color=black,fill=tomato},
mor3/.style={line width=0.75,color=black,fill=spinach},
usual/.style={line width=1.2,color=black},
crossline/.style={preaction={draw=white,line width=5.0pt,-},preaction={draw=black,line width=0.9pt,-}},
dot/.style = {
decoration={markings,
post length=0.25mm,
pre length=0.25mm,
mark=at position #1 with {\node[circle,radius=0.15cm,inner sep=-1.2pt,color=black,fill=black]{};}
},
postaction={decorate}
},
dot/.default=1,
dotg/.style = {
decoration={markings,
post length=0.25mm,
pre length=0.25mm,
mark=at position #1 with {\node[circle,radius=0.5cm,inner sep=-2.5pt,color=neon,fill=neon]{};}
},
postaction={decorate}
},
dotg/.default=0.5,
}
\tikzstyle directed=[postaction={decorate,decoration={markings,
mark=at position #1 with {\arrow[line width=0.25mm, black]{>}}}}]

\usepackage{pgfplots}
\pgfplotsset{compat=1.18}

\usepackage{aliascnt,etoolbox}
\def\NewTheorem#1{%
\newaliascnt{#1}{equation}%
\newtheorem{#1}[#1]{#1}%
\aliascntresetthe{#1}%
\expandafter\def\csname #1autorefname\endcsname{#1}%
}
\def\equationautorefname~#1\null{(#1)\null}

\numberwithin{equation}{subsection}

\NewTheorem{Proposition}
\NewTheorem{Theorem}
\NewTheorem{Corollary}
\AtEndEnvironment{Corollary}{\null\hfill$\square$}%
\NewTheorem{Lemma}
\NewTheorem{Conjecture}
\NewTheorem{Speculation}
\NewTheorem{Observation}
\NewTheorem{Assumption}
\NewTheorem{Condition}
\theoremstyle{definition}
\NewTheorem{Definition}
\AtEndEnvironment{Definition}{\null\hfill$\Diamond$}%
\NewTheorem{Classification Problem}
\AtEndEnvironment{Classification Problem}{\null\hfill$\Diamond$}%
\NewTheorem{Notation}
\AtEndEnvironment{Notation}{\null\hfill$\Diamond$}%
\NewTheorem{Example}
\AtEndEnvironment{Example}{\null\hfill$\Diamond$}%
\NewTheorem{Examples}
\AtEndEnvironment{Examples}{\vskip-10mm\null\hfill$\Diamond$}%

\theoremstyle{remark}
\NewTheorem{Remark}
\AtEndEnvironment{Remark}{\null\hfill$\Diamond$}%
\NewTheorem{Question}
\AtEndEnvironment{Question}{\null\hfill$\Diamond$}%

\usepackage[stretch=16]{microtype}
\usepackage[T1]{fontenc}
\usepackage{geometry}
\usepackage{graphicx}
\usepackage{footmisc}
\usepackage{datetime}
\usepackage[sort&compress,numbers]{natbib}
\usepackage{tikz}
\usepackage{tikz-3dplot}
\usetikzlibrary{automata,positioning,math,svg.path,tikzmark,arrows,trees,external,decorations.pathreplacing,calligraphy,patterns,patterns.meta,shapes.geometric}
\usepackage{tikz-cd}
\usepackage{float}
\usepackage{subcaption}
\usepackage{appendix}
\usepackage{faktor}
\usetikzlibrary{decorations.markings}
\usepackage{etoolbox}
\usepackage{pdflscape}
\usepackage{afterpage}
\usepackage{rotating}
\usepackage{totcount}
\usepackage{suffix}
\usepackage{enumitem}
\usepackage{multirow}
\usepackage{tabularray}
\usepackage{multicol}
\usepackage{scalerel}
\usepackage{array}
\usepackage[numbered,framed]{matlab-prettifier}

\usepackage{amsmath} 
\usepackage{amsfonts} 
\usepackage{mathdots}
\usepackage{amssymb}
\usepackage{physics} 
\usepackage{mathrsfs} 
\usepackage{mathtools}
\usepackage{amsthm}
\usepackage{thmtools}
\usepackage{mathdots}
\usepackage{bbm}
\usepackage{pifont} 
\usepackage{esint}
\usepackage{mleftright} 
\usepackage{dsfont} 

\newcounter{proofs}
\newcounter{proofLevel}

\newcounter{mcases}
\newcounter{msubcases}[mcases]
\renewcommand{\themsubcases}{\themcases.\arabic{msubcases}}
\newcounter{msubsubcases}[msubcases]
\renewcommand{\themsubsubcases}{\themsubcases.\arabic{msubsubcases}}
\newcounter{msubsubsubcases}[msubsubcases]
\renewcommand{\themsubsubsubcases}{\themsubsubcases.\arabic{msubsubsubcases}}
\newcounter{msubsubsubsubcases}[msubsubsubcases]
\renewcommand{\themsubsubsubsubcases}{\themsubsubsubcases.\arabic{msubsubsubsubcases}}
\newcommand{\case}{}

\newcommand{\nomspace}{\vspace{-3.25ex plus 1ex minus .2ex}} 

\makeatletter
\newenvironment{mcases}
{%
\setcounter{mcases}{0}
\stepcounter{proofLevel}%
\renewcommand{\case}{%
\stepcounter{proofs}%
\refstepcounter{mcases}%
\paragraph*{Case~\themcases}%
\hypertarget{\theproofs}{}%
}%
\addtolength{\leftskip}{\parindent}
}
{%
\addtocounter{proofLevel}{-1}%
\vspace{3.25ex plus 1ex minus .2ex}%
\addtolength{\leftskip}{-\parindent}%
\aftergroup\@afterindentfalse
\aftergroup\@afterheading
}%

\newenvironment{msubcases}
{%
\setcounter{msubcases}{0}
\addtocounter{proofLevel}{2}%
\renewcommand{\case}{%
\stepcounter{proofs}%
\refstepcounter{msubcases}%
\paragraph*{Case~\themsubcases}%
\hypertarget{\theproofs}{}%
}%
\addtolength{\leftskip}{\parindent}
}
{%
\addtocounter{proofLevel}{-2}%
\vspace{3.25ex plus 1ex minus .2ex}%
\addtolength{\leftskip}{-\parindent}
\aftergroup\@afterindentfalse 
\aftergroup\@afterheading 
}%

\newcommand{\by}[1]{\quad\text{(#1)}}

\newcommand{\key}[2][1]{
\quad%
\setlength{\fboxrule}{0.5pt}%
\setlength{\fboxsep}{0.4em}%
\ifnum#1=1
\fbox{%
\begin{tblr}{
columns = {c},
rows = {mode=math}
}%
\SetCell[c=2]{c}\text{\underline{Key}} & \\[1.5ex]%
#2
\end{tblr}%
}%
\else
\ifnum#1=2
\fbox{%
\begin{tblr}{
columns = {c},
rows = {mode=math}
}%
\SetCell[c=4]{c}\text{\underline{Key}} & & & \\[1.5ex]%
#2
\end{tblr}%
}%
\else
\fbox{%
\begin{tblr}{
columns = {c},
rows = {mode=math}
}%
\SetCell[c=6]{c}\text{\underline{Key}} & & & & & \\[1.5ex]%
#2
\end{tblr}%
}%
\fi%
\fi%
\,%
}

\let\oldoverrightarrow\overrightarrow
\renewcommand{\overrightarrow}[1]{\oldoverrightarrow{#1\rule{0pt}{1.9ex}}}

\mleftright 

\let\originalmiddle=\middle
\def\middle#1{\mathrel{}\originalmiddle#1\mathrel{}}

\makeatletter
\newsavebox{\@brx}
\newcommand{\llangle}[1][]{\savebox{\@brx}{\(\m@th{#1\langle}\)}%
\mathopen{\copy\@brx\mkern2mu\kern-0.9\wd\@brx\usebox{\@brx}}}
\newcommand{\rrangle}[1][]{\savebox{\@brx}{\(\m@th{#1\rangle}\)}%
\mathclose{\copy\@brx\mkern2mu\kern-0.9\wd\@brx\usebox{\@brx}}}
\makeatother

\newlength{\LactHeight}

\newlength{\halfEquals}
\usepackage[hypertexnames=false]{hyperref}
\usepackage{bookmark}
\hypersetup{
pdftoolbar=true,
pdfmenubar=true,
pdffitwindow=false,
pdfstartview={FitH},
pdftitle={Generalized diagram categories and monoids, and their representations},
pdfauthor={Matthias Fresacher, Willow Stewart and Daniel Tubbenhauer},
pdfsubject={},
pdfcreator={Matthias Fresacher, Willow Stewart and Daniel Tubbenhauer},
pdfproducer={Matthias Fresacher, Willow Stewart and Daniel Tubbenhauer},
pdfkeywords={},
pdfnewwindow=true,
colorlinks=true,
linkcolor=mydarkblue,
citecolor=teal,
filecolor=magenta,
urlcolor=orchid,
linkbordercolor=lava,
citebordercolor=teal,
urlbordercolor=orchid,
linktocpage=true
}

\def\makeautorefname#1#2{\csdef{#1autorefname}{#2}}

\makeautorefname{section}{Section}%
\makeautorefname{subsection}{Section}%
\makeautorefname{subsubsection}{Section}%

\usepackage[capitalise,nameinlink,noabbrev]{cleveref}
\crefname{equation}{}{}
\crefname{enumi}{}{}
\crefname{corol}{Corollary}{Corollaries}
\Crefname{corol}{Corollary}{Corollaries}
\crefname{corolTOC}{Corollary}{Corollaries}
\Crefname{corolTOC}{Corollary}{Corollaries}
\crefname{task}{Task}{Tasks}
\Crefname{task}{Task}{Tasks}
\crefname{question}{Question}{Questions}
\Crefname{question}{Question}{Questions}
\crefname{comment}{Comment}{Comments}
\Crefname{comment}{Comment}{Comments}
\crefname{ass}{Assumption}{Assumptions}
\Crefname{ass}{Assumption}{Assumptions}
\crefname{assTOC}{Assumption}{Assumptions}
\Crefname{assTOC}{Assumption}{Assumptions}
\crefname{defi}{Definition}{Definitions}
\Crefname{defi}{Definition}{Definitions}
\crefname{chpDefi}{Chapter Definition}{Chapter Definitions}
\Crefname{chpDefi}{Chapter Definition}{Chapter Definitions}
\crefname{prop}{Proposition}{Propositions}
\Crefname{prop}{Proposition}{Propositions}
\crefname{propNoTOC}{Proposition}{Propositions}
\Crefname{propNoTOC}{Proposition}{Propositions}
\crefname{thm}{Theorem}{Theorems}
\Crefname{thm}{Theorem}{Theorems}
\crefname{thmNoTOC}{Theorem}{Theorems}
\Crefname{thmNoTOC}{Theorem}{Theorems}
\crefname{con}{Condition}{Conditions}
\Crefname{con}{Condition}{Conditions}
\crefname{stat}{Statement}{Statements}
\Crefname{stat}{Statement}{Statements}
\crefname{nota}{Notation}{Notations}
\Crefname{nota}{Notation}{Notations}
\crefname{mcases}{Case}{Cases}
\Crefname{mcases}{Case}{Cases}
\crefalias{msubcases}{mcases}
\crefalias{msubsubcases}{mcases}
\crefalias{msubsubsubcases}{mcases}
\crefalias{msubsubsubsubcases}{mcases}

\let\oldcref\cref
\renewcommand{\cref}[1]{%
  \begingroup\hypersetup{linkcolor=black}\oldcref{#1}\endgroup
}
\let\oldCref\Cref
\renewcommand{\Cref}[1]{%
  \begingroup\hypersetup{linkcolor=black}\oldCref{#1}\endgroup
}

\usepackage{fancyvrb}
\SaveVerb{webpage}'staff.cdms.westernsydney.edu.au/~mfresacher/'

\title[Generalized diagram categories and monoids, and their representations]{Generalized diagram categories and monoids, and their representations}
\author[M. Fresacher, W. Stewart, D. Tubbenhauer]{Matthias Fresacher, Willow Stewart, Daniel Tubbenhauer}
\address{M.F.: Western Sydney University, Centre for Research in Mathematics and Data Science, Locked Bag 1797, Penrith NSW 2751, Australia, \href{https://staff.cdms.westernsydney.edu.au/~mfresacher/}{\UseVerb[fontfamily=cmr]{webpage}}, \href{https://orcid.org/0000-0003-0677-3701}{ORCID: 0000-0003-0677-3701}}
\email{M.Fresacher@westernsydney.edu.au}
\address{W.S.: The University of Sydney, School of Mathematics and Statistics F07, Office Carslaw 807, NSW 2006, Australia, \href{https://www.maths.usyd.edu.au/ut/people?who=W_Stewart}{www.maths.usyd.edu.au/ut/people?who=W\_Stewart}, \href{https://orcid.org/0009-0000-2854-2256}{ORCID: 0009-0000-2854-2256}}
\email{willow.stewart@sydney.edu.au}
\address{D.T.: The University of Sydney, School of Mathematics and Statistics F07, Office Carslaw 827, NSW 2006, Australia, \href{http://www.dtubbenhauer.com}{www.dtubbenhauer.com}, \href{https://orcid.org/0000-0001-7265-5047}{ORCID 0000-0001-7265-5047}}
\email{daniel.tubbenhauer@sydney.edu.au}

\begin{document}

\begin{abstract}
Classical diagram categories and monoids, including the Temperley--Lieb, Brauer, and partition cases, arise as special instances of the category of two dimensional cobordisms and admit additional twists that produce a large new family of diagram categories and monoids. In this paper we introduce this family and develop a unified approach to their representation theory.
\end{abstract}

\subjclass[2020]{Primary: 20G05, 05E10, Secondary: 20M30, 05E16}
\keywords{Diagram categories, monoid/semigroup representations, cellular algebras, Schur--Weyl duality}

\maketitle

\tableofcontents

\section{Introduction}\label{intro}

\emph{Diagram categories/algebras/monoids} package algebraic operations into pictures.
In their simplest form, a morphism is a collection of strands and nodes drawn in a rectangle:
composition is stacking, and the monoidal product is juxtaposition.
\begin{gather*}
\begin{tikzpicture}[anchorbase]
\draw[usual] (0,0) to[out=45,in=180] (0.75,0.20) to[out=0,in=135] (1.5,0);
\draw[usual] (0.5,0) to[out=90,in=180] (0.75,0.1) to[out=0,in=90] (1,0);
\draw[usual] (0,0.5) to[out=270,in=180] (0.25,0.3) to[out=0,in=270] (0.5,0.5);
\draw[usual] (1,0.5) to[out=270,in=180] (1.25,0.3) to[out=0,in=270] (1.5,0.5);
\end{tikzpicture}
\circ
\begin{tikzpicture}[anchorbase]
\draw[usual] (0,0) to (0,0.5);
\draw[usual] (1.5,0) to (0.5,0.5);
\draw[usual] (0.5,0) to[out=90,in=180] (0.75,0.2) to[out=0,in=90] (1,0);
\draw[usual] (1,0.5) to[out=270,in=180] (1.25,0.3) to[out=0,in=270] (1.5,0.5);
\end{tikzpicture}
=
\begin{tikzpicture}[anchorbase]
\draw[usual] (0,0) to (0,0.5);
\draw[usual] (1.5,0) to (0.5,0.5);
\draw[usual] (0.5,0) to[out=90,in=180] (0.75,0.2) to[out=0,in=90] (1,0);
\draw[usual] (1,0.5) to[out=270,in=180] (1.25,0.3) to[out=0,in=270] (1.5,0.5);
\draw[usual] (0,0.5) to[out=45,in=180] (0.75,0.7) to[out=0,in=135] (1.5,0.5);
\draw[usual] (0.5,0.5) to[out=90,in=180] (0.75,0.6) to[out=0,in=90] (1,0.5);
\draw[usual] (0,1) to[out=270,in=180] (0.25,0.8) to[out=0,in=270] (0.5,1);
\draw[usual] (1,1) to[out=270,in=180] (1.25,0.8) to[out=0,in=270] (1.5,1);
\end{tikzpicture}
,\quad
\begin{tikzpicture}[anchorbase]
\draw[usual] (0,0) to[out=45,in=180] (0.75,0.20) to[out=0,in=135] (1.5,0);
\draw[usual] (0.5,0) to[out=90,in=180] (0.75,0.1) to[out=0,in=90] (1,0);
\draw[usual] (0,0.5) to[out=270,in=180] (0.25,0.3) to[out=0,in=270] (0.5,0.5);
\draw[usual] (1,0.5) to[out=270,in=180] (1.25,0.3) to[out=0,in=270] (1.5,0.5);
\end{tikzpicture}
\otimes
\begin{tikzpicture}[anchorbase]
\draw[usual] (0,0) to (0,0.5);
\draw[usual] (1.5,0) to (0.5,0.5);
\draw[usual] (0.5,0) to[out=90,in=180] (0.75,0.2) to[out=0,in=90] (1,0);
\draw[usual] (1,0.5) to[out=270,in=180] (1.25,0.3) to[out=0,in=270] (1.5,0.5);
\end{tikzpicture}
=
\begin{tikzpicture}[anchorbase]
\draw[usual] (0,0) to[out=45,in=180] (0.75,0.20) to[out=0,in=135] (1.5,0);
\draw[usual] (0.5,0) to[out=90,in=180] (0.75,0.1) to[out=0,in=90] (1,0);
\draw[usual] (0,0.5) to[out=270,in=180] (0.25,0.3) to[out=0,in=270] (0.5,0.5);
\draw[usual] (1,0.5) to[out=270,in=180] (1.25,0.3) to[out=0,in=270] (1.5,0.5);
\draw[usual] (2,0) to (2,0.5);
\draw[usual] (3.5,0) to (2.5,0.5);
\draw[usual] (2.5,0) to[out=90,in=180] (2.75,0.2) to[out=0,in=90] (3,0);
\draw[usual] (3,0.5) to[out=270,in=180] (3.25,0.3) to[out=0,in=270] (3.5,0.5);
\end{tikzpicture}
,
\\[0.2cm]
\begin{tikzpicture}[anchorbase]
\draw[usual] (0,0) to (0,0.5);
\draw[usual] (1.5,0) to (0.5,0.5);
\draw[usual] (0.5,0) to[out=90,in=180] (0.75,0.2) to[out=0,in=90] (1,0);
\draw[usual] (1,0.5) to[out=270,in=180] (1.25,0.3) to[out=0,in=270] (1.5,0.5);
\end{tikzpicture}
\circ
\begin{tikzpicture}[anchorbase]
\draw[usual] (0,0) to[out=45,in=180] (0.75,0.20) to[out=0,in=135] (1.5,0);
\draw[usual] (0.5,0) to[out=90,in=180] (0.75,0.1) to[out=0,in=90] (1,0);
\draw[usual] (0,0.5) to[out=270,in=180] (0.25,0.3) to[out=0,in=270] (0.5,0.5);
\draw[usual] (1,0.5) to[out=270,in=180] (1.25,0.3) to[out=0,in=270] (1.5,0.5);
\end{tikzpicture}
=
\begin{tikzpicture}[anchorbase]
\draw[usual] (0,0) to[out=45,in=180] (0.75,0.20) to[out=0,in=135] (1.5,0);
\draw[usual] (0.5,0) to[out=90,in=180] (0.75,0.1) to[out=0,in=90] (1,0);
\draw[usual] (0,0.5) to[out=270,in=180] (0.25,0.3) to[out=0,in=270] (0.5,0.5);
\draw[usual] (1,0.5) to[out=270,in=180] (1.25,0.3) to[out=0,in=270] (1.5,0.5);
\draw[usual] (0,0.5) to (0,1);
\draw[usual] (1.5,0.5) to (0.5,1);
\draw[usual] (0.5,0.5) to[out=90,in=180] (0.75,0.7) to[out=0,in=90] (1,0.5);
\draw[usual] (1,1) to[out=270,in=180] (1.25,0.8) to[out=0,in=270] (1.5,1);
\end{tikzpicture}
,\quad
\begin{tikzpicture}[anchorbase]
\draw[usual] (0,0) to (0,0.5);
\draw[usual] (1.5,0) to (0.5,0.5);
\draw[usual] (0.5,0) to[out=90,in=180] (0.75,0.2) to[out=0,in=90] (1,0);
\draw[usual] (1,0.5) to[out=270,in=180] (1.25,0.3) to[out=0,in=270] (1.5,0.5);
\end{tikzpicture}
\otimes
\begin{tikzpicture}[anchorbase]
\draw[usual] (0,0) to[out=45,in=180] (0.75,0.20) to[out=0,in=135] (1.5,0);
\draw[usual] (0.5,0) to[out=90,in=180] (0.75,0.1) to[out=0,in=90] (1,0);
\draw[usual] (0,0.5) to[out=270,in=180] (0.25,0.3) to[out=0,in=270] (0.5,0.5);
\draw[usual] (1,0.5) to[out=270,in=180] (1.25,0.3) to[out=0,in=270] (1.5,0.5);
\end{tikzpicture}
=
\begin{tikzpicture}[anchorbase]
\draw[usual] (0,0) to (0,0.5);
\draw[usual] (1.5,0) to (0.5,0.5);
\draw[usual] (0.5,0) to[out=90,in=180] (0.75,0.2) to[out=0,in=90] (1,0);
\draw[usual] (1,0.5) to[out=270,in=180] (1.25,0.3) to[out=0,in=270] (1.5,0.5);
\draw[usual] (2,0) to[out=45,in=180] (2.75,0.20) to[out=0,in=135] (3.5,0);
\draw[usual] (2.5,0) to[out=90,in=180] (2.75,0.1) to[out=0,in=90] (3,0);
\draw[usual] (2,0.5) to[out=270,in=180] (2.25,0.3) to[out=0,in=270] (2.5,0.5);
\draw[usual] (3,0.5) to[out=270,in=180] (3.25,0.3) to[out=0,in=270] (3.5,0.5);
\end{tikzpicture}
,
\\[0.2cm]
\begin{tikzpicture}[anchorbase]
\draw[usual] (0,0) to[out=45,in=180] (0.75,0.20) to[out=0,in=135] (1.5,0);
\draw[usual] (0.5,0) to[out=90,in=180] (0.75,0.1) to[out=0,in=90] (1,0);
\draw[usual] (0,0.5) to[out=270,in=180] (0.25,0.3) to[out=0,in=270] (0.5,0.5);
\draw[usual] (1,0.5) to[out=270,in=180] (1.25,0.3) to[out=0,in=270] (1.5,0.5);
\end{tikzpicture}
\circ
\begin{tikzpicture}[anchorbase]
\draw[usual] (0,0) to[out=45,in=180] (0.75,0.20) to[out=0,in=135] (1.5,0);
\draw[usual] (0.5,0) to[out=90,in=180] (0.75,0.1) to[out=0,in=90] (1,0);
\draw[usual] (0,0.5) to[out=270,in=180] (0.25,0.3) to[out=0,in=270] (0.5,0.5);
\draw[usual] (1,0.5) to[out=270,in=180] (1.25,0.3) to[out=0,in=270] (1.5,0.5);
\end{tikzpicture}
=
\begin{tikzpicture}[anchorbase]
\draw[usual] (0,0) to[out=45,in=180] (0.75,0.20) to[out=0,in=135] (1.5,0);
\draw[usual] (0.5,0) to[out=90,in=180] (0.75,0.1) to[out=0,in=90] (1,0);
\draw[usual] (0,0.5) to[out=270,in=180] (0.25,0.3) to[out=0,in=270] (0.5,0.5);
\draw[usual] (1,0.5) to[out=270,in=180] (1.25,0.3) to[out=0,in=270] (1.5,0.5);
\draw[usual] (0,0.5) to[out=45,in=180] (0.75,0.7) to[out=0,in=135] (1.5,0.5);
\draw[usual] (0.5,0.5) to[out=90,in=180] (0.75,0.6) to[out=0,in=90] (1,0.5);
\draw[usual] (0,1) to[out=270,in=180] (0.25,0.8) to[out=0,in=270] (0.5,1);
\draw[usual] (1,1) to[out=270,in=180] (1.25,0.8) to[out=0,in=270] (1.5,1);
\end{tikzpicture}
,\quad
\begin{tikzpicture}[anchorbase]
\draw[usual] (0,0) to[out=45,in=180] (0.75,0.20) to[out=0,in=135] (1.5,0);
\draw[usual] (0.5,0) to[out=90,in=180] (0.75,0.1) to[out=0,in=90] (1,0);
\draw[usual] (0,0.5) to[out=270,in=180] (0.25,0.3) to[out=0,in=270] (0.5,0.5);
\draw[usual] (1,0.5) to[out=270,in=180] (1.25,0.3) to[out=0,in=270] (1.5,0.5);
\end{tikzpicture}
\otimes
\begin{tikzpicture}[anchorbase]
\draw[usual] (0,0) to[out=45,in=180] (0.75,0.20) to[out=0,in=135] (1.5,0);
\draw[usual] (0.5,0) to[out=90,in=180] (0.75,0.1) to[out=0,in=90] (1,0);
\draw[usual] (0,0.5) to[out=270,in=180] (0.25,0.3) to[out=0,in=270] (0.5,0.5);
\draw[usual] (1,0.5) to[out=270,in=180] (1.25,0.3) to[out=0,in=270] (1.5,0.5);
\end{tikzpicture}
=
\begin{tikzpicture}[anchorbase]
\draw[usual] (0,0) to[out=45,in=180] (0.75,0.20) to[out=0,in=135] (1.5,0);
\draw[usual] (0.5,0) to[out=90,in=180] (0.75,0.1) to[out=0,in=90] (1,0);
\draw[usual] (0,0.5) to[out=270,in=180] (0.25,0.3) to[out=0,in=270] (0.5,0.5);
\draw[usual] (1,0.5) to[out=270,in=180] (1.25,0.3) to[out=0,in=270] (1.5,0.5);
\draw[usual] (2,0) to[out=45,in=180] (2.75,0.20) to[out=0,in=135] (3.5,0);
\draw[usual] (2.5,0) to[out=90,in=180] (2.75,0.1) to[out=0,in=90] (3,0);
\draw[usual] (2,0.5) to[out=270,in=180] (2.25,0.3) to[out=0,in=270] (2.5,0.5);
\draw[usual] (3,0.5) to[out=270,in=180] (3.25,0.3) to[out=0,in=270] (3.5,0.5);
\end{tikzpicture}
.
\end{gather*}

This philosophy sits behind several classical structures in representation theory and low dimensional topology.
They have been around for donkey's years, and many diagram categories have been discovered (often independently) in different guises over the last century,
remaining central objects at the interface of algebra and combinatorics.
A key example is the \emph{Temperley--Lieb category}:
it arises as a Schur--Weyl dual for $\mathfrak{sl}_2$ \cite{RTW-Valenztheorie},
in statistical mechanics (notably the Potts model) \cite{TL},
and in subfactor theory \cite{Jo}, among many other contexts.

Other important examples include
\emph{Brauer-type categories}, related to orthogonal and symplectic symmetries \cite{Br-brauer-algebra-original},
and the \emph{partition category}, connected to symmetric groups and their interpolations \cite{Jones,Martin,De-cat-st}.
These objects have been studied from many viewpoints, including structural approaches,
combinatorial methods, and tensor-categorical or Schur--Weyl perspectives.

The diagram categories we focus on in this paper are subcategories of the partition category
(top to bottom, left column first; details are given in \autoref{S:SchurWeyl}):
planar partition, Motzkin, Temperley--Lieb, planar rook, planar symmetric (trivial), partition,
rook Brauer, Brauer, rook, and symmetric. They are summarized in the following table. (This table is a cheat-sheet; each row is proved later in \autoref{S:SchurWeyl}, with parameters specified there.)
\begin{gather*}
\scalebox{0.98}{\begin{tabular}{c|c|c|c|c||c|c|c|c|c}
\arrayrulecolor{tomato}
Symbol & Diagrams & $G$ & $V$ & $\Z$?
& Symbol & Diagrams & $G$ & $V$ & $\Z$?
\\
\hline
\hline
$\ppamon[\mathbf{c}](t)$ & \begin{tikzpicture}[anchorbase]
\draw[usual] (0.5,0) to[out=90,in=180] (1.25,0.45) to[out=0,in=90] (2,0);
\draw[usual] (0.5,0) to[out=90,in=180] (1,0.35) to[out=0,in=90] (1.5,0);
\draw[usual] (0.5,1) to[out=270,in=180] (1,0.55) to[out=0,in=270] (1.5,1);
\draw[usual] (1.5,1) to[out=270,in=180] (2,0.55) to[out=0,in=270] (2.5,1);
\draw[usual] (0,0) to (0,1);
\draw[usual] (2.5,0) to (2.5,1);
\draw[usual,dot] (1,0) to (1,0.2);
\draw[usual,dot] (1,1) to (1,0.8);
\draw[usual,dot] (2,1) to (2,0.8);
\end{tikzpicture} & $SL2$ & $\C^2{\otimes}\C^2$ & Y
& $\pamon[\mathbf{c}](t)$ & \begin{tikzpicture}[anchorbase]
\draw[usual] (0.5,0) to[out=90,in=180] (1.25,0.45) to[out=0,in=90] (2,0);
\draw[usual] (0.5,0) to[out=90,in=180] (1,0.35) to[out=0,in=90] (1.5,0);
\draw[usual] (0,1) to[out=270,in=180] (0.75,0.55) to[out=0,in=270] (1.5,1);
\draw[usual] (1.5,1) to[out=270,in=180] (2,0.55) to[out=0,in=270] (2.5,1);
\draw[usual] (0,0) to (0.5,1);
\draw[usual] (1,0) to (1,1);
\draw[usual] (2.5,0) to (2.5,1);
\draw[usual,dot] (2,1) to (2,0.8);
\end{tikzpicture} & $S(t)$ & $\C^t$ & ?
\\
\hline
$\momon[\mathbf{c}](t)$ & \begin{tikzpicture}[anchorbase]
\draw[usual] (0.5,0) to[out=90,in=180] (1.25,0.5) to[out=0,in=90] (2,0);
\draw[usual] (1,0) to[out=90,in=180] (1.25,0.25) to[out=0,in=90] (1.5,0);
\draw[usual] (2,1) to[out=270,in=180] (2.25,0.75) to[out=0,in=270] (2.5,1);
\draw[usual] (0,0) to (1,1);
\draw[usual,dot] (2.5,0) to (2.5,0.2);
\draw[usual,dot] (0,1) to (0,0.8);
\draw[usual,dot] (0.5,1) to (0.5,0.8);
\draw[usual,dot] (1.5,1) to (1.5,0.8);
\end{tikzpicture} & $SL2$ & $\C^2{\oplus}\C$ & Y
& $\robrmon[\mathbf{c}](t)$ & \begin{tikzpicture}[anchorbase]
\draw[usual] (1,0) to[out=90,in=180] (1.25,0.25) to[out=0,in=90] (1.5,0);
\draw[usual] (1,1) to[out=270,in=180] (1.75,0.55) to[out=0,in=270] (2.5,1);
\draw[usual] (0,0) to (0.5,1);
\draw[usual] (2.5,0) to (2,1);
\draw[usual,dot] (0.5,0) to (0.5,0.2);
\draw[usual,dot] (2,0) to (2,0.2);
\draw[usual,dot] (0,1) to (0,0.8);
\draw[usual,dot] (1.5,1) to (1.5,0.8);
\end{tikzpicture} & $Ot,SPt$ & $\C^t$ & Y${}^\ast$
\\
\hline
$\tlmon[\mathbf{c}](t)$ & \begin{tikzpicture}[anchorbase]
\draw[usual] (0.5,0) to[out=90,in=180] (1.25,0.5) to[out=0,in=90] (2,0);
\draw[usual] (1,0) to[out=90,in=180] (1.25,0.25) to[out=0,in=90] (1.5,0);
\draw[usual] (0,1) to[out=270,in=180] (0.25,0.75) to[out=0,in=270] (0.5,1);
\draw[usual] (2,1) to[out=270,in=180] (2.25,0.75) to[out=0,in=270] (2.5,1);
\draw[usual] (0,0) to (1,1);
\draw[usual] (2.5,0) to (1.5,1);
\end{tikzpicture} & $SL2$ & $\C^2$ & Y
& $\brmon[\mathbf{c}](t)$ & \begin{tikzpicture}[anchorbase]
\draw[usual] (0.5,0) to[out=90,in=180] (1.25,0.45) to[out=0,in=90] (2,0);
\draw[usual] (1,0) to[out=90,in=180] (1.25,0.25) to[out=0,in=90] (1.5,0);
\draw[usual] (0,1) to[out=270,in=180] (0.75,0.55) to[out=0,in=270] (1.5,1);
\draw[usual] (1,1) to[out=270,in=180] (1.75,0.55) to[out=0,in=270] (2.5,1);
\draw[usual] (0,0) to (0.5,1);
\draw[usual] (2.5,0) to (2,1);
\end{tikzpicture} & $Ot,SPt$ & $\C^t$ & Y${}^\ast$
\\
\hline
$\promon[\mathbf{c}](t)$ & \begin{tikzpicture}[anchorbase]
\draw[usual] (0,0) to (0.5,1);
\draw[usual] (0.5,0) to (1,1);
\draw[usual] (2,0) to (1.5,1);
\draw[usual] (2.5,0) to (2.5,1);
\draw[usual,dot] (1,0) to (1,0.2);
\draw[usual,dot] (1.5,0) to (1.5,0.2);
\draw[usual,dot] (0,1) to (0,0.8);
\draw[usual,dot] (2,1) to (2,0.8);
\end{tikzpicture} & $GL2$ & $\C{\oplus}\C^{det}$ & Y
& $\romon[\mathbf{c}](t)$ & \begin{tikzpicture}[anchorbase]
\draw[usual] (0,0) to (1,1);
\draw[usual] (0.5,0) to (0,1);
\draw[usual] (2,0) to (2,1);
\draw[usual] (2.5,0) to (0.5,1);
\draw[usual,dot] (1,0) to (1,0.2);
\draw[usual,dot] (1.5,0) to (1.5,0.2);
\draw[usual,dot] (1.5,1) to (1.5,0.8);
\draw[usual,dot] (2.5,1) to (2.5,0.8);
\end{tikzpicture} & $GLt$ & $\C^t{\oplus}\C$ & Y
\\
\hline
$1$ & \begin{tikzpicture}[anchorbase]
\draw[usual] (0,0) to (0,1);
\draw[usual] (0.5,0) to (0.5,1);
\draw[usual] (1,0) to (1,1);
\draw[usual] (1.5,0) to (1.5,1);
\draw[usual] (2,0) to (2,1);
\draw[usual] (2.5,0) to (2.5,1);
\end{tikzpicture} & $SL2$ & $\C$ & Y
& $\sym[\mathbf{c}](t)$ & \begin{tikzpicture}[anchorbase]
\draw[usual] (0,0) to (1,1);
\draw[usual] (0.5,0) to (0,1);
\draw[usual] (1,0) to (1.5,1);
\draw[usual] (1.5,0) to (2.5,1);
\draw[usual] (2,0) to (2,1);
\draw[usual] (2.5,0) to (0.5,1);
\end{tikzpicture} & $GLt$ & $\C^t$ & Y
\end{tabular}}
\end{gather*}
We refer to these as the \emph{classical diagram categories}.
They admit Schur--Weyl dualities with a group (or group scheme) $G$ acting on tensor powers $V^{\otimes n}$
of a fixed $G$-representation $V$.

There are two complementary ways to study their representation theory.
On the one hand, one can exploit Schur--Weyl duality and transfer problems to the group side.
(Historically, this direction is in a sense opposite to the original viewpoint of Schur \cite{Sc-ClassicSchurWeyl}
and Brauer \cite{Br-brauer-algebra-original}, who used the dual algebras to access representation theory of groups).
On the other hand, Brown \cite{Br-gen-matrix-algebras} advocated a unifying internal approach,
which in modern language leads to 
\emph{cellularity} \cite{GrLe-cellular}, or more precisely, \emph{sandwich cellularity}. This, in plain language, boils down to a factorization into ``shirt-belt-pants'' (here for Brauer diagrams):
\begin{gather*}
\begin{tikzpicture}[anchorbase]
\draw[usual] (0,1) to[out=270,in=180] (0.25,0.75) 
to[out=0,in=270] (0.5,1);
\draw[usual] (1,1) to[out=270,in=180] (1.75,0.75) to[out=0,in=270] (2.5,1);
\draw[usual] (0.5,0) to[out=90,in=180] (1,0.25) 
to[out=0,in=90] (1.5,0);
\draw[usual] (2.5,0) to[out=90,in=180] (2.75,0.25) 
to[out=0,in=90] (3,0);
\draw[usual] (0,0) to[out=90,in=180] (0.5,0.5) to (3,0.5) to[out=0,in=270] (3.5,1);
\draw[usual] (1,0) to (1.5,1);
\draw[usual] (2,0) to (2,1);
\end{tikzpicture}
=
\underbrace{\begin{tikzpicture}[anchorbase]
\draw[usual] (0,0) to[out=270,in=180] (0.25,-0.25) 
to[out=0,in=270] (0.5,0);
\draw[usual] (1,0) to[out=270,in=180] (1.75,-0.5) to[out=0,in=270] (2.5,0);
\draw[usual] (1.5,0) to (1.5,-1);
\draw[usual] (2,0) to (2,-1);
\draw[usual] (3,0) to (3,-1);
\end{tikzpicture}}_{\text{shirt}}
\;\circ\;
\underbrace{\begin{tikzpicture}[anchorbase]
\draw[usual] (0,0) to (1,1);
\draw[usual] (0.5,0) to (0,1);
\draw[usual] (1,0) to (0.5,1);
\end{tikzpicture}}_{\text{belt}}
\;\circ\;
\underbrace{\begin{tikzpicture}[anchorbase]
\draw[usual] (-2.5,0) to (-2.5,1);
\draw[usual] (-1.5,0) to (-1.5,1);
\draw[usual] (-2,0) to[out=90,in=180] (-1.5,0.25) 
to[out=0,in=90] (-1,0);
\draw[usual] (-0.5,0) to (-0.5,1);
\draw[usual] (0,0) to[out=90,in=180] (0.25,0.25) 
to[out=0,in=90] (0.5,0);
\end{tikzpicture}}_{\text{pants}}
,\quad
\begin{aligned}
\begin{tikzpicture}[anchorbase,scale=1]
\draw[line width=0.75,color=black,fill=cream] (0,1) to (0.25,0.5) to (0.75,0.5) to (1,1) to (0,1);
\node at (0.5,0.75){$T$};
\end{tikzpicture}
&\text{ shirt,}
\\
\begin{tikzpicture}[anchorbase,scale=1]
\draw[white,ultra thin] (0,0) to (1,0);
\draw[line width=0.75,color=black,fill=cream] (0.25,0) to (0.25,0.5) to (0.75,0.5) to (0.75,0) to (0.25,0);
\node at (0.5,0.25){$m$};
\end{tikzpicture}
&\text{ belt,}
\\
\begin{tikzpicture}[anchorbase,scale=1]
\draw[line width=0.75,color=black,fill=cream] (0,-0.5) to (0.25,0) to (0.75,0) to (1,-0.5) to (0,-0.5);
\node at (0.5,-0.25){$B$};
\end{tikzpicture}
&\text{ pants.}
\end{aligned}
\end{gather*}
This perspective applies uniformly to many diagram algebras and provides a robust toolset for constructing
and organizing their simple modules; see, for example, \cite{FG,Tu-sandwich,Scrimshaw} for the classical diagram algebras, and \cite{KoXi-affine-cellular,GuWi-almost-cellular,LRMD,MaTu,Liu25,HeTu} for other examples of sandwich structures on diagram categories.

In this paper we enlarge the classical families to a much broader class of diagram categories
(which we call \emph{generalized diagram categories}).
We develop both a Brown-style treatment of their representation theory,
and, whenever available, a Schur--Weyl-style approach.
Along the way we also supply formulations, mild generalizations and random-walk-type statistics of several classical Schur--Weyl dualities
that are often not in the literature.

\subsection{Generalized diagram categories}

Our starting point is the \emph{two dimensional cobordism category} $\cob[\infty]$ (in the usual ``spine'' picture with ``handle equals dot''),
whose morphisms are surfaces built from pairs of pants, caps/cups, and swaps.
Its algebraic models are Frobenius algebras \cite{Dijkgraaf,Abrams}; see, for instance, \cite{Ko-tqfts}.
The reader may keep the following basic generators in mind:
\begin{gather*}
\begin{tikzpicture}[tqft/cobordism/.style={draw},tqft/every lower boundary component/.style={draw=black},anchorbase]
\pic[tqft/reverse pair of pants,genus=1];
\end{tikzpicture}
\leftrightsquigarrow
\begin{tikzpicture}[anchorbase,yscale=-1]
\draw[usual] (0,0) to (0.5,0.5) to (1,0);
\draw[usual,dotg] (0.5,0.5) to (0.5,1);
\end{tikzpicture}
=
\begin{tikzpicture}[anchorbase,yscale=-1]
\draw[usual,dotg=0.25] (0,0) to (0.5,0.5) to (1,0);
\draw[usual] (0.5,0.5) to (0.5,1);
\end{tikzpicture}
=
\begin{tikzpicture}[anchorbase,yscale=-1]
\draw[usual,dotg=0.75] (0,0) to (0.5,0.5) to (1,0);
\draw[usual] (0.5,0.5) to (0.5,1);
\end{tikzpicture}
.
\end{gather*}
Following \cite{Jones,Comes,KhSa-cobordisms}, $\cob[\infty]$ gives rise to the partition category:
by passing to suitable quotients (imposing local relations and, crucially, evaluating closed components) one
recovers many familiar diagram categories.
A key feature of the cobordism viewpoint is that it naturally produces \emph{families}:
one may choose evaluation data controlling how closed surfaces are replaced by scalars,
and one may impose relations on handles.
Already at this level one obtains diagram categories, algebras, and monoids beyond the ``classical'' ones.
This is closely related in spirit to other topological and diagrammatic extensions of standard diagram algebras,
such as blob \cite{MaSa-blob}, decorated Temperley--Lieb \cite{Gr-gen-tl-algebra}, and cyclotomic Hecke algebras \cite{ArKo-hecke-algebra,BrMa-hecke,Ch-gelfandtzetlin}, to name a few.

Concretely, the construction in \autoref{S:DiagCat} produces, in a uniform way, planar and symmetric diagram monoids
(Temperley--Lieb, Motzkin, planar rook, Brauer, rook Brauer, partition, \dots), together with parameterized variants.
The first message of the paper is that these generalizations are not cosmetic:
they genuinely enlarge the landscape, but they do so in a controlled way that can still be studied systematically.

As a \emph{second layer} of generalization:
A different, very flexible way to modify a monoid is to keep track of a statistic that appears under multiplication.
For diagram concatenation, the most basic such statistic counts the number of floating components created in the middle.
Motivated by this, we study \emph{twisted products of monoids} (twists for short), which have been particularly active in the semigroup literature in recent years; see, for example, \cite{east2025twisted}.
Given a monoid $S$, a map $\Phi:S\times S\to \mathbb{N}$ satisfying a cocycle-type identity,
an additive commutative monoid $M$, and a distinguished element $q\in M$,
one forms a twisted product $M\times_{\Phi}^{q} S$ whose multiplication records $\Phi$.
For example, $\Phi$ could be the floating-component count (refined by surface type in the parameterized setting).

At first glance, twistings look representation-theoretically dangerous: they change multiplication and can create new phenomena.
One of the main conclusions of this paper is that, for the diagram monoids under consideration,
a large class of twistings is surprisingly mild from the viewpoint of simple modules.
In particular, most twisted diagram monoids we consider still fall under the umbrella of Brown's unified approach
via sandwich cellularity.

\subsection{What we prove}

The guiding goal is to show that, for a broad and natural class of \emph{generalized} diagram monoids
(coming from parameters and from twistings), the core representation theory is essentially the same
as for their classical counterparts.
\begin{gather*}
\fcolorbox{orchid!50}{spinach!10}{\mystrut\parbox{0.88\linewidth}{\centering
The representation theory is, to a large extent, agnostic to parameters and twists.}}
\end{gather*}
This is made precise through three complementary strands of results.

\medskip
\noindent\textbf{(1) Cellular structure and a uniform representation-theoretic toolkit.}
Diagram algebras are often cellular (or close to it), and this structure is one of the main mechanisms that makes
their simple modules accessible.
Building on sandwich cellularity, we set up in \autoref{S:DiagCat}
a general framework in which our cobordism-quotient diagram categories and monoids admit compatible cellular data.
This provides a common language for indexing simples, comparing parameter choices,
and linking monoid-theoretic structure (Green's relations) to representation theory.

\medskip
\noindent\textbf{(2) A reduction theorem for tight twistings.}
In \autoref{S:Twistings} we specialize the general theory of \cite{east2025twisted} to the diagrammatic twistings relevant here. Upon one condition \autoref{con: exists implies for all idempotents},
for planar diagram monoids (and, more generally, finite monoids with trivial $\mathcal{H}$-classes),
we prove that any tight twisted version has the same simple module indexing and the same simple dimensions
as either the untwisted monoid or the corresponding ``$0$-twisted'' monoid
(\autoref{T:Gram2}; see also the structural analysis in \autoref{SS:Cellular}).
In other words, within this wide class, twistings do not create new simple dimensions: they only decide
whether floating components behave as they do classically or are annihilated.
(We emphasize that this statement is not asserted for so-called loose twistings; see the remarks after \autoref{T:Gram2}.)

\medskip
\noindent\textbf{(3) Schur--Weyl dualities and (expected) dimensions.}
To turn structural results into concrete numerics one needs access to dimensions.
A powerful bridge is provided by Schur--Weyl-type dualities:
commuting actions of a diagram monoid/algebra and a classical group (or Lie algebra)
on a tensor space, generating mutual centralizers.
Such dualities are classical in the symmetric group case \cite{Sc-ClassicSchurWeyl}
and have many diagrammatic incarnations (Temperley--Lieb, Motzkin, Brauer, partition, rook, \dots);
see, for example, \cite{BH-MotzkinAlgebras,LZ-SecondInvOrtho,HaRa-partition-algebras,AM-SchurRook} for more recent references.
In \autoref{S:SchurWeyl} we collect and (re)prove a suite of these dualities in a form that is as uniform as possible,
including statements that are scattered in the literature or appear mostly as folk knowledge.

Once these dualities are in place, dimensions of simple diagrammatic modules can be read off from the ``group side'':
decomposition multiplicities and classical character theory translate into closed formulas.
We carry this out in two regimes.
In the semisimple case (cf. \autoref{S:Resultssemi}) we recover and streamline known formulas and asymptotics
for sums of dimensions (and related size measures).
A crucial input is a \emph{typical highest weight} (or concentration) phenomenon, developed in \autoref{S:TypicalHW}:
for large tensor powers, most of the mass is carried by a narrow window of highest weights,
so that global growth can be captured by a typical truncation without changing the exponential rate.
For Temperley--Lieb this gives the following binomial-type picture of the dimensions of simple representations (the bottom is a log plot):
\begin{gather*}
\begin{tikzpicture}[anchorbase]
\node at (0,0) {\includegraphics[height=5cm]{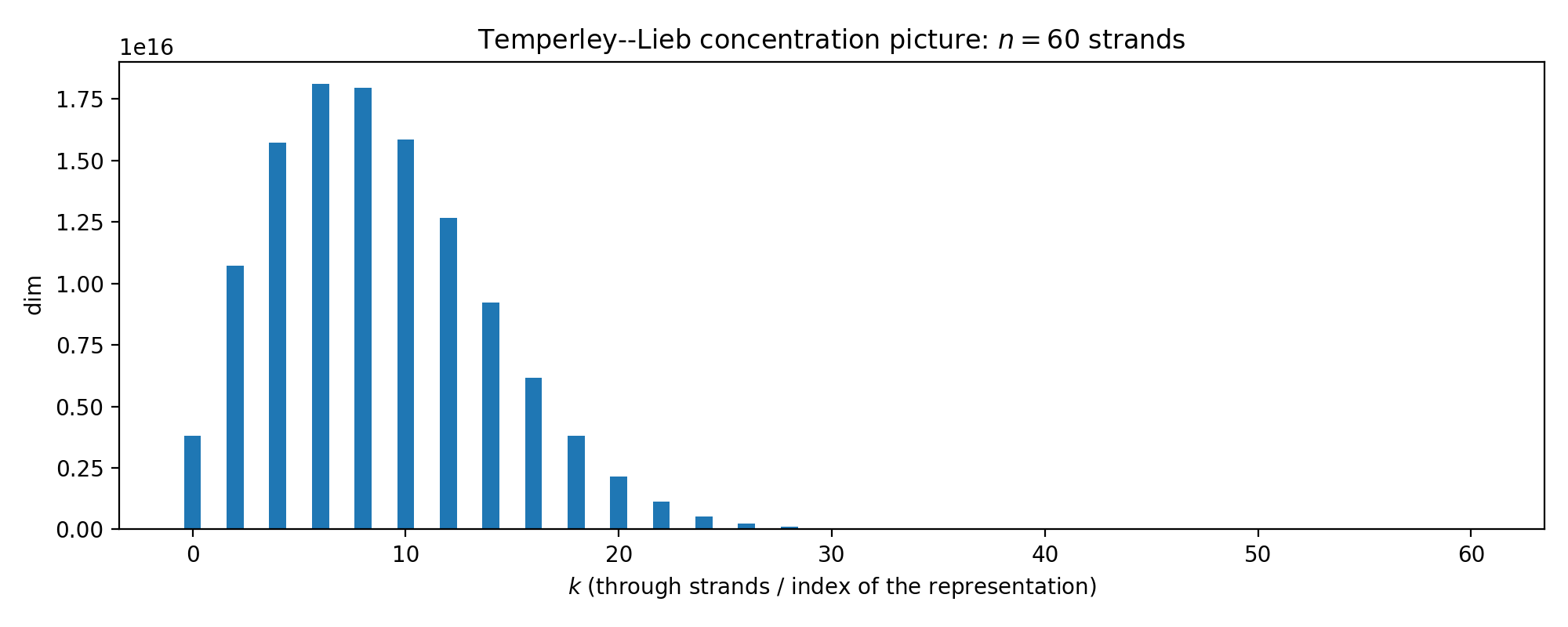}};
\end{tikzpicture}
,\\
\begin{tikzpicture}[anchorbase]
\node at (0,0) {\includegraphics[height=5cm]{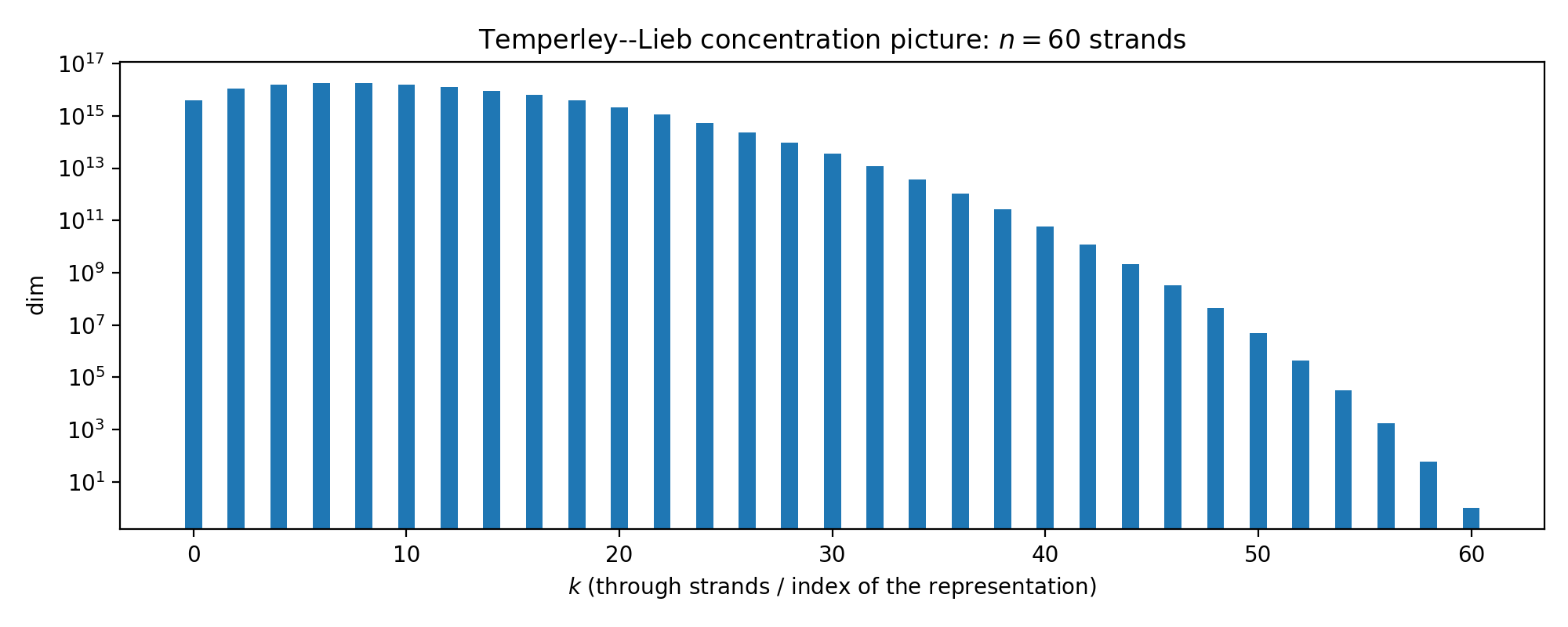}};
\end{tikzpicture}
.
\end{gather*}
In the nonsemisimple case (\autoref{S:Results}) we show that, for the diagram families treated here,
nonsemisimplicity changes the fine structure (extensions, truncations, exceptional small ranks),
but in many cases the dominant growth of simple dimensions remains close to the semisimple prediction.

\noindent\textbf{Acknowledgments.}
We thank James East for many helpful conversations on twistings and insights into his joint work \cite{east2025twisted}. DT is also grateful to Ben Elias for the shirt–belt–pants picture.

WS was supported by the Postgraduate Research Scholarship in Mathematics and Statistics (SC4238). DT is supported by the ARC Future Fellowship FT230100489, and they leave any remaining errors as an exercise for the reader (with no hints).

\section{Generalized diagram monoids: cobordisms}\label{S:DiagCat}

We assume some familiarity with monoidal categories, see e.g. \cite{egno-tensor-2015} for background. Our conventions follow \cite{Tu-qtop}. Let $\mathbb{K}$ be any field.

\subsection{Cobordism categories}

Our starting point is a prototypical class of \emph{diagram categories}: quotients $\cob[{p/q}]$ of the (symmetric monoidal) \emph{cobordism category} $\cob[\infty]$, see e.g. \cite{Ko-tqfts} and \cite{KhSa-cobordisms}. Here the objects are one dimensional compact manifolds (circles) and the morphisms are two dimensional cobordisms (pants), which we will draw using their spines with handles indicated by dots:
\begin{gather*}
\begin{tikzpicture}[tqft/cobordism/.style={draw},tqft/every lower boundary component/.style={draw=black},anchorbase]
\pic[tqft/pair of pants];
\end{tikzpicture}
\leftrightsquigarrow
\begin{tikzpicture}[anchorbase]
\draw[usual] (0,0) to (0.5,0.5) to (1,0);
\draw[usual] (0.5,0.5) to (0.5,1);
\end{tikzpicture}
,\quad
\begin{tikzpicture}[tqft/cobordism/.style={draw},tqft/every lower boundary component/.style={draw=black},anchorbase]
\pic[tqft/cap];
\end{tikzpicture}
\leftrightsquigarrow
\begin{tikzpicture}[anchorbase]
\draw[usual,dot] (0,0) to (0,0.5);
\end{tikzpicture}
,\quad
\begin{tikzpicture}[tqft/cobordism/.style={draw},tqft/every lower boundary component/.style={draw=black},anchorbase]
\pic[tqft/cylinder to prior, at={(0,0)}];
\pic[tqft/cylinder to next, at={(-1,0)}];
\end{tikzpicture}
\leftrightsquigarrow
\begin{tikzpicture}[anchorbase]
\draw[usual] (0,0) to (1,1);
\draw[usual] (1,0) to (0,1);
\end{tikzpicture}
,\quad
\begin{tikzpicture}[tqft/cobordism/.style={draw},tqft/every lower boundary component/.style={draw=black},anchorbase]
\pic[tqft/cylinder,genus=1];
\end{tikzpicture}
\leftrightsquigarrow
\begin{tikzpicture}[anchorbase]
\draw[usual,dotg] (0,0) to (0,1);
\end{tikzpicture}
,\quad
\text{\etc}
\end{gather*}
Composition $\circ$ is stacking and the monoidal product $\otimes$ is juxtaposition with the convention:
\begin{gather*}
E\circ D
=
\begin{tikzpicture}[anchorbase,smallnodes,rounded corners]
\node[rectangle,draw,minimum width=0.5cm,minimum height=0.5cm,ultra thick] at(0,0){\raisebox{-0.05cm}{$D$}};
\node[rectangle,draw,minimum width=0.5cm,minimum height=0.5cm,ultra thick] at(0,0.5){\raisebox{-0.05cm}{$E$}};
\end{tikzpicture}
,\quad
E\otimes D
=
\begin{tikzpicture}[anchorbase,smallnodes,rounded corners]
\node[rectangle,draw,minimum width=0.5cm,minimum height=0.5cm,ultra thick] at(0,0){\raisebox{-0.05cm}{$D$}};
\node[rectangle,draw,minimum width=0.5cm,minimum height=0.5cm,ultra thick] at(-0.5,0){\raisebox{-0.05cm}{$E$}};
\end{tikzpicture}
.
\end{gather*}
The swap cobordism as above gives this category the structure of a \emph{symmetric monoidal category} and the sideways tubes (cap and cup) give it the structure of a \emph{pivotal category}.
\begin{gather*}
\text{Symmetric: }
\begin{tikzpicture}[anchorbase]
\draw[usual] (0,0) to (1,1);
\draw[usual] (1,0) to (0,1);
\end{tikzpicture}
,\quad
\text{pivotal: }
\begin{tikzpicture}[anchorbase]
\draw[usual] (0,0) to[out=90,in=180] (0.5,0.5) to[out=0,in=90] (1,0);
\end{tikzpicture}
=
\begin{tikzpicture}[anchorbase]
\draw[usual] (0,0) to (0.5,0.5) to (1,0);
\draw[usual,dot] (0.5,0.5) to (0.5,1);
\end{tikzpicture}
.
\end{gather*}
We will always use these structures when referring to $\cob[\infty]$ or its quotients and their subcategories.

The quotient categories $\cob[{p/q}]$ depend on the choice of a generating function given as a quotient of two polynomials $p/q$, and we capture them using $\mathbf{a}=\{a_i\mid i\in\N\}$, $\mathbf{b}=\{b_i\mid i\in\Z_{>1}\}$ and $k=\min\{\deg p+1,\deg q\}\in\N$.

This works in the following way:
The coefficients $a_i$ of the Taylor expansion of $p/q$ are used to evaluate closed surfaces (we stop drawing the surfaces from now on):
\begin{gather}\label{Eq:Evaluate}
\begin{tikzpicture}[baseline ={([yshift=-.5ex]current bounding box.center)}, scale=0.2]
\draw (0,0) circle (2cm);
\draw (-2,0) arc (180:360:2 and 0.6);
\draw[dashed] (2,0) arc (0:180:2 and 0.6);
\end{tikzpicture}=a_0,\quad
\begin{tikzpicture}[scale=0.2, baseline ={([yshift=-.5ex]current bounding box.center)}]
\draw (0,0) ellipse (3cm and 2cm);
\draw (-1.5, 0.2) arc (180:360:1.5 and 0.6);
\draw (1,-0.2) arc (0:180:1 and 0.6);
\end{tikzpicture}=a_1,\quad
\begin{tikzpicture}[scale=0.2, baseline ={([yshift=-.5ex]current bounding box.center)}]
\draw (0,0) ellipse (4.8cm and 2cm);
\draw (-3.5, 0.2) arc (180:360:1.5 and 0.6);
\draw (-1,-0.2) arc (0:180:1 and 0.6);
\draw (0.5, 0.2) arc (180:360:1.5 and 0.6);
\draw (3,-0.2) arc (0:180:1 and 0.6);
\end{tikzpicture}=a_2,\quad\dots\leftrightsquigarrow
\begin{tikzpicture}[anchorbase]
\draw[usual,dot=0] (0,0) to (0,0.1);
\draw[usual] (0,0.1) to (0,0.9);
\draw[usual,dot=1] (0,0.9) to (0,1);
\end{tikzpicture}
=a_{0},\quad
\begin{tikzpicture}[anchorbase]
\draw[usual,dot=0] (0,0) to (0,0.1);
\draw[usual,dotg] (0,0.1) to (0,0.9);
\draw[usual,dot=1] (0,0.9) to (0,1);
\end{tikzpicture}
=a_{1},\quad
\begin{tikzpicture}[anchorbase]
\draw[usual,dot=0] (0,0) to (0,0.1);
\draw[usual,dotg=0.25,dotg=0.75] (0,0.1) to (0,0.9);
\draw[usual,dot=1] (0,0.9) to (0,1);
\end{tikzpicture}
=a_{2},\quad
\dots,
\end{gather}
and, after $\K$-linear extension (we will not distinguish the categories and their $\K$-linear extensions in the notation), these relations can be interpreted over any field $\K$ in which $a_i\in\K$ make sense. Write $q=1-b_1x+b_2x^2\mp\dots$, and add the additional relation (the label next to the dots represent their number)
\begin{gather}\label{Eq:Evaluate2}
\begin{tikzpicture}[anchorbase]
\draw[usual,dotg=0.5] (0,0) to (0,0.5)node[left]{$k$} to (0,1);
\end{tikzpicture}
-b_1\cdot 
\begin{tikzpicture}[anchorbase]
\draw[usual,dotg=0.5] (0,0) to (0,0.5)node[left]{$(k-1)$} to (0,1);
\end{tikzpicture}
+b_2\cdot 
\begin{tikzpicture}[anchorbase]
\draw[usual,dotg=0.5] (0,0) to (0,0.5)node[left]{$(k-2)$} to (0,1);
\end{tikzpicture}
-\dots+(-1)^{\deg q}
\begin{tikzpicture}[anchorbase]
\draw[usual,dotg=0.5] (0,0) to (0,0.5)node[left]{$(k-\deg q)$} to (0,1);
\end{tikzpicture}
=0
,
\end{gather}
called the \emph{handle relation}.

\begin{Example}
The case $p=t,q=1-x$ is the \emph{partition category} $\pacat[{t/(1-x)}]$, in which all handles disappear:
\begin{gather*}
\dots
=
\begin{tikzpicture}[anchorbase]
\draw[usual,dot=0] (0,0) to (0,0.1);
\draw[usual,dotg=0.25,dotg=0.75] (0,0.1) to (0,0.9);
\draw[usual,dot=1] (0,0.9) to (0,1);
\end{tikzpicture}
=
\begin{tikzpicture}[anchorbase]
\draw[usual,dotg,dot=0,dot=1] (0,0) to (0,0.5) to (0,1);
\end{tikzpicture}
=
\begin{tikzpicture}[anchorbase]
\draw[usual,dot=0,dot=1] (0,0) to (0,1);
\end{tikzpicture}
=t\cdot\emptyset
,\quad
\begin{tikzpicture}[anchorbase]
\draw[usual,dotg=0.5] (0,0) to (0,0.5) to (0,1);
\end{tikzpicture}
= 
\begin{tikzpicture}[anchorbase]
\draw[usual] (0,0) to (0,0.5) to (0,1);
\end{tikzpicture}
,
\end{gather*}
and $t=a_{0}=a_{1}=\dots\in\K$ is the common value of floating components; this is precisely the category studied in \cite{Jones,Martin,De-cat-st}.
\end{Example}

\begin{Lemma}\label{L:MonoidPara}
For $p\in\{0,1\}[x]$ (a binary polynomial) and $q=1+(-1)^rx^r$ for odd $r\in\N$, $\cob[{p/q}]$ can be viewed as a settheoretic category (i.e. without fixing an underlying field).
\end{Lemma}

\begin{proof}
Under these assumptions \autoref{Eq:Evaluate2} becomes
\begin{gather*}
\begin{tikzpicture}[anchorbase]
\draw[usual,dotg=0.5] (0,0) to (0,0.5)node[left]{$k$} to (0,1);
\end{tikzpicture}
= 
\begin{tikzpicture}[anchorbase]
\draw[usual,dotg=0.5] (0,0) to (0,0.5)node[left]{$(k-r)$} to (0,1);
\end{tikzpicture},
\end{gather*}
and all $a_i\in\{0,1\}$. Now we artificially adjoin a zero to hom-spaces, and, in \autoref{Eq:Evaluate}, interpret $1$ as removing the floating component and $0$ as setting the diagram equal to the zero element.
\end{proof}

\begin{Definition}
We call the choice in \autoref{L:MonoidPara}  
\emph{monoid parameters}. We denote them using $\mathbf{c}=\{\mathbf{a},r\}$ only. The choice $\mathbf{a}=(1,1,\dots)$ and $r=1$ is \emph{classical}.
\end{Definition}

\subsection{Some first generalized diagram monoids}

From this we also obtain subcategories of the partition category that are related to the classical \emph{diagram monoids}, defining their diagrams. And, notably, we obtain generalizations of the diagram categories and monoids.

These monoids are the endomorphism monoids of the respective categories for classical monoid parameters, but we can also define them for all monoid parameters $\mathbf{c}$; the subscript denotes the parameter, while the number in brackets indicates the object $n$, meaning that one takes certain endomorphisms of $n\in\N$ circles $\bullet$):
\begin{enumerate}[label=$\bullet$]

\item The (decorated) \emph{partition monoid} $\pamon[\mathbf{c}](n)$ consisting of all diagrams of partitions of a $2n$-element set (plus dots).

\item The (decorated) \emph{rook Brauer monoid} $\robrmon[{a_0,a_1,r}](n)$ consisting of all diagrams with components of size $1$ or $2$ (plus dots).

\item The (decorated) \emph{Brauer monoid} $\brmon[{a_1,r}](n)$ consisting of all diagrams with components of size $2$ (plus dots).

\item The (decorated) \emph{rook monoid} $\romon[{a_0,r}](n)$ consisting of all diagrams with components of size $1$ or $2$, such that all partitions have at most one component
at the bottom and at most one at the top (plus dots).

\item The (decorated) \emph{symmetric group} $\sym[r](n)$ consisting of all matchings with components of size $2$ (plus dots).

\item \emph{Planar} versions of these: $\ppamon[\mathbf{c}](n)$, $\momon[{a_0,a_1,r}](n)$, $\tlmon[{a_1,r}](n)$, $\promon[{a_0,r}](n)$ and $\psym[r](n)\cong 1$ (the latter denotes the
trivial monoid). The planar rook Brauer monoid is also called the \emph{Motzkin monoid}, the planar Brauer monoid is also known as the \emph{Temperley--Lieb monoid}, and the planar symmetric group is trivial.

\item Additionally, there are oriented versions of these. The one we will very briefly use is the \emph{oriented Brauer monoid} $\obrmon[{a_1,r}](n)$, whose elements are Brauer diagrams equipped with orientations, and their decorated versions.

\end{enumerate}
We can also consider more general diagram categories, namely the corresponding subcategories of 
$\pacat[{p/q}]$ for any choice of polynomials and their endomorphism algebras.

\begin{Remark}
Above, $\pacat[\mathbf{c}]$ and $\ppacat[\mathbf{c}]$ depend on the full set of parameters $\mathbf{c}$, $\robrcat[{a_0,a_1,r}]$ and $\mocat[{a_0,a_1,r}]$ depend only on two ($a_0$ for a sphere and $a_1$ for a torus) and $r$, $\tlcat[{a_1,r}]$ and $\brcat[{a_1,r}]$ depend only on $a_1$ (and similarly for $\obrcat[{a_1,r}]$) and $r$,  
$\procat[{a_0,r}]$ and $\rocat[{a_0,r}]$ depend only on $a_0$ and $r$, and the remaining ones are parameter-free but still have $r$.
\end{Remark}

\begin{Notation}\label{N:ClassicalDiagramMonoids}
Fix $\delta\in\K$. If $r=1$, then denote $\mathbf{M}_{\mathbf{a}} = \mathbf{M}_{\mathbf{a},1}$.
If, further, the parameters $\mathbf{a}$ are such that $a_i = \delta$ for all $i$ and 
$r=1$, then denote $M_{\delta}(n)=M_{\delta/(1-x)}(n)$, and similarly for their categories. These are then the \emph{classical diagram algebras} and \emph{categories}. 
\end{Notation}

\begin{Remark}
The terminology for these objects is not completely uniform in the literature, so we briefly record our conventions. Our partition algebras $\pamon[\delta](n)$ are the usual partition algebras of Jones and Martin. The rook Brauer algebras $\robrmon[\delta](n)$ are also known as \emph{partial Brauer algebras}, and hence $\momon[\delta](n)$ is often called the \emph{planar partial Brauer algebra}.
The rook algebras $\romon[\delta](n)$ are often referred to in semigroup theory as the \emph{symmetric inverse monoids}. The algebra $\tlmon[\delta](n)$ is by far the most well-known one and has been rediscovered several times; it appears under the names \emph{Temperley--Lieb}, \emph{Jones}, \emph{Kauffman}, \emph{Rumer}, and \emph{Rumer--Teller--Weyl algebra} in the literature. Finally, the oriented Brauer algebras $\obrmon[\delta](n)$ are often referred to as \emph{walled Brauer} algebras.
\end{Remark}

\begin{Remark}\label{R:PlanarPartition}
By \cite[(1.5)]{HaRa-partition-algebras}, $pPa_{\delta^2}(n) \cong TL_{\delta}(2n)$ for $\delta \neq 0$. Therefore, our results for $pPa_1(n)$ will be covered in the Temperley--Lieb section. For general $\mathbf{c}$, however, these algebras are not necessarily isomorphic.
\end{Remark}

We will always use the structures in \autoref{L:RigidBraided} below for these categories.

\begin{Lemma}\label{L:RigidBraided}
For $p,q\in\K[x]$, all the above categories are $\K$-linear and (strict) monoidal, and for monoid parameters all of these are (strict) monoidal. Moreover:
\begin{enumerate}
\item All the non-planar ones are symmetric, and 
$\tlcat[{p/q}]$ and $\mocat[{p/q}]$ are braided (under a mild assumption specified below).
\item All categories are pivotal except for
$\rocat[{p/q}]$, $\procat[{p/q}]$, the symmetric category, and its planar version.
\end{enumerate}
\end{Lemma}

\begin{proof}
The universal construction \cite{BHMV} guarantees associativity of the compositions $\circ$ and $\otimes$ as well as the interchange law. The rest follows essentially by definition, except the statement that 
$\tlcat[{p/q}]$ and $\mocat[{p/q}]$ are braided, which follows from the Kauffman bracket formula. That is, find $v\in\K$ with a square root $v^{1/2}$ such that $-v-v^{-1}=a_1$ for $\tlcat[{p/q}]$, and $1-v-v^{-1}=a_1$ for $\mocat[{p/q}]$. That such $v$ exist is the assumption mentioned in the statement of the lemma above. Then
\begin{gather*}
\beta_{\bullet\bullet}=
v^{1/2}\cdot\big(\,
\begin{tikzpicture}[anchorbase]
\draw[usual] (0,0) to (0,0.5);
\draw[usual] (0.5,0) to (0.5,0.5);
\end{tikzpicture}
-
\begin{tikzpicture}[anchorbase]
\draw[usual] (0,0) to (0,0.5);
\draw[usual,dot] (0.5,0) to (0.5,0.2);
\draw[usual,dot] (0.5,0.5) to (0.5,0.3);
\end{tikzpicture}
-
\begin{tikzpicture}[anchorbase]
\draw[usual,dot] (0,0) to (0,0.2);
\draw[usual,dot] (0,0.5) to (0,0.3);
\draw[usual] (0.5,0) to (0.5,0.5);
\end{tikzpicture}
\,\big)
+
v^{-1/2}\cdot\big(\,
\begin{tikzpicture}[anchorbase]
\draw[usual] (0,0) to[out=90,in=90] (0.5,0);
\draw[usual] (0,0.5) to[out=270,in=270] (0.5,0.5);
\end{tikzpicture}
-
\begin{tikzpicture}[anchorbase]
\draw[usual] (0,0) to[out=90,in=90] (0.5,0);
\draw[usual,dot] (0,0.5) to (0,0.3);
\draw[usual,dot] (0.5,0.5) to (0.5,0.3);
\end{tikzpicture}
-
\begin{tikzpicture}[anchorbase]
\draw[usual] (0,0.5) to[out=270,in=270] (0.5,0.5);
\draw[usual,dot] (0,0) to (0,0.2);
\draw[usual,dot] (0.5,0) to (0.5,0.2);
\end{tikzpicture}
\,\big)
+
\begin{tikzpicture}[anchorbase]
\draw[usual] (0,0) to (0.5,0.5);
\draw[usual,dot] (0,0.5) to (0,0.3);
\draw[usual,dot] (0.5,0) to (0.5,0.2);
\end{tikzpicture}
+
\begin{tikzpicture}[anchorbase]
\draw[usual] (0.5,0) to (0,0.5);
\draw[usual,dot] (0.5,0.5) to (0.5,0.3);
\draw[usual,dot] (0,0) to (0,0.2);
\end{tikzpicture}
+(v^{1/2}+v^{-1/2}-1)
\cdot
\begin{tikzpicture}[anchorbase]
\draw[usual,dot] (0,0.5) to (0,0.3);
\draw[usual,dot] (0.5,0) to (0.5,0.2);
\draw[usual,dot] (0.5,0.5) to (0.5,0.3);
\draw[usual,dot] (0,0) to (0,0.2);
\end{tikzpicture}
,
\end{gather*}
defines a braiding on these categories. The displayed diagrams are in $\mocat[{p/q}]$; for $\tlcat[{p/q}]$ one uses the respective diagrams that are part of this category.
\end{proof}

\begin{Remark}
There are several other methods to have even more diagram monoids. For example, as in \cite{Liu25}, one could also keep floating components around and only evaluate them if a fixed number, say $m$, have piled up. Or, as in \cite{KMY}, one could look at monoids ``in between'' the planar and the symmetric monoids. Or quantizations, or more subalgebras as e.g. in \cite{ST-RepGapRigidPlanar}.
\end{Remark}

\subsection{Cellularity}\label{SS:Cellular}

We now assume familiarity with \emph{sandwich cellular algebras} as in \cite{TuVa-handlebody,Tu-sandwich}. See also \cite{MaTu} for a correction of a typo. 

\begin{Remark}
The reader familiar with \emph{Green's relations} (a.k.a. cells) and representation theory of monoids, see e.g. \cite{Gr-structure-semigroups,St-rep-monoid}, will feel at home: in the monoid parameter case, there is no crucial difference between sandwich cellular algebras on the one hand, and Green's relations and the Clifford--Munn--Ponizovski\u{\i} theorem on the other hand.
\end{Remark}

The \emph{diagrammatic anti-involution} ${}^{\ast}$ flips a cobordism upside-down.
We call a diagram a \emph{merge diagram} if it contains only multiplications, counits, a minimal number of crossings, and dots (= handles) on nonthrough components. 
(For readers familiar with \cite{GT-GrowthDiagCat}, this corrects a typo in that paper.)
A \emph{split diagram} is a ${}^{\ast}$-flipped merge diagram. A \emph{dotted permutation diagram} contains only dots and crossings. Here are some examples (flipping the left illustration gives a split diagram):
\begin{gather*}
\text{Merge diagram: }
\begin{tikzpicture}[anchorbase]
\draw[usual] (0,0) to (0.5,0.5) to (1,0);
\draw[usual] (0.5,0.5) to (0.5,1);
\draw[usual] (2,0) to (0.5,0.75);
\draw[usual] (0.5,0) to (1,1);
\draw[usual,dotg=0.5,dot] (-0.5,0) to (-0.5,0.5);
\end{tikzpicture}
,\quad
\text{dotted permutation diagram: }
\begin{tikzpicture}[anchorbase]
\draw[usual,dotg=0.85,dotg=0.95] (0,0) to (3,1);
\draw[usual,dotg=0.85] (1,0) to (1,1);
\draw[usual] (2,0) to (0,1);
\draw[usual] (3,0) to (2,1);
\end{tikzpicture}
.
\end{gather*}
Recalling sandwich cellularity, we have the following picture for $\cob[{p/q}]$:
\begin{gather}\label{Eq:CellsCob}
\begin{tikzpicture}[anchorbase,scale=1]
\draw[line width=0.75,color=black,fill=cream] (0,1) to (0.25,0.5) to (0.75,0.5) to (1,1) to (0,1);
\node at (0.5,0.75){$T$};
\draw[line width=0.75,color=black,fill=cream] (0.25,0) to (0.25,0.5) to (0.75,0.5) to (0.75,0) to (0.25,0);
\node at (0.51,0.25){$m$};
\draw[line width=0.75,color=black,fill=cream] (0,-0.5) to (0.25,0) to (0.75,0) to (1,-0.5) to (0,-0.5);
\node at (0.5,-0.25){$B$};
\end{tikzpicture}
\quad\text{where}\quad
\begin{aligned}
\begin{tikzpicture}[anchorbase,scale=1]
\draw[line width=0.75,color=black,fill=cream] (0,1) to (0.25,0.5) to (0.75,0.5) to (1,1) to (0,1);
\node at (0.5,0.75){$T$};
\end{tikzpicture}
&\text{ a split diagram,}
\\
\begin{tikzpicture}[anchorbase,scale=1]
\draw[white,ultra thin] (0,0) to (1,0);
\draw[line width=0.75,color=black,fill=cream] (0.25,0) to (0.25,0.5) to (0.75,0.5) to (0.75,0) to (0.25,0);
\node at (0.5,0.25){$m$};
\end{tikzpicture}
&\text{ a dotted permutation,}
\\
\begin{tikzpicture}[anchorbase,scale=1]
\draw[line width=0.75,color=black,fill=cream] (0,-0.5) to (0.25,0) to (0.75,0) to (1,-0.5) to (0,-0.5);
\node at (0.5,-0.25){$B$};
\end{tikzpicture}
&\text{ a merge diagram.}
\end{aligned}
\quad
\text{${}^{\ast}$: }
\left(\,
\begin{tikzpicture}[anchorbase,scale=1]
\draw[line width=0.75,color=black,fill=cream] (0,1) to (0.25,0.5) to (0.75,0.5) to (1,1) to (0,1);
\node at (0.5,0.75){$T$};
\draw[line width=0.75,color=black,fill=cream] (0.25,0) to (0.25,0.5) to (0.75,0.5) to (0.75,0) to (0.25,0);
\node at (0.51,0.25){$m$};
\draw[line width=0.75,color=black,fill=cream] (0,-0.5) to (0.25,0) to (0.75,0) to (1,-0.5) to (0,-0.5);
\node at (0.5,-0.25){$B$};
\end{tikzpicture}\,
\right)^{\ast}
=
\begin{tikzpicture}[anchorbase,scale=1]
\draw[line width=0.75,color=black,fill=cream] (0,1) to (0.25,0.5) to (0.75,0.5) to (1,1) to (0,1);
\node at (0.5,0.75){\reflectbox{\rotatebox{180}{$B$}}};
\draw[line width=0.75,color=black,fill=cream] (0.25,0) to (0.25,0.5) to (0.75,0.5) to (0.75,0) to (0.25,0);
\node at (0.51,0.25){\reflectbox{\rotatebox{180}{$m$}}};
\draw[line width=0.75,color=black,fill=cream] (0,-0.5) to (0.25,0) to (0.75,0) to (1,-0.5) to (0,-0.5);
\node at (0.5,-0.25){\reflectbox{\rotatebox{180}{$T$}}};
\end{tikzpicture}
.
\end{gather}
We call the existence of a spanning set with a decomposition as above a \emph{precell structure}, and the respective algebras \emph{presandwich}.
Let us list the presandwich cellular bases for the remaining diagram categories. For the planar partition it is the same but without crossings, and:
\begin{gather}\label{Eq:CellsCob2}
\begin{gathered}
\mocat[{p/q}],\\
\tlcat[{p/q}],\\
\procat[{p/q}]
\end{gathered}
\colon
\begin{aligned}
\begin{tikzpicture}[anchorbase,scale=1]
\draw[line width=0.75,color=black,fill=cream] (0,1) to (0.25,0.5) to (0.75,0.5) to (1,1) to (0,1);
\node at (0.5,0.75){$T$};
\end{tikzpicture}
&\text{ a cup and unit diagram,}
\\
\begin{tikzpicture}[anchorbase,scale=1]
\draw[white,ultra thin] (0,0) to (1,0);
\draw[line width=0.75,color=black,fill=cream] (0.25,0) to (0.25,0.5) to (0.75,0.5) to (0.75,0) to (0.25,0);
\node at (0.5,0.25){$m$};
\end{tikzpicture}
&\text{ a dotted identity,}
\\
\begin{tikzpicture}[anchorbase,scale=1]
\draw[line width=0.75,color=black,fill=cream] (0,-0.5) to (0.25,0) to (0.75,0) to (1,-0.5) to (0,-0.5);
\node at (0.5,-0.25){$B$};
\end{tikzpicture}
&\text{ cap and counit diagram,}
\end{aligned}
\quad
\begin{gathered}
\robrcat[{p/q}],\\
\brcat[{p/q}],\\
\rocat[{p/q}]
\end{gathered}
\colon
\begin{aligned}
\begin{tikzpicture}[anchorbase,scale=1]
\draw[line width=0.75,color=black,fill=cream] (0,1) to (0.25,0.5) to (0.75,0.5) to (1,1) to (0,1);
\node at (0.5,0.75){$T$};
\end{tikzpicture}
&\text{ a cup and unit diagram,}
\\
\begin{tikzpicture}[anchorbase,scale=1]
\draw[white,ultra thin] (0,0) to (1,0);
\draw[line width=0.75,color=black,fill=cream] (0.25,0) to (0.25,0.5) to (0.75,0.5) to (0.75,0) to (0.25,0);
\node at (0.5,0.25){$m$};
\end{tikzpicture}
&\text{ a dotted permutation,}
\\
\begin{tikzpicture}[anchorbase,scale=1]
\draw[line width=0.75,color=black,fill=cream] (0,-0.5) to (0.25,0) to (0.75,0) to (1,-0.5) to (0,-0.5);
\node at (0.5,-0.25){$B$};
\end{tikzpicture}
&\text{ a cap and counit diagram.}
\end{aligned}
\end{gather}
For the symmetric category and its planar counterpart the presandwich structure is trivial. For completeness:
\begin{gather}\label{Eq:CellsCob3}
\obrcat[{p/q}]\colon
\begin{aligned}
\begin{tikzpicture}[anchorbase,scale=1]
\draw[line width=0.75,color=black,fill=cream] (0,1) to (0.25,0.5) to (0.75,0.5) to (1,1) to (0,1);
\node at (0.5,0.75){$T$};
\end{tikzpicture}
&\text{ a cup diagram,}
\\
\begin{tikzpicture}[anchorbase,scale=1]
\draw[white,ultra thin] (0,0) to (1,0);
\draw[line width=0.75,color=black,fill=cream] (0.25,0) to (0.25,0.5) to (0.75,0.5) to (0.75,0) to (0.25,0);
\node at (0.5,0.25){$m$};
\end{tikzpicture}
&\text{ a $\uparrow$ dotted permutation and a $\downarrow$ dotted permutation,}
\\
\begin{tikzpicture}[anchorbase,scale=1]
\draw[line width=0.75,color=black,fill=cream] (0,-0.5) to (0.25,0) to (0.75,0) to (1,-0.5) to (0,-0.5);
\node at (0.5,-0.25){$B$};
\end{tikzpicture}
&\text{ a cap diagram.}
\end{aligned}
\end{gather}

As before, write $\infty$ for $p/q$ in the case we do not impose any relations of the form \autoref{Eq:Evaluate} or \autoref{Eq:Evaluate2}.

\begin{Lemma}\label{L:Cell}
The endomorphism algebras of $\cob[{\infty}]$ are involutive presandwich cellular with precell structure as in \autoref{Eq:CellsCob}. Ditto for the other diagram categories, using the presandwich structure in \autoref{Eq:CellsCob2} and \autoref{Eq:CellsCob3}, and the diagrammatic antiinvolution.
\end{Lemma}

\begin{proof}
Easy and omitted. Hint: Use \cite[Section 1.4.16]{Ko-tqfts}.
\end{proof}

For the next statement we assume familiarity with the usual tableaux combinatorics, see, for example, \cite{Ma-hecke-schur}.
Write \autoref{Eq:Evaluate2} as
\begin{gather*}
\prod_{i=0}^k
\bigg(
\begin{tikzpicture}[anchorbase]
\draw[usual,dotg=0.5] (0,0) to (0,0.5) to (0,1);
\end{tikzpicture}
-u_i\cdot
\begin{tikzpicture}[anchorbase]
\draw[usual] (0,0) to (0,0.5) to (0,1);
\end{tikzpicture}\,
\bigg)
=0
\end{gather*}
for some $u_i$ in the algebraic closure of $\K$. Let $\ell$ be the number of distinct values of $u_i$.

\begin{Theorem}\label{T:Cob}
Let $\K$ be of characteristic $e$, and assume we have fixed $p,q$. The presandwich structure from \autoref{L:Cell} descends to a sandwich cellular basis after removing all diagrams with $\geq k$ dots and:
\begin{enumerate}

\item Every simple $\cob[{p/q}](n)$-representation is associated to some \emph{apex} being the largest cell not annihilating it.

\item The set of apexes of $\cob[{p/q}](n)$ is the subset of $\{0,\dots,n\}$ (the number of strands in the middle) given by
\begin{gather*}
\begin{cases}
\{0,\dots,n\} & \text{if }a_0\neq 0,
\\
\{0,\dots,n\}\setminus\{n-1\} & \text{if }a_0=0\text{ but some }a_r\neq 0,r\leq n,
\\
\{0,\dots,n\}\setminus\{0,n-1\} & \text{if }\mathbf{a}=0.
\end{cases}
\end{gather*}

\item The finite dimensional simple $\cob[{p/q}](n)$-modules of apex $m$, up to equivalence, are indexed by \emph{$e$-restricted $\ell$-multipartitions} of $m$.

\item The dimension of the \emph{cell representation} for an $e$-restricted $k$-multipartition $\lambda$ of $m$ is
\begin{gather*}
\#\{\text{merge diagrams with $m$ top strands}\}
\cdot
\#\{\text{standard tableaux of shape $\lambda$}\},
\\
\#\{\text{merge diagrams with $m$ top strands}\}
=
\sum_{i=0}^{\lfloor k/2\rfloor}
\ell^{k-2i}
\sum_{b=k}^{n}
\begin{Bmatrix}n\\b\end{Bmatrix}
\binom{b}{k}
\binom{k}{2i}
(2i-1)!!
,
\end{gather*}
where $\begin{Bsmallmatrix}n\\t\end{Bsmallmatrix}$ denotes the Stirling numbers of the second kind.

\item If $\cob[{p/q}](n)$ is semisimple, then the cell representations are simple.

\end{enumerate}
\end{Theorem}

\begin{proof}
The basis claim follows from \autoref{L:Cell} and \cite[Proposition 2]{KhSa-cobordisms}.
Most of the other claims are then formal consequences of having a sandwich structure, except (b) and (c).

For (b), note that $n$ is an apex as the identity diagram gives an idempotent. Moreover, we can use
\begin{gather*}
\begin{tikzpicture}[anchorbase]
\draw[usual,dot] (0,0.5) to (0,0.3);
\draw[usual,dot] (0,0) to (0,0.2);
\end{tikzpicture}
\circ
\begin{tikzpicture}[anchorbase]
\draw[usual,dot] (0,0.5) to (0,0.3);
\draw[usual,dot] (0,0) to (0,0.2);
\end{tikzpicture}
=a_0\cdot
\begin{tikzpicture}[anchorbase]
\draw[usual,dot] (0,0.5) to (0,0.3);
\draw[usual,dot] (0,0) to (0,0.2);
\end{tikzpicture}
,
\end{gather*}
to create pseudo-idempotents whenever $a_0$ is not zero, which works for any other number of through strands. Otherwise, as long as there is some through strand,
\begin{gather*}
\begin{tikzpicture}[anchorbase]
\draw[usual] (0,0) to[out=90,in=90] (0.5,0);
\draw[usual] (0.5,0.5) to[out=270,in=270] (1,0.5);
\draw[usual] (1,0) to (0,0.5);
\end{tikzpicture}
\circ
\begin{tikzpicture}[anchorbase]
\draw[usual] (0,0) to[out=90,in=90] (0.5,0);
\draw[usual] (0.5,0.5) to[out=270,in=270] (1,0.5);
\draw[usual] (1,0) to (0,0.5);
\end{tikzpicture}
=
\begin{tikzpicture}[anchorbase]
\draw[usual] (0,0) to[out=90,in=90] (0.5,0);
\draw[usual] (0.5,0.5) to[out=270,in=270] (1,0.5);
\draw[usual] (1,0) to (0,0.5);
\end{tikzpicture}
,\quad
\begin{tikzpicture}[anchorbase]
\draw[usual] (0,0) to[out=90,in=180] (0.5,0.15) to[out=0,in=90] (1,0);
\draw[usual] (0.5,0) to (0.5,0.15);
\draw[usual,dot] (0.5,0.5) to (0.5,0.3);
\draw[usual] (1,0.5) to[out=270,in=270] (1.5,0.5);
\draw[usual] (1.5,0) to[out=90,in=270] (0,0.5);
\end{tikzpicture}
\circ
\begin{tikzpicture}[anchorbase]
\draw[usual] (0,0) to[out=90,in=180] (0.5,0.15) to[out=0,in=90] (1,0);
\draw[usual] (0.5,0) to (0.5,0.15);
\draw[usual,dot] (0.5,0.5) to (0.5,0.3);
\draw[usual] (1,0.5) to[out=270,in=270] (1.5,0.5);
\draw[usual] (1.5,0) to[out=90,in=270] (0,0.5);
\end{tikzpicture}
=
\begin{tikzpicture}[anchorbase]
\draw[usual] (0,0) to[out=90,in=180] (0.5,0.15) to[out=0,in=90] (1,0);
\draw[usual] (0.5,0) to (0.5,0.15);
\draw[usual,dot] (0.5,0.5) to (0.5,0.3);
\draw[usual] (1,0.5) to[out=270,in=270] (1.5,0.5);
\draw[usual] (1.5,0) to[out=90,in=270] (0,0.5);
\end{tikzpicture}
,
\end{gather*}
can be used to create idempotents.
For the remaining cell with zero through strands, say $a_3\neq 0$, then
\begin{gather*}
\begin{tikzpicture}[anchorbase]
\draw[usual] (0,0) to[out=90,in=180] (0.5,0.15) to (1,0.15) to[out=0,in=90] (1.5,0);
\draw[usual] (0.5,0) to (0.5,0.15);
\draw[usual] (1,0) to (1,0.15);
\draw[usual] (0,0.5) to[out=270,in=180] (0.5,0.35) to (1,0.35) to[out=0,in=270] (1.5,0.5);
\draw[usual] (0.5,0.5) to (0.5,0.35);
\draw[usual] (1,0.5) to (1,0.35);
\end{tikzpicture}
\circ
\begin{tikzpicture}[anchorbase]
\draw[usual] (0,0) to[out=90,in=180] (0.5,0.15) to (1,0.15) to[out=0,in=90] (1.5,0);
\draw[usual] (0.5,0) to (0.5,0.15);
\draw[usual] (1,0) to (1,0.15);
\draw[usual] (0,0.5) to[out=270,in=180] (0.5,0.35) to (1,0.35) to[out=0,in=270] (1.5,0.5);
\draw[usual] (0.5,0.5) to (0.5,0.35);
\draw[usual] (1,0.5) to (1,0.35);
\end{tikzpicture}
=
a_3\cdot
\begin{tikzpicture}[anchorbase]
\draw[usual] (0,0) to[out=90,in=180] (0.5,0.15) to (1,0.15) to[out=0,in=90] (1.5,0);
\draw[usual] (0.5,0) to (0.5,0.15);
\draw[usual] (1,0) to (1,0.15);
\draw[usual] (0,0.5) to[out=270,in=180] (0.5,0.35) to (1,0.35) to[out=0,in=270] (1.5,0.5);
\draw[usual] (0.5,0.5) to (0.5,0.35);
\draw[usual] (1,0.5) to (1,0.35);
\end{tikzpicture},
\end{gather*}
is a pseudo idempotent. If no such $a_r$ exits, then everything in this cell squares to zero.

For the count in (d), we consider symmetric diagrams, which is the same as counting half-diagrams.
Fix integers $k$ and $i$ and look at such diagrams with 
$k$ through strands, $i$ propagating nonidentity blocks (connected components), and hence $k-2i$ fixed propagating blocks. 
Since we allow $\ell$ possible handle decorations on each nonthrough strand, the number of such diagrams is
\begin{gather*}
\ell^{k-2i}
\sum_{b=k}^{n}
\begin{Bmatrix}n\\b\end{Bmatrix}
\binom{b}{k}
\binom{k}{2i}
(2i-1)!!.
\end{gather*}
This is justified as follows. First partition the top and bottom rows of a symmetric diagram identically into $b$ blocks, in 
$\begin{Bsmallmatrix}n\\b\end{Bsmallmatrix}$ ways. Then choose $k$ of these $b$ blocks to be propagating, and from those $k$ blocks choose $2i$ of them to correspond to transpositions, matching them up in $(2i-1)!!$ ways. The remaining $k-2i$ propagating blocks are fixed and can each be decorated with one of $\ell$ handles, contributing the factor $\ell^{k-2i}$. Finally, we sum over $b$ (and elsewhere over the admissible values of $i$, or equivalently of $k-2i$).

For (c), observe that the $\K$-linear span of dotted permutation diagrams form an algebra isomorphic to the dequantized Ariki--Koike algebra of the symmetric group, with defining relation \autoref{Eq:Evaluate2}. The claim then follows from the theory of sandwich cellular algebras and \cite[Theorem 3.7]{Mat}.
\end{proof}

\begin{Remark}
For all other versions of the diagram categories essentially the same works. Replace:
\begin{itemize}
\item The apexes of $pPa_{\mathbf{c}}(n)$ are the same as for $Pa_{\mathbf{c}}(n)$; the apexes of Motzkin and rook Brauer are also the same but (some $a_r\neq 0,r\leq n$) needs to be replaced by ($a_1\neq 0$); Temperley--Lieb and Brauer have $\{\text{0 or 1},\dots,n-2,n\}$, for $a_1\neq 0$, and $\{\text{0 or 1},\dots,n-2,n\}\setminus\{0\}$, for $a_1=0$; planar rook and rook have $\{0,\dots,n\}$, for $a_0\neq 0$, and $\{n\}$, for $a_0=0$.
\item For the planar diagram algebras there are $\ell$ finite dimensional simples per apex, and the sandwiched algebra is $(\K[x]/\eqref{Eq:Evaluate2})^{\otimes k}$. For the symmetric ones the sandwiched algebras stay the same as above.
\item In the dimension formula, replace ``merge diagrams'' with the respective analog, and the contribution of $\#\{\text{standard tableaux of shape $\lambda$}\}$ is removed for planar versions. 
\end{itemize}
More details below are given in \autoref{P:PlanCellModuleDims} and \autoref{P:SymmCellMod} below.
However, Green's relations (cells) can change, and we elaborate on this in \autoref{S:Twistings}, and a more detailed analysis of their representations is given in \autoref{S:Resultssemi} and \autoref{S:Results}.
\end{Remark}

Denote by $\mathbf{s}$ a choice of $p/q$ so that the associated diagram algebras are semisimple. This happens, for example, when all $a_i$ are transcendental (they do not satisfy a polynomial equation with integer coefficients) and \autoref{Eq:Evaluate2} has distinct roots. Note that under these assumptions, the cell modules $\Delta$ (subscript = number of through strands, superscript = cell, bracket = additional label from the sandwich if needed) are the simples. When the cell modules are not simple, we replace $\Delta$ with $L$ to distinguish the simple modules from the cell modules. Below $\K$ is of characteristic $e$, as before.

\begin{Proposition}\label{P:PlanCellModuleDims}
The representations of the planar diagram algebras have the following properties:
\begin{enumerate}
\item The simple representations in $pPa_\mathbf{s}(n)$ are indexed $1:1$ (up to isomorphism) by the number of through strands $k\in\{n,n-1,\dots,1,0\}$ and $\{1,\dots,\ell^k\}$, and, for $\ell=1$, have dimensions as for $TL_{\mathbf{s}'}(2n)$ below.

\item The simple representations in $Mo_\mathbf{s}(n)$ are indexed $1:1$ (up to isomorphism) by the number of through strands $k\in\{n,n-1,\dots,1,0\}$ and $\{1,\dots,\ell^k\}$, and, for $\ell=1$, have dimensions
\[
\dim\Delta^n_k  =\sum_{t=0}^{n}\frac{k+1}{k+t+1}\binom{n}{k+2t}
\binom{k+2t}{t}.
\]

\item The simple representations in $TL_\mathbf{s}(n)$ are indexed $1:1$ (up to isomorphism) by the number of through strands $k \in \{ n, n-2, \dots ,0\}$ if $n$ is even and $n \in \{ n, n-2, \dots ,1\}$ if $n$ is odd, and $\{1,\dots,\ell^k\}$, and, for $\ell=1$, have dimensions
\[
\dim\Delta^n_k =\tfrac{n-2c(k)+1}{n-c(k)+1}\binom{n}{c(k)},
\]
where $c(k) = \tfrac{n-k}{2}$.

\item The simple representations in $pRo_\mathbf{s}(n)$ are indexed $1:1$ (up to isomorphism) by the number of through strands $k \in \{ n, n-1, \dots 0\}$ and $\{1,\dots,\ell^k\}$, and, for $\ell=1$, have dimensions
\[
\dim\Delta^n_k =\binom{n}{k}.
\]
\end{enumerate}
\end{Proposition}

\begin{proof}
Our sandwiched algebra is $(\K[x]/\eqref{Eq:Evaluate2})^{\otimes k}$, and the Chinese remainder theorem implies that the number of its simple representations is $\ell^k$

For everything except $pPa_\mathbf{s}(n)$, the proposition then follows from \cite[Proposition 4B.3, Proposition 4F.4]{khovanov-monoidal-2024} and \autoref{T:Cob}. For $pPa_\mathbf{s}(n)$, we take $pPa_\delta(n)$ (see the isomorphism from \autoref{R:PlanarPartition}) as an isomorphism of sets. Then semisimplicity implies that the planar partition case is a ``doubling'' of Temperley--Lieb.
\end{proof}

For symmetric diagram algebras we have additional representations indexed by $e$-restricted $\ell$ multipartitions partitions $\lambda\vdash_\ell^e k$ and we write $\Delta^n_k(\lambda)$ for these.

\begin{Proposition}\label{P:SymmCellMod}
The representations of the symmetric diagram algebras have the following properties:
\begin{enumerate}
\item The simple representations in $Pa_\mathbf{s}(n)$ are indexed $1:1$ (up to isomorphism) by $\lambda \vdash_\ell^e k$ for each number of through strands $k \in \{ n, n-1, \dots, 1,0 \}$, and, for $\ell=1$, have dimensions
\[
\dim\Delta^n_k(\lambda) = \sum_{t=0}^{n}\begin{Bmatrix}n\\t\end{Bmatrix}\binom{t}{k}\cdot \#\{\text{standard tableaux of shape $\lambda$}\},
\]
where $\begin{Bsmallmatrix}n\\t\end{Bsmallmatrix}$ denotes the Stirling number of the second kind.

\item The simple representations in $RoBr_\mathbf{s}(n)$ are indexed $1:1$ (up to isomorphism) by $\lambda \vdash_\ell^e k$ for each number of through strands $k \in \{ n, n-1, \dots, 1,0 \}$, and, for $\ell=1$, have dimensions
\[
\dim\Delta^n_k(\lambda) = \sum_{t=0}^{n}\binom{n}{k}\binom{n-k}{2t}(2t-1)!!\cdot \#\{\text{standard tableaux of shape $\lambda$}\}.
\]

\item The simple representations in $Br_\mathbf{s}(n)$ are indexed $1:1$ (up to isomorphism) by $\lambda \vdash_\ell^e k$ for each number of through strands $k \in \{ n, n-2, \dots 0\}$ if $n$ is even and $n \in \{ n, n-2, \dots 1\}$ if $n$ is odd, and, for $\ell=1$, have dimensions
\[
\dim\Delta^n_k(\lambda) = \binom{n}{k}(n-k-1)!!\cdot \#\{\text{standard tableaux of shape $\lambda$}\}.
\]

\item The simple representations in $Ro_\mathbf{s}(n)$ are indexed $1:1$ (up to isomorphism) by $\lambda \vdash_\ell^e k$ for each number of through strands $k \in \{ n, n-1, \dots 0\}$, and, for $\ell=1$, have dimensions
\[
\dim\Delta^n_k(\lambda) = \binom{n}{k}\cdot \#\{\text{standard tableaux of shape $\lambda$}\}.
\]
\end{enumerate}
\end{Proposition}

\begin{proof}
This follows from \cite[Proposition 5B.2, Proposition 5F.2]{khovanov-monoidal-2024} and \autoref{T:Cob}. 
\end{proof}

\section{Generalized diagram monoids: twists}\label{S:Twistings}

This section defines and studies an even larger class of diagram monoids.
Everything we say below works for semigroups, but we often say monoid for convenience. 
Moreover, as before, we assume familiarity with \emph{Green's relations} \cite{Gr-structure-semigroups}. For the algebraically minded reader, these are the monoid version of \emph{cells}. The following summarizes our conventions, with notation in line with \cite{Tu-sandwich}:
\begin{gather*}
\begin{tikzpicture}[baseline=(A.center),every node/.style=
{anchor=base,minimum width=1.4cm,minimum height=1cm}]
\matrix (A) [matrix of math nodes,ampersand replacement=\&] 
{
\hcell_{11} \& \hcell_{12} 
\& \hcell_{13} \& \hcell_{14}
\\
\hcell_{21} \& \hcell_{22} 
\& \hcell_{23} \& \hcell_{24}
\\
\hcell_{31} \& \hcell_{32} 
\& \hcell_{33} \& \hcell_{34}
\\
};
\draw[fill=blue,opacity=0.25] (A-3-1.north west) node[blue,left,xshift=0.15cm,yshift=-0.5cm,opacity=1] 
{$\rcell$} to (A-3-4.north east) 
to (A-3-4.south east) to (A-3-1.south west) to (A-3-1.north west);
\draw[fill=red,opacity=0.25] (A-1-3.north west) node[red,above,xshift=0.7cm,opacity=1] 
{$\lcell$} to (A-3-3.south west) 
to (A-3-3.south east) to (A-1-3.north east) to (A-1-3.north west);
\draw[very thick,black,dotted] (A-1-2.north west) to (A-3-2.south west);
\draw[very thick,black,dotted] (A-1-3.north west) to (A-3-3.south west);
\draw[very thick,black,dotted] (A-1-4.north west) to (A-3-4.south west);
\draw[very thick,black,dotted] (A-2-1.north west) to (A-2-4.north east);
\draw[very thick,black,dotted] (A-3-1.north west) to (A-3-4.north east);
\draw[very thick,black] (A-1-1.north west) node[black,above,xshift=-0.5cm] {$\jcell$} to 
(A-1-4.north east) to (A-3-4.south east) 
to (A-3-1.south west) to (A-1-1.north west);
\draw[very thick,black,->] ($(A-1-1.north west)+(-0.4,0.4)$) to (A-1-1.north west);
\draw[very thick,black,->] ($(A-3-4.south east)+(0.5,0.2)$) 
node[right]{$\hcell(\lcell,\rcell)=\hcell_{33}$} 
to[out=180,in=0] ($(A-3-3.south east)+(0,0.2)$);
\end{tikzpicture}
.
\end{gather*}
This is often called an \emph{egg box diagram} or \emph{cell}. There are also the associated orders denoted, for example, $<_\mathcal{L}$.

\subsection{Representation theoretic motivation for what follows}

We recall the following way to compute the dimensions of simple representations of a monoid, or, somewhat analogously, of sandwich cellular algebras. We restrict to the case where the sandwich part is trivial. To get started, given a subset \( X \) of a semigroup \( S \), define the \emph{set of idempotents} \( E \left( X \right) \) to be
\begin{gather*}
E \left( X \right) = \left\{ e \in X \mid e^2 = e \right\} \,.
\end{gather*}

\begin{Definition}
Given a \( \mathcal{J} \)-class \( J \), define a matrix \( \operatorname{G} \left( J \right) \) as follows: let \( J \) have \( l \) many \( \mathcal{R} \)-classes and \( m \) many \( \mathcal{L} \)-classes, both of which are ordered. For all \( j \in \left\{ 1, 2, \dotsc, l \right\} \) and for all \( k \in \left\{ 1, 2, \dotsc, m \right\} \) define the \( \mathcal{H} \)-class \( H_{j, k} \) to be
\begin{gather*}
H_{jk} = R_j \cap L_k,
\end{gather*}
where \( R_j \) is the \( j \)th \( \mathcal{R} \)-class and \( L_k \) is the \( k \)th \( \mathcal{L} \)-class. Then \( \operatorname{G} \left( J \right) \) is an \( l \times m \)-matrix with entries
\begin{gather*}
\operatorname{G} \left( J \right)_{jk} =
\begin{cases}
1 & \text{if } E \left( H_{jk} \right) \neq \emptyset, \\
0 & \text{otherwise.}
\end{cases}
\end{gather*}
The matrix \( \operatorname{G} \left( J \right) \) is called the \emph{Gram matrix}. More generally, in the sandwich cellular case, we replace the entry \( 1 \) by the eigenvalue of a basis element that is a pseudo-idempotent (meaning: if \( e^2 = a \cdot e \) for some $a\in\K$, then we record \( a \)).
\end{Definition}

Denote by \( \textnormal{rk}\bigl(\operatorname{G} \left( J \right)\bigr) \) the rank of the Gram matrix. An \emph{apex} is a \( \mathcal{J} \)-class such that the Gram matrix is not zero.

\begin{Theorem}\label{T:Gram}
Let \( A \) be a finite monoid with trivial \( \mathcal{H} \)-classes, or any of the planar diagram algebras from \autoref{S:DiagCat} such that \autoref{Eq:Evaluate2} has only two summands.
Then every simple \( A \)-module is 1:1 associated to an apex \(J\), and this correspondence is 1:1. Moreover, the dimension of the module is
\( \textnormal{rk}\bigl(\operatorname{G} \left( J \right)\bigr) \).
\end{Theorem}

\begin{proof}
Well known; see e.g.\ \cite{GaMaSt-irreps-semigroups} or \cite{Tu-sandwich} for newer expositions.
\end{proof}

There is also a more general version of \autoref{T:Gram} with a more complicated matrix (often called the \emph{sandwich matrix}) that is ``essentially obtained from the Gram matrix by induction'', cf.\ \cite{St-rep-monoid,Tu-sandwich} for some additional details. In any case, \autoref{T:Gram} motivates us to look at idempotents in monoids and how they arrange themselves into cells.

\begin{Example}\label{E:TLcells1}
For Temperley--Lieb and four strands we get:
\begin{gather*}
\xy
(0,0)*{\begin{gathered}
\begin{tabular}{C|C}
\arrayrulecolor{tomato}
\cellcolor{mydarkblue!25}
\begin{tikzpicture}[anchorbase]
\draw[usual] (0,0) to[out=90,in=180] (0.25,0.2) to[out=0,in=90] (0.5,0);
\draw[usual] (0,0.5) to[out=270,in=180] (0.25,0.3) to[out=0,in=270] (0.5,0.5);
\draw[usual] (1,0) to[out=90,in=180] (1.25,0.2) to[out=0,in=90] (1.5,0);
\draw[usual] (1,0.5) to[out=270,in=180] (1.25,0.3) to[out=0,in=270] (1.5,0.5);
\end{tikzpicture} &
\cellcolor{mydarkblue!25}
\begin{tikzpicture}[anchorbase]
\draw[usual] (0,0) to[out=45,in=180] (0.75,0.20) to[out=0,in=135] (1.5,0);
\draw[usual] (0.5,0) to[out=90,in=180] (0.75,0.1) to[out=0,in=90] (1,0);
\draw[usual] (0,0.5) to[out=270,in=180] (0.25,0.3) to[out=0,in=270] (0.5,0.5);
\draw[usual] (1,0.5) to[out=270,in=180] (1.25,0.3) to[out=0,in=270] (1.5,0.5);
\end{tikzpicture}
\\
\hline
\cellcolor{mydarkblue!25}
\begin{tikzpicture}[anchorbase]
\draw[usual] (0,0) to[out=90,in=180] (0.25,0.2) to[out=0,in=90] (0.5,0);
\draw[usual] (1,0) to[out=90,in=180] (1.25,0.2) to[out=0,in=90] (1.5,0);
\draw[usual] (0,0.5) to[out=315,in=180] (0.75,0.3) to[out=0,in=225] (1.5,0.5);
\draw[usual] (0.5,0.5) to[out=270,in=180] (0.75,0.4) to[out=0,in=270] (1,0.5);
\end{tikzpicture} &
\cellcolor{mydarkblue!25}
\begin{tikzpicture}[anchorbase]
\draw[usual] (0,0) to[out=45,in=180] (0.75,0.20) to[out=0,in=135] (1.5,0);
\draw[usual] (0,0.5) to[out=315,in=180] (0.75,0.3) to[out=0,in=225] (1.5,0.5);
\draw[usual] (0.5,0) to[out=90,in=180] (0.75,0.1) to[out=0,in=90] (1,0);
\draw[usual] (0.5,0.5) to[out=270,in=180] (0.75,0.4) to[out=0,in=270] (1,0.5);
\end{tikzpicture}
\end{tabular}
\\[3pt]
\begin{tabular}{C|C|C}
\arrayrulecolor{tomato}
\cellcolor{mydarkblue!25}
\begin{tikzpicture}[anchorbase]
\draw[usual] (0,0) to[out=90,in=180] (0.25,0.2) to[out=0,in=90] (0.5,0);
\draw[usual] (0,0.5) to[out=270,in=180] (0.25,0.3) to[out=0,in=270] (0.5,0.5);
\draw[usual] (1,0) to (1,0.5);
\draw[usual] (1.5,0) to (1.5,0.5);
\end{tikzpicture} & 
\cellcolor{mydarkblue!25}
\begin{tikzpicture}[anchorbase]
\draw[usual] (0.5,0) to[out=90,in=180] (0.75,0.2) to[out=0,in=90] (1,0);
\draw[usual] (0,0.5) to[out=270,in=180] (0.25,0.3) to[out=0,in=270] (0.5,0.5);
\draw[usual] (0,0) to (1,0.5);
\draw[usual] (1.5,0) to (1.5,0.5);
\end{tikzpicture} &
\begin{tikzpicture}[anchorbase]
\draw[usual] (1,0) to[out=90,in=180] (1.25,0.2) to[out=0,in=90] (1.5,0);
\draw[usual] (0,0.5) to[out=270,in=180] (0.25,0.3) to[out=0,in=270] (0.5,0.5);
\draw[usual] (0,0) to (1,0.5);
\draw[usual] (0.5,0) to (1.5,0.5);
\end{tikzpicture}
\\
\hline
\cellcolor{mydarkblue!25}
\begin{tikzpicture}[anchorbase]
\draw[usual] (0,0) to[out=90,in=180] (0.25,0.2) to[out=0,in=90] (0.5,0);
\draw[usual] (0.5,0.5) to[out=270,in=180] (0.75,0.3) to[out=0,in=270] (1,0.5);
\draw[usual] (1,0) to (0,0.5);
\draw[usual] (1.5,0) to (1.5,0.5);
\end{tikzpicture} & 
\cellcolor{mydarkblue!25}
\begin{tikzpicture}[anchorbase]
\draw[usual] (0,0) to (0,0.5);
\draw[usual] (0.5,0) to[out=90,in=180] (0.75,0.2) to[out=0,in=90] (1,0);
\draw[usual] (0.5,0.5) to[out=270,in=180] (0.75,0.3) to[out=0,in=270] (1,0.5);
\draw[usual] (1.5,0) to (1.5,0.5);
\end{tikzpicture} &
\cellcolor{mydarkblue!25}
\begin{tikzpicture}[anchorbase]
\draw[usual] (0,0) to (0,0.5);
\draw[usual] (0.5,0) to (1.5,0.5);
\draw[usual] (1,0) to[out=90,in=180] (1.25,0.2) to[out=0,in=90] (1.5,0);
\draw[usual] (0.5,0.5) to[out=270,in=180] (0.75,0.3) to[out=0,in=270] (1,0.5);
\end{tikzpicture}
\\
\hline
\begin{tikzpicture}[anchorbase]
\draw[usual] (0,0) to[out=90,in=180] (0.25,0.2) to[out=0,in=90] (0.5,0);
\draw[usual] (1,0.5) to[out=270,in=180] (1.25,0.3) to[out=0,in=270] (1.5,0.5);
\draw[usual] (1,0) to (0,0.5);
\draw[usual] (1.5,0) to (0.5,0.5);
\end{tikzpicture} & 
\cellcolor{mydarkblue!25}
\begin{tikzpicture}[anchorbase]
\draw[usual] (0,0) to (0,0.5);
\draw[usual] (1.5,0) to (0.5,0.5);
\draw[usual] (0.5,0) to[out=90,in=180] (0.75,0.2) to[out=0,in=90] (1,0);
\draw[usual] (1,0.5) to[out=270,in=180] (1.25,0.3) to[out=0,in=270] (1.5,0.5);
\end{tikzpicture} &
\cellcolor{mydarkblue!25}
\begin{tikzpicture}[anchorbase]
\draw[usual] (0,0) to (0,0.5);
\draw[usual] (0.5,0) to (0.5,0.5);
\draw[usual] (1,0) to[out=90,in=180] (1.25,0.2) to[out=0,in=90] (1.5,0);
\draw[usual] (1,0.5) to[out=270,in=180] (1.25,0.3) to[out=0,in=270] (1.5,0.5);
\end{tikzpicture}
\end{tabular}
\\[3pt]
\begin{tabular}{C}
\arrayrulecolor{tomato}
\cellcolor{mydarkblue!25}
\begin{tikzpicture}[anchorbase]
\draw[usual] (0,0) to (0,0.5);
\draw[usual] (0.5,0) to (0.5,0.5);
\draw[usual] (1,0) to (1,0.5);
\draw[usual] (1.5,0) to (1.5,0.5);
\end{tikzpicture}
\end{tabular}
\end{gathered}};
(-45,13)*{\jcell_{0}};
(-45,-2.5)*{\jcell_{2}};
(-45,-16.5)*{\jcell_{4}};
(45,13)*{\sand[0]\cong\onemon};
(45,-2.5)*{\sand[2]\cong\onemon};
(45,-16.5)*{\sand[4]\cong\onemon};
(-51,0)*{\phantom{a}};
\endxy
\end{gather*}
as the cell picture, with pseudo-idempotents in colored boxes. The two main monoid examples we have are $a_1=1$ and $a_1=0$, and generically we could take $a_1=\pi$, for which
\begin{gather*}
a_1=1\colon
\begin{gathered}
\begin{pmatrix}
1 & 1\\ 1 & 1
\end{pmatrix}
\\
\begin{pmatrix}
1 & 1 & 0 \\ 1& 1 & 1\\ 0& 1 & 1
\end{pmatrix}
\\
\begin{pmatrix}
1
\end{pmatrix}
\end{gathered}
,\quad
a_1=0\colon
\begin{gathered}
\begin{pmatrix}
0
\end{pmatrix}
\\
\begin{pmatrix}
0 & 0\\ 0 & 0
\end{pmatrix}
\\
\begin{pmatrix}
0 & 1 & 0 \\ 1& 0 & 1\\ 0& 1 & 0
\end{pmatrix}
\\
\begin{pmatrix}
1
\end{pmatrix}
\end{gathered}
,\quad
a_1=\pi\colon
\begin{gathered}
\begin{pmatrix}
\pi^2 & \pi\\ \pi & \pi^2
\end{pmatrix}
\\
\begin{pmatrix}
\pi & 1 & 0 \\ 1& \pi & 1\\ 0& 1 & \pi
\end{pmatrix}
\\
\begin{pmatrix}
1
\end{pmatrix}
\end{gathered}
,
\end{gather*}
are then the Gram matrices. (Note that for $a_1=0$ we attach a zero element which corresponds to the top matrix. It does nothing important and can be ignored.) The associated simple dimensions are therefore (bottom to top) $1,3,1$, and $1,2$, and $1,3,2$.
\end{Example}

\subsection{Equivalences}

Below we allow infinite monoids and semigroups for which \( \mathcal{D} \)- and \( \mathcal{J} \)-classes may differ. For finite monoids and semigroups they are the same and are also called two-sided cells.

\begin{Definition} \label{def: D class preserving}
Define the binary relation \( \simeq \) on the collection of all \( \mathcal{D} \)-classes of all semigroups as follows. Two \( \mathcal{D} \)-classes \( D_1 \) and \( D_2 \) are related if there exists a function \( \Theta : D_1 \to D_2 \) such that all of the following conditions hold:
\begin{enumerate}
\item\label{def: D class preserving i} Bijection: \( \Theta \) is a bijection;

\item\label{def: D class preserving ii} \( \mathcal{L} \)-class preserving: for all \( a, b \in D_1 \) it holds that \( a \mathrel{\mathcal{L}} b \iff \Theta \left( a \right) \mathrel{\mathcal{L}} \Theta \left( b \right) \);

\item\label{def: D class preserving iii} \( \mathcal{R} \)-class preserving: for all \( a, b \in D_1 \) it holds that \( a \mathrel{\mathcal{R}} b \iff \Theta \left( a \right) \mathrel{\mathcal{R}} \Theta \left( b \right) \).\qedhere
\end{enumerate}
(For readers familiar with cells: for example, \( a \mathrel{\mathcal{L}} b \) is often written \( a \sim_L b \), or similarly.)
\end{Definition}

\begin{Lemma} \label{lem: D class preserving}
The binary relation \( \simeq \) is an equivalence relation.
\end{Lemma}

\begin{proof}
Boring and omitted.
\end{proof}

Using this, two \( \mathcal{D} \)-classes \( D_1 \) and \( D_2 \) are \emph{\( \mathcal{D} \)-class preserving} if \( D_1 \simeq D_2 \). Note that, and that is key, \( D_1 \) and \( D_2 \) are allowed to be from different semigroups. The notion is however not quite what we need, as we want to keep track of idempotents.

\begin{Definition} \label{def: D class idempotent preserving}
Define the binary relation \( \cong \) on the collection of all \( \mathcal{D} \)-classes of all semigroups as follows. Two \( \mathcal{D} \)-classes \( D_1 \) and \( D_2 \) are related if all of the following conditions hold:
\begin{enumerate}
\item\label{def: D class idempotent preserving i} \( \mathcal{D} \)-class preserving: \( D_1 \simeq D_2 \) with function \( \Theta : D_1 \to D_2 \) from \autoref{def: D class preserving};

\item\label{def: D class idempotent preserving ii} Idempotent preserving: for all \( \mathcal{H} \)-classes \( H \subseteq D_1 \) it holds that
\[
E \left( H \right) = \emptyset \iff  E \big( \Theta \left( H \right) \big) = \emptyset \,.\qedhere
\]
\end{enumerate}
With \autoref{lem: D class idempotent preserving} below, \( D_1 \) and \( D_2 \) are \emph{\( \mathcal{D} \)-class idempotent preserving} if \( D_1 \cong D_2 \).
\end{Definition}

It is clear from~\cref{def: D class preserving ii,def: D class preserving iii} of \autoref{def: D class preserving} that \( \mathcal{H} \)-classes are preserved by \( \Theta \), and hence it makes sense to talk about idempotents being preserved in~\cref{def: D class idempotent preserving ii} of the above definition.

\begin{Lemma} \label{lem: D class idempotent preserving}
The binary relation \( \cong \) is an equivalence relation.
\end{Lemma}

\begin{proof}
By \autoref{lem: D class preserving} the relation \( \simeq \) is an equivalence relation, and that~\cref{def: D class idempotent preserving ii} also defines an equivalence relation again follows straightforwardly from the definitions.
\end{proof}

The motivation behind \autoref{def: D class idempotent preserving} is that the egg box diagrams of the two \( \mathcal{D} \)-classes, where the idempotents are indicated in some way, are identical and hence the ranks of their corresponding Gram matrices are the same. Formally:

\begin{Lemma} \label{lem: implies same rank}
Given two \( \mathcal{D} \)-classes \( D_1 \) and \( D_2 \), it holds that
\begin{gather*}
D_1 \cong D_2 \implies  \textnormal{rk}\bigl( \operatorname{G} \left( D_1 \right) \bigr) = \textnormal{rk} \bigl( \operatorname{G} \left( D_2 \right) \bigr) \,.
\end{gather*}
\end{Lemma}

\begin{proof}
As \( D_1 \cong D_2 \), it is clear that the number of \( \mathcal{H} \)-classes that contain an idempotent is the same between \( D_1 \) and \( D_2 \). Therefore it remains to be shown that the change in arrangement of these idempotents from \( D_1 \) to \( D_2 \) is a permutation of rows and columns, and hence the ranks of the corresponding Gram matrices are the same by elementary linear algebra.

We prove this by contradiction. Suppose, WLOG (in the sense that the dual argument can be used if there exist two such \( \mathcal{L} \)-classes and appropriately dualized hypotheses), there exist two \( \mathcal{R} \)-classes \( R_1, R_2 \subseteq D_1 \) that have different permutations of their columns, or \( \mathcal{H} \)-classes, when mapped under \( \Theta \). Then, WLOG (up to relabeling of the indices), suppose \( H_{1, 1} \subseteq R_1 \) and \( H_{2, 1} \subseteq R_2 \) are such that
\begin{gather}
H_{1, 1}, H_{2, 1} \subseteq L_1 \,, \label{eq: implies same rank 1}
\end{gather}
for some \( \mathcal{L} \)-class \( L_1 \subseteq D_1 \), but \( \Theta \left( H_{1, 1} \right) \) and \( \Theta \left( H_{2, 1} \right) \) belong to different \( \mathcal{L} \)-classes in \( D_2 \).

For example, consider the following relevant parts of \( D_1 \) being mapped to \( D_2 \) by \( \Theta \):
\begin{gather*}
\begin{tblr}{c|c|c|c|c|}
\cline{2-5}
R_1 & H_{1, 1} & H_{1, 2} & H_{1, 3} & H_{1, 4} \\
\cline{2-5}
R_2 & H_{2, 1} & H_{2, 2} & H_{2, 3} & H_{2, 4} \\
\cline{2-5}
\end{tblr}%
\quad%
\to%
\begin{tblr}{c|c|c|c|c|}
\cline{2-5}
\Theta \left( R_1 \right) & \Theta \left( H_{1, 1} \right) & \Theta \left( H_{1, 2} \right) & \Theta \left( H_{1, 3} \right) & \Theta \left( H_{1, 4} \right) \\
\cline{2-5}
\Theta \left( R_2 \right) & \Theta \left( H_{2, 2} \right) & \Theta \left( H_{2, 1} \right) & \Theta \left( H_{2, 3} \right) & \Theta \left( H_{2, 4} \right) \\
\cline{2-5}
\end{tblr} \,.
\end{gather*}
Here it can be seen that the permutation \( \Theta \) applied to \( R_1 \) is the identity, whilst for \( R_2 \) it swaps the first and second \( \mathcal{H} \)-classes.

Let \( a \in H_{1, 1} \) and let \( b \in H_{2, 1} \). By~\cref{eq: implies same rank 1} it follows that \( a \mathrel{\mathcal{L}} b \), and hence by~\cref{def: D class preserving ii} of \autoref{def: D class preserving} it follows that \( \Theta \left( a \right) \mathrel{\mathcal{L}} \Theta \left( b \right) \). But this contradicts that \( \Theta \left( H_{1, 1} \right) \) and \( \Theta \left( H_{2, 1} \right) \) belong to different \( \mathcal{L} \)-classes in \( D_2 \). Hence all \( \mathcal{R} \)-classes of \( D_1 \) have the same permutation of their \( \mathcal{H} \)-classes when mapped under \( \Theta \), and by the dual argument, all \( \mathcal{L} \)-classes of \( D_1 \) have the same permutation of their \( \mathcal{H} \)-classes when mapped under \( \Theta \), as required.
\end{proof}

All of the above, by keeping track of the eigenvalue, applies to the case of having pseudo-idempotents as well.

\begin{Example}\label{E:TLcells2}
In \autoref{E:TLcells1}, the bottom $\mathcal{D}$ class is $\cong$ between all of the examples given, and it is the only $\mathcal{D}$ class with that property.
\end{Example}

\subsection{Twistings}

Recall twistings as in \cite{east2025twisted}.

\begin{Definition}
Given a semigroup \(S\), define a \emph{twisting} to be a map \(\Phi : S \times S \to \mathbb{N}\) such that for all \(a,b,c \in S\) it holds that
\begin{gather}
\Phi(a,b) + \Phi(ab,c) = \Phi(a,bc) + \Phi(b,c) \,.
\end{gather}
With a slight abuse of notation, \(\Phi(a,b,c)\) denotes this common value.
\end{Definition}

\begin{Definition}\label{def: tight}
Given a semigroup \(S\) and associated twisting \(\Phi\), the twisting is \emph{tight} if all of the following conditions hold:
\begin{enumerate}
\item for all \(a,b \in S\) there exists \(a' \in S\) such that \(ab = a'b\) and \(\Phi(a',b) = 0\);

\item for all \(a,b \in S\) there exists \(b' \in S\) such that \(ab = ab'\) and \(\Phi(a,b') = 0\).
\end{enumerate}
Otherwise, \(\Phi\) is \emph{loose}.
\end{Definition}

The following proposition describes a canonical twisting for diagram monoids.

\begin{Proposition}\label{P:CanonicalTwisting}
Let \(\mathbf{a} = \{a_i \mid i \in \mathbb{N}\}\) be a choice of monoid parameters, and let \(\monoid_1(n)\) be any of the classical diagram monoids in \autoref{S:DiagCat}. Let \(x,y \in \monoid_{1}(n)\), and classify each floating component of the product \(xy \in \monoid_{\mathbf{a}}(n)\) by its corresponding closed surface \(a_i^{x,y}\) as in \autoref{Eq:Evaluate}, where in particular each \(a_i^{x,y}\) has a corresponding parameter \(a_i \in \mathbf{a}\). 
Define the twisting
\[
\Phi_{\mathbf{a}} : \monoid_{1}(n)\times \monoid_{1}(n) \longrightarrow \mathbb{N}
\]
by
\[
\Phi_{\mathbf{a}}(x,y) = 
\begin{cases}
0 & \textnormal{if } a_i = 0 \textnormal{ for some } a_i^{x,y}, \\[0.2em]
\#\{\textnormal{floating components}\} & \textnormal{otherwise.}
\end{cases}
\]
Then we have the following facts:
\begin{enumerate}
\item For \(\monoid \in \{Pa,Br,pPa,TL\}\), the twisting \(\Phi_{\mathbf{a}}\) is tight for all applicable \(\mathbf{a}\).
\item For \(\monoid \in \{RoBr,Mo\}\), the twisting \(\Phi_{\mathbf{a}}\) is tight for \(a_0 = 1\) and loose for \(a_0 = 0\). 
\item For \(\monoid \in \{ Ro,pRo\}\), the twisting \(\Phi_{\mathbf{a}}\) is tight for \(a_0 = 1\) and loose for \(a_0 = 0\).
\end{enumerate}
\end{Proposition}

\begin{proof}
We will prove this in a follow-up paper. The curious reader might try this themselves, as it is not difficult (but involved to write down formally, so postponed by us).
\end{proof}

\begin{Definition} \label{def: twisted S}
Given a semigroup \(S\), an associated twisting \(\Phi\), a commutative monoid \(M\) (additively written), and a (distinguished) fixed element \(q \in M\), define the \emph{twisted (product) semigroup}
\[
T = M \times_\Phi^q S = (M \times S,\cdot)
\]
where
\begin{gather*}
(j,a)(k,b) = (j + k + \Phi(a,b)q,\,ab)
\end{gather*}
and \(\Phi(a,b)q\) means \(q + q + \dotsb + q\) with \(\Phi(a,b)\) summands.
\end{Definition}

\begin{Lemma}
The twisted product semigroup is a semigroup (i.e.\ the composition is associative).
\end{Lemma}

\begin{proof}
A key result of \cite{east2025twisted}.
\end{proof}

For such $T$ we will call $S$ the \emph{untwisted} semigroup.

\begin{Remark}
For a diagram monoid \(S\) and a twisted product monoid \(T = M \times_{\Phi}^q S\), we can define a twisted diagram category \(\mathbf{T}\) with a sandwich cellular structure in the same way as in \autoref{S:DiagCat}. The data is essentially equivalent, so we stay with the monoid version for convenience.
\end{Remark}

\begin{Example}\label{E}
Twistings give generally new examples of diagram monoids and categories that we have not seen in the tensor category and related literature. A nice example is the so-called \emph{Sylvester twisting} that works for all the diagram monoids of \autoref{S:DiagCat}. (To make this work one copies \cite[Section 4.4]{east2025twisted}.)
\end{Example}

Twisting may reduce the number of idempotents:

\begin{Lemma} \label{lem: E(T) <= E(S)}
Given a twisted product semigroup \(T = M \times_\Phi^q S\), it holds that
\begin{gather*}
p_2\big(E(T)\big) \subseteq E(S)
\end{gather*}
where \(p_2\) is the projection onto the second coordinate.
\end{Lemma}

\begin{proof}
Immediate from \autoref{def: twisted S}.
\end{proof}

\begin{Definition}
A twisted semigroup \(T = M \times_\Phi^q S\) is tight if \(\Phi\) is tight, and loose otherwise.
\end{Definition}

\begin{Definition} \label{def: 0 twisted}
Given a semigroup \(S = (X,\cdot)\) and associated twisting \(\Phi\), define the \emph{\(\mathbf{0}\)-twisted semigroup}
\[
T_\Phi^0 = (X \sqcup \{0\}, *)
\]
where
\begin{gather*}
a * b =
\begin{cases}
ab & \text{if } a,b \neq 0 \text{ and } \Phi(a,b) = 0, \\
0 & \text{otherwise}.
\end{cases}
\end{gather*}
(This can be thought of as ``setting parameters to zero'', but is more general.)
\end{Definition}

\begin{Example}
For a diagram monoid \(\monoid_{\mathbf{a}}(n)\) (e.g. the ones from \autoref{S:DiagCat}) such that \(\mathbf{a}\neq \mathbf{1}\) (see \autoref{N:ClassicalDiagramMonoids}), let \(D_1,D_2,\dots,D_k\) be the \(\mathcal{D}\)-classes of \(\monoid_{\mathbf{a}}(n)\). Then the \(0\)-twisted monoid associated with the counting floating components twisting \(\Phi_{\mathbf{a}}\) from \autoref{P:CanonicalTwisting} has \(\mathcal{D}\)-classes \(D_1^0,D_2^0,\dots,D_k^0\) such that \(D_i \cong D_i^0\) (defined in \autoref{def: D class idempotent preserving}) for all \(1\leq i \leq k\).
\end{Example}

\begin{Lemma} \label{lem: E(T0) <= E(S)}
Given a semigroup \(S\) and an associated \(0\)-twisted semigroup \(T_\Phi^0\), it holds that
\begin{gather*}
E\big(T_\Phi^0 \setminus \{0\}\big) \subseteq E(S) \,.
\end{gather*}
\end{Lemma}

\begin{proof}
Immediate from \autoref{def: 0 twisted}.
\end{proof}

\subsection{Green's relations}

We now study the cells of a twisted semigroup in relation to the original cells of the parent untwisted semigroup.

\begin{Lemma} \label{lem: 0 <=K x}
Given a semigroup \(S\) with zero, for all \(x \in S\) and for all \(\mathcal{K} \in \{\mathcal{R},\mathcal{L},\mathcal{J}\}\) it holds that \(0 \leq_{\mathcal{K}} x\).
\end{Lemma}

\begin{proof}
Boring.
\end{proof}

\begin{Lemma} \label{lem: x <=K 0 iff x=0}
Given a semigroup \(S\) with zero, for all \(x \in S\) and for all \(\mathcal{K} \in \{\mathcal{R},\mathcal{L},\mathcal{J}\}\) it holds that
\begin{gather*}
x \leq_{\mathcal{K}} 0 \iff x = 0 \,.
\end{gather*}
\end{Lemma}

\begin{proof}
Trivial.
\end{proof}

\begin{Proposition} \label{lem: twisted greens}
Given a semigroup \(S = (X,\cdot)\) and associated tight twisting \(\Phi\), and given \(a,b \in X\), for all
\[
\xi \in \{\mathcal{R},\mathcal{L},\mathcal{D},\mathcal{J},\mathcal{H},\leq_{\mathcal{R}},\leq_{\mathcal{L}},\leq_{\mathcal{J}}\}
\]
it holds that
\begin{gather*}
a \mathrel{\xi^{T_\Phi^0}} b \iff a \mathrel{\xi^S} b \,.
\end{gather*}
In other words, tight twistings preserve cells.
\end{Proposition}

\begin{proof}
Consider the following cases.

\begin{mcases}

\case  \label{cs: twisted greens 7}

\(\xi = \:\leq_{\mathcal{J}}\).

Suppose \(a \leq_{\mathcal{J}}^{T_\Phi^0} b\). The proof is trivial if \(a = b\), so suppose \(a \neq b\). Then there exist \(c,d \in X \sqcup \{0\}\) such that
\begin{gather}
a = c * b * d \,. \label{eq: twisted greens 1a}
\end{gather}
By \autoref{def: 0 twisted}, as \(a \neq 0\), it follows that \(c,d \neq 0\) and
\begin{subequations}
\begin{align*}
c b d &= c * b * d \\
&= a \by{by~\cref{eq: twisted greens 1a}} \,.
\end{align*}
\end{subequations}
Hence \(a \leq_{\mathcal{J}}^S b\).

Conversely, suppose \(a \leq_{\mathcal{J}}^S b\). The proof is trivial if \(a = b\), so suppose \(a \neq b\). Then there exist \(c,d \in X\) such that
\begin{gather}
a = c b d \,. \label{eq: twisted greens 3a}
\end{gather}
By tightness applied to \(b,d\) there exists \(d' \in X\) such that
\begin{gather}
b d = b * d' \,,  \label{eq: twisted greens 2a}
\end{gather}
and by tightness applied to \(c, b * d'\) there exists \(c' \in X\) such that \(c (b * d') = c' * (b * d')\). Hence
\begin{subequations}
\begin{align*}
c' * b * d' &= c (b * d') \\
&= c b d \by{by~\cref{eq: twisted greens 2a}} \\
&= a \by{by~\cref{eq: twisted greens 3a}} \,.
\end{align*}
\end{subequations}
Thus \(a \leq_{\mathcal{J}}^{T_\Phi^0} b\).

\case

\(\xi = \mathcal{J}\): immediate from \cref{cs: twisted greens 7}.

\case \label{cs: twisted greens 1}

\(\xi = \:\leq_{\mathcal{R}}\): analogous to \cref{cs: twisted greens 7}.

\case \label{cs: twisted greens 2}

\(\xi = \mathcal{R}\): immediate from \cref{cs: twisted greens 1}.

\case \label{cs: twisted greens 3}

\(\xi \in \{\mathcal{L},\leq_{\mathcal{L}}\}\): dual to \cref{cs: twisted greens 1,cs: twisted greens 2}.

\case \label{cs: twisted greens 4}

\(\xi \in \{\mathcal{D},\mathcal{H}\}\): immediate from \cref{cs: twisted greens 2,cs: twisted greens 3}.
\end{mcases}

\noindent The proof is complete.
\end{proof}

\begin{Lemma} \label{lem: twisted greens equivalences}
Given a semigroup \(S\) and associated tight twisting \(\Phi\), for all \(\mathcal{K} \in \{\mathcal{R},\mathcal{L},\mathcal{D},\mathcal{J},\mathcal{H}\}\) it holds that
\begin{gather*}
\mathcal{K}^{T_\Phi^0} = \mathcal{K}^S \cup \{(0,0)\} \,.
\end{gather*}
\end{Lemma}

\begin{proof}
Immediate from \autoref{lem: 0 <=K x} and \autoref{lem: x <=K 0 iff x=0}.
\end{proof}

\begin{Lemma}\label{thm: tight => Green's simple}
Given a twisted product semigroup \(T = M \times_\Phi^q S\) and given \((j,a),(k,b) \in T\), for all
\[
\xi \in \{\mathcal{R},\mathcal{L},\mathcal{D},\mathcal{J},\mathcal{H},\leq_{\mathcal{R}},\leq_{\mathcal{L}},\leq_{\mathcal{J}}\}
\]
if \(T\) is tight then
\begin{gather*}
(j,a) \mathrel{\xi^T} (k,b) \iff (j \mathrel{\xi^M} k \text{ and } a \mathrel{\xi^S} b) \,.
\end{gather*}
\end{Lemma}

\begin{proof}
This is \cite[Theorem~5.1]{east2025twisted}.
\end{proof}

\begin{Lemma} \label{lem: tight => Green's classes}
Given a tight twisted product semigroup \(T = M \times_\Phi^q S\), and given \((j,a) \in T\), for all \(\mathcal{K} \in \{\mathcal{R},\mathcal{L},\mathcal{D},\mathcal{J},\mathcal{H}\}\) it holds that
\begin{gather*}
K_{(j,a)}^T = K_j^M \times K_a^S = H_j^M \times K_a^S \,.
\end{gather*}
\end{Lemma}

\begin{proof}
This is \cite[Corollary~5.2(i)]{east2025twisted}.
\end{proof}

\subsection{Main theorem}

We now come back to idempotents.

\begin{Lemma} \label{lem: E(T)}
Given a tight twisted product semigroup \(T = M \times_\Phi^q S\) where \(M\) is \(\mathcal{D}\)-trivial, it holds that
\begin{gather*}
E(T) = \{(i,e) \in E(M) \times E(S) \mid i = i + \Phi(e,e)q\} \,.
\end{gather*}
\end{Lemma}

\begin{proof}
In the notation of \cite{east2025twisted} as used in their proof (see the proof of \cite[Proposition~6.4]{east2025twisted}), as \(M\) is \(\mathcal{D}\)-trivial (meaning all \(\mathcal{D}\)-classes are trivial) and commutative, we have that all Green's relations are trivial. It follows that \(\varepsilon(i,e) = (i,e)\) where \(i = i + \Phi(e,e)q\). This implies that \(i \leq_{\mathcal{J}}^M \Phi(e,e)q\) and hence \(\Omega = E(M) \times E(S)\).
\end{proof}

\begin{Lemma} \label{lem: xy = xyy}
Given a commutative additive semigroup \(M\), given \(i \in E(M)\) and given \(j \in M\), it holds that
\begin{gather*}
i + j = i + j + j \iff i + j \in E(M) \,.
\end{gather*}
\end{Lemma}

\begin{proof}
As \(i \in E(M)\) it holds that
\begin{subequations}
\begin{align*}
i + j = i + j + j &\iff i + j = i + i + j + j \\
&\iff i + j = i + j + i + j \by{commutativity} \\
&\iff i + j \in E(M) \,.
\end{align*}
\end{subequations}
The proof is complete.
\end{proof}

\begin{Lemma} \label{lem: xy^n = xy^m}
Given an additive semigroup \(M\) and \(j,k \in M\), for all \(n,m \in \mathbb{Z}_{>0}\) if \(j + k = j + k + k\) then
\begin{gather}
j + nk = j + mk \,.
\end{gather}
\end{Lemma}

\begin{proof}
Trivial.
\end{proof}

In order to state \autoref{thm: D classes not twisted or 0-twist}, we need the following condition. This is stated in its most general formulation. A more tractable but stronger hypothesis that implies this condition is given by \autoref{lem: exists implies for all idempotents}.

\begin{Condition} \label{con: exists implies for all idempotents}
Given a tight twisted product semigroup \(T = M \times_\Phi^q S\), and given a \(\mathcal{D}\)-class \(D_{(j,a)}^T\) of~\(T\), define the following condition:
if there exists \(e \in E(D_a^S)\) such that \(\Phi(e,e) > 0\) and \(j = j + \Phi(e,e)q\) then for all \(f \in E(D_a^S)\) it holds that
\begin{gather*}
\Phi(f,f) > 0 \implies j = j + \Phi(f,f)q \,.
\end{gather*}
\end{Condition}

\begin{Lemma} \label{lem: exists implies for all idempotents}
\autoref{con: exists implies for all idempotents} holds if \(E(D_{(j,a)}^T) \neq \emptyset\) and \(j + q \in E(M)\).
\end{Lemma}

\begin{proof}
Let \(f \in E(D_a^S)\) such that \(\Phi(f,f) > 0\). Trivially \(E(D_{(j,a)}^T) \neq \emptyset\) and \autoref{lem: E(T)} imply that \(j \in E(M)\). Hence, together with \(j + q \in E(M)\), the hypothesis of \autoref{lem: xy = xyy} is satisfied. This in turn satisfies the hypothesis of \autoref{lem: xy^n = xy^m} as \(\Phi(e,e),\Phi(f,f) > 0\). Hence
\begin{subequations}
\begin{align*}
j + \Phi(f,f)q &= j + \Phi(e,e)q \\
&= j \by{by assumption} \,.
\end{align*}
\end{subequations}
We are done.
\end{proof}

\begin{Theorem} \label{thm: D classes not twisted or 0-twist}
Given a tight twisted product semigroup \(T = M \times_\Phi^q S\) where \(M\) is \(\mathcal{D}\)-trivial, for all \(\mathcal{D}\)-classes \(D_{(j,a)}^T\) of \(T\), if \(E(D_{(j,a)}^T) \neq \emptyset\) and \autoref{con: exists implies for all idempotents} holds then one (possibly both) of the following holds:
\begin{enumerate}
\item \(D_{(j,a)}^T \cong D_a^S\), or

\item \(D_{(j,a)}^T \cong D_a^{T_\Phi^0}\).
\end{enumerate}
Precisely one (not both) of the above holds if and only if there exists \(e \in E(D_a^S)\) such that \(\Phi(e,e) > 0\).
\end{Theorem}

\begin{proof}
By \autoref{lem: tight => Green's classes}, it holds that \(D_{(j,a)}^T = D_j^M \times D_a^S\). As \(M\) is \(\mathcal{D}\)-trivial, it follows that \(D_j^M = \{j\}\) and hence \(D_{(j,a)}^T = \{j\} \times D_a^S\). Define \(\Theta : D_{(j,a)}^T \to D_a^S\) by
\begin{gather}
\Theta(j,u) = u \,. \label{eq: D classes not twisted or 0-twist 1}
\end{gather}
It is clear that \(\Theta\) is a bijection.

Let \((j,u),(j,v) \in D_{(j,a)}^T\). Then by \autoref{thm: tight => Green's simple}
\begin{subequations}
\begin{align*}
(j,u) \mathrel{\mathcal{L}^T} (j,v) 
&\iff (j \mathrel{\mathcal{L}^M} j \text{ and } u \mathrel{\mathcal{L}^S} v) \\
&\iff u \mathrel{\mathcal{L}^S} v \\
&\iff \Theta(j,u) \mathrel{\mathcal{L}^S} \Theta(j,v) \by{by~\cref{eq: D classes not twisted or 0-twist 1}} \,.
\end{align*}
\end{subequations}
Dually \((j,u) \mathrel{\mathcal{R}^T} (j,v) \iff \Theta(j,u) \mathrel{\mathcal{R}^S} \Theta(j,v)\), and hence \(D_{(j,a)}^T \simeq D_a^S\).

By \autoref{lem: twisted greens equivalences} it follows that \(D_a^S = D_a^{T_\Phi^0}\). Let \(u,v \in D_a^S\). By \autoref{lem: twisted greens equivalences} it is clear that \(u \mathrel{\mathcal{L}^S} v \iff u \mathrel{\mathcal{L}^{T_\Phi^0}} v\), and similarly \(u \mathrel{\mathcal{R}^S} v \iff u \mathrel{\mathcal{R}^{T_\Phi^0}} v\). Hence using the identity map gives \(D_a^S \simeq D_a^{T_\Phi^0}\), and by \autoref{lem: D class preserving} it follows that
\[
D_{(j,a)}^T \simeq D_a^S \simeq D_a^{T_\Phi^0} \,.
\]

As \(E(D_{(j,a)}^T) \neq \emptyset\), by \autoref{lem: E(T)} it follows that \(j \in E(M)\), and by \autoref{lem: E(T) <= E(S)} it is clear that \(E(D_a^S) \neq \emptyset\). Additionally, by \autoref{lem: E(T0) <= E(S)} it is clear that choosing \(e \in E(D_a^S)\) captures all possible idempotents in \(D_{(j,a)}^T\), \(D_a^{T_\Phi^0}\), and \(D_a^S\) (in the sense that any idempotent in \(D_{(j,a)}^T\) is of the form \((j,e)\), any idempotent in \(D_a^{T_\Phi^0}\) is of the form \(e\), and any idempotent in \(D_a^S\) is of the form \(e\), where in all cases \(e \in E(D_a^S)\)). This allows us to negate \cref{def: D class idempotent preserving ii} of \autoref{def: D class idempotent preserving} for \(D_a^S\) or \(D_a^{T_\Phi^0}\), possibly both. Consider the following cases.

\begin{mcases}

\case \label{cs: D classes not twisted or 0-twist 1}

Suppose \(\Phi(e,e) = 0\). Then trivially \(j = j + \Phi(e,e)q\) is satisfied and hence, by \autoref{lem: E(T)}, \((j,e)\) is an idempotent. Thus \(E(H_{(j,e)}^T) \neq \emptyset\) and
\begin{subequations}
\begin{align*}
E\big(\Theta(H_{(j,e)}^T)\big) &= E(H_e^S) \\
&= E(H_e^{T_\Phi^0}) \by{by \autoref{lem: twisted greens equivalences} and \(\Phi(e,e) = 0\)} \\
&\neq \emptyset \,.
\end{align*}
\end{subequations}
Hence for these \(\mathcal{H}\)-classes, \cref{def: D class idempotent preserving ii} is satisfied for both \(D_a^S\) and \(D_a^{T_\Phi^0}\).

\case \label{cs: D classes not twisted or 0-twist 2}

Suppose \(\Phi(e,e) > 0\). Then \(E(H_e^{T_\Phi^0}) = \emptyset\) but \(E(H_e^S) \neq \emptyset\), and hence \(D_a^S \not\cong D_a^{T_\Phi^0}\). This explains why the existence of such an \(e\) is equivalent to precisely one of the two conclusions of the theorem holding. Consider the following subcases.

\begin{msubcases}

\case \label{cs: D classes not twisted or 0-twist 2.1}

Suppose \(j = j + \Phi(e,e)q\). Hence by \autoref{lem: E(T)} it follows that \((j,e)\) is an idempotent. Thus \(E(H_{(j,e)}^T) \neq \emptyset\) and \(E\big(\Theta(H_{(j,e)}^T)\big) = E(H_e^S) \neq \emptyset\). Hence for these \(\mathcal{H}\)-classes, \cref{def: D class idempotent preserving ii} is satisfied for \(D_a^S\) but not for \(D_a^{T_\Phi^0}\), since \(E\big(\Theta(H_{(j,e)}^T)\big) = E(H_e^{T_\Phi^0}) = \emptyset\).

\case \label{cs: D classes not twisted or 0-twist 2.2}

Suppose \(j \neq j + \Phi(e,e)q\). Hence by \autoref{lem: E(T)} it follows that \((j,e)\) is not an idempotent. Thus \(E(H_{(j,e)}^T) = \emptyset\) and \(E\big(\Theta(H_{(j,e)}^T)\big) = E(H_e^{T_\Phi^0}) = \emptyset\). Hence for these \(\mathcal{H}\)-classes, \cref{def: D class idempotent preserving ii} is satisfied for \(D_a^{T_\Phi^0}\) but not for \(D_a^S\), as \(E\big(\Theta(H_{(j,e)}^T)\big) = E(H_e^S) \neq \emptyset\).

\end{msubcases}
\nomspace
\end{mcases}
\nomspace

By \autoref{con: exists implies for all idempotents} it is clear that, for all choices of \(e \in E(D_a^S)\), if \cref{cs: D classes not twisted or 0-twist 2} is triggered, then only one of \cref{cs: D classes not twisted or 0-twist 2.1} or \cref{cs: D classes not twisted or 0-twist 2.2} can occur (never a mixture for different choices of \(e\)). This is important: if only \cref{cs: D classes not twisted or 0-twist 2.1} occurs (in addition to possibly \cref{cs: D classes not twisted or 0-twist 1}), then all \(\mathcal{H}\)-classes satisfy \cref{def: D class idempotent preserving ii} for \(D_a^S\), and hence \(D_{(j,a)}^T \cong D_a^S\). Alternatively, if only \cref{cs: D classes not twisted or 0-twist 2.2} occurs (again possibly together with \cref{cs: D classes not twisted or 0-twist 1}), then all \(\mathcal{H}\)-classes satisfy \cref{def: D class idempotent preserving ii} for \(D_a^{T_\Phi^0}\), and hence \(D_{(j,a)}^T \cong D_a^{T_\Phi^0}\).
\end{proof}

\subsection{Representation theoretic interpretation}

We now have the following, slightly disappointing (from the viewpoint of representation theory), theorem.

\begin{Theorem}\label{T:Gram2}
Let $S$ be a planar diagram monoid as in \autoref{S:DiagCat} (or any finite monoid with trivial $\mathcal{H}$-classes). Assume \autoref{con: exists implies for all idempotents}. Any tight twisted version of it with $\mathcal{D}$-trivial $M$ has the same indexing set and dimensions for its simple representations as $S$ itself or the 0-twisted version.
\end{Theorem}

\begin{proof}
Apply \autoref{T:Gram}, \autoref{lem: implies same rank} and \autoref{thm: D classes not twisted or 0-twist}.
\end{proof}

Note that \autoref{T:Gram2} when combined with \autoref{SS:Cellular} determines a lot of the representation theory of these monoids.
Below we will use this to restrict to the study of the classical diagram monoids, and their counterparts where floating components are evaluated to zero.

\begin{Remark}
We currently know little from a representation-theoretic perspective about loose twistings, as they can disrupt the cellular structure.
\end{Remark}

\section{Diagrammatic Schur--Weyl dualities}\label{S:SchurWeyl}

We now collect a number of \emph{Schur--Weyl-type dualities} (all for $\ell=1$) in a form that is as uniform and general as possible. 
Many of them are not so easily found in the literature or only folk knowledge, in particular, in the generality we state, prove, and use them.

\begin{Notation}
In this paper, all representations are assumed to be finite dimensional over a given field and all modules are left modules, unless stated otherwise. Whenever appropriate, we restrict attention to type 1 and rational representations (or their natural analogs).
\end{Notation}

To set up terminology, we start with the basic example. Let $G = GL_m(\C)$ and, on the diagrammatic side, consider the symmetric group $S(n)$, with its usual wire diagram presentation, that we see as a partition of $\{1,\dots,n\}\cup\{1',\dots,n'\}$ with block sizes 2 containing one element from each of the two sets:
\begin{gather*}
\begin{tikzpicture}[anchorbase]
\draw[usual] (0,0) to (3,1);
\draw[usual] (1,0) to (1,1);
\draw[usual] (2,0) to (0,1);
\draw[usual] (3,0) to (2,1);
\end{tikzpicture}
\leftrightsquigarrow
\big\{\{1,4'\},\{2,2'\},\{3,1'\},\{4,3'\}\big\}.
\end{gather*}
We obtain this for $p=1,q=1-x$ in the language of \autoref{S:DiagCat}.

Let $V = \C^m$. Then:
\begin{enumerate}[label=(\Alph*)]

\item There are commuting actions of $G$ and $S(n)$ on $V^{\otimes n}$, and each action generates the full centralizer of the other. Concretely, we have surjective algebra homomorphisms (we will only use the latter)
\begin{gather*}
g \colon G \twoheadrightarrow \End_{S(n)}(V^{\otimes n}),\quad
f \colon \C[S(n)] \twoheadrightarrow \End_{G}(V^{\otimes n}),
\end{gather*}
which in particular imply that $\End_{G}(V^{\otimes n})$ admits a diagrammatic description in terms of a quotient of $S(n)$. This is often referred to as the \emph{first fundamental theorem}. Categorically, this can be expressed as the existence of a full functor
\begin{gather*}
F \colon \mathbf{S} \to \mathbf{Rep}(G)
\end{gather*}
from the category of symmetric groups $\mathbf{S}$ to the category of $G$-representations $\mathbf{Rep}(G)$
with all the expected properties (in particular, $F$ is monoidal).

\item Identifying the kernel of $f$ gives an explicit description of $\End_{G}(V^{\otimes n})$ as $S(n)/\ker f$. This is often referred to as the \emph{second fundamental theorem}. Categorically, this amounts to the existence of a fully faithful functor
\begin{gather*}
F \colon \mathbf{S} / \mathbf{ker}\,F \to \mathbf{Rep}(G)
\end{gather*}
again with all the desired properties (e.g.\ monoidal).

\item There is a decomposition of $V^{\otimes n}$ as a $\big(G\times S(n)\big)$-bimodule of the form
\begin{gather*}
V^{\otimes n}\cong\bigoplus_{\lambda}V(\lambda)\otimes L(\lambda),
\end{gather*}
where $V(\lambda)$ are simple $G$-modules, $L(\lambda)$ are simple $S(n)$-modules, and the sum runs over a suitable indexing set (making appearing modules pairwise distinct, e.g. $V(\lambda)\cong V(\mu)$ if and only if $\lambda=\mu$).
In particular, and we will use this below, the number of times 
$V(\lambda)$ appears as a summand in $V^{\otimes n}$ is $\dim[\C]L(\lambda)$, and vice versa. (This match of multiplicities and dimensions, as far as we know, has no name, so we just say \emph{by Schur--Weyl duality}.)

\end{enumerate}
With this terminology in place, we now explain the following Schur--Weyl duality table (with parameters specified later):
\begin{gather*}
\scalebox{0.98}{\begin{tabular}{c|c|c|c|c||c|c|c|c|c}
\arrayrulecolor{tomato}
Symbol & Diagrams & $G$ & $V$ & $\Z$?
& Symbol & Diagrams & $G$ & $V$ & $\Z$?
\\
\hline
\hline
$\ppamon[\mathbf{c}](t)$ & \begin{tikzpicture}[anchorbase]
\draw[usual] (0.5,0) to[out=90,in=180] (1.25,0.45) to[out=0,in=90] (2,0);
\draw[usual] (0.5,0) to[out=90,in=180] (1,0.35) to[out=0,in=90] (1.5,0);
\draw[usual] (0.5,1) to[out=270,in=180] (1,0.55) to[out=0,in=270] (1.5,1);
\draw[usual] (1.5,1) to[out=270,in=180] (2,0.55) to[out=0,in=270] (2.5,1);
\draw[usual] (0,0) to (0,1);
\draw[usual] (2.5,0) to (2.5,1);
\draw[usual,dot] (1,0) to (1,0.2);
\draw[usual,dot] (1,1) to (1,0.8);
\draw[usual,dot] (2,1) to (2,0.8);
\end{tikzpicture} & $SL2$ & $\C^2{\otimes}\C^2$ & Y
& $\pamon[\mathbf{c}](t)$ & \begin{tikzpicture}[anchorbase]
\draw[usual] (0.5,0) to[out=90,in=180] (1.25,0.45) to[out=0,in=90] (2,0);
\draw[usual] (0.5,0) to[out=90,in=180] (1,0.35) to[out=0,in=90] (1.5,0);
\draw[usual] (0,1) to[out=270,in=180] (0.75,0.55) to[out=0,in=270] (1.5,1);
\draw[usual] (1.5,1) to[out=270,in=180] (2,0.55) to[out=0,in=270] (2.5,1);
\draw[usual] (0,0) to (0.5,1);
\draw[usual] (1,0) to (1,1);
\draw[usual] (2.5,0) to (2.5,1);
\draw[usual,dot] (2,1) to (2,0.8);
\end{tikzpicture} & $S(t)$ & $\C^t$ & ?
\\
\hline
$\momon[\mathbf{c}](t)$ & \begin{tikzpicture}[anchorbase]
\draw[usual] (0.5,0) to[out=90,in=180] (1.25,0.5) to[out=0,in=90] (2,0);
\draw[usual] (1,0) to[out=90,in=180] (1.25,0.25) to[out=0,in=90] (1.5,0);
\draw[usual] (2,1) to[out=270,in=180] (2.25,0.75) to[out=0,in=270] (2.5,1);
\draw[usual] (0,0) to (1,1);
\draw[usual,dot] (2.5,0) to (2.5,0.2);
\draw[usual,dot] (0,1) to (0,0.8);
\draw[usual,dot] (0.5,1) to (0.5,0.8);
\draw[usual,dot] (1.5,1) to (1.5,0.8);
\end{tikzpicture} & $SL2$ & $\C^2{\oplus}\C$ & Y
& $\robrmon[\mathbf{c}](t)$ & \begin{tikzpicture}[anchorbase]
\draw[usual] (1,0) to[out=90,in=180] (1.25,0.25) to[out=0,in=90] (1.5,0);
\draw[usual] (1,1) to[out=270,in=180] (1.75,0.55) to[out=0,in=270] (2.5,1);
\draw[usual] (0,0) to (0.5,1);
\draw[usual] (2.5,0) to (2,1);
\draw[usual,dot] (0.5,0) to (0.5,0.2);
\draw[usual,dot] (2,0) to (2,0.2);
\draw[usual,dot] (0,1) to (0,0.8);
\draw[usual,dot] (1.5,1) to (1.5,0.8);
\end{tikzpicture} & $Ot,SPt$ & $\C^t$ & Y${}^\ast$
\\
\hline
$\tlmon[\mathbf{c}](t)$ & \begin{tikzpicture}[anchorbase]
\draw[usual] (0.5,0) to[out=90,in=180] (1.25,0.5) to[out=0,in=90] (2,0);
\draw[usual] (1,0) to[out=90,in=180] (1.25,0.25) to[out=0,in=90] (1.5,0);
\draw[usual] (0,1) to[out=270,in=180] (0.25,0.75) to[out=0,in=270] (0.5,1);
\draw[usual] (2,1) to[out=270,in=180] (2.25,0.75) to[out=0,in=270] (2.5,1);
\draw[usual] (0,0) to (1,1);
\draw[usual] (2.5,0) to (1.5,1);
\end{tikzpicture} & $SL2$ & $\C^2$ & Y
& $\brmon[\mathbf{c}](t)$ & \begin{tikzpicture}[anchorbase]
\draw[usual] (0.5,0) to[out=90,in=180] (1.25,0.45) to[out=0,in=90] (2,0);
\draw[usual] (1,0) to[out=90,in=180] (1.25,0.25) to[out=0,in=90] (1.5,0);
\draw[usual] (0,1) to[out=270,in=180] (0.75,0.55) to[out=0,in=270] (1.5,1);
\draw[usual] (1,1) to[out=270,in=180] (1.75,0.55) to[out=0,in=270] (2.5,1);
\draw[usual] (0,0) to (0.5,1);
\draw[usual] (2.5,0) to (2,1);
\end{tikzpicture} & $Ot,SPt$ & $\C^t$ & Y${}^\ast$
\\
\hline
$\promon[\mathbf{c}](t)$ & \begin{tikzpicture}[anchorbase]
\draw[usual] (0,0) to (0.5,1);
\draw[usual] (0.5,0) to (1,1);
\draw[usual] (2,0) to (1.5,1);
\draw[usual] (2.5,0) to (2.5,1);
\draw[usual,dot] (1,0) to (1,0.2);
\draw[usual,dot] (1.5,0) to (1.5,0.2);
\draw[usual,dot] (0,1) to (0,0.8);
\draw[usual,dot] (2,1) to (2,0.8);
\end{tikzpicture} & $GL2$ & $\C{\oplus}\C^{det}$ & Y
& $\romon[\mathbf{c}](t)$ & \begin{tikzpicture}[anchorbase]
\draw[usual] (0,0) to (1,1);
\draw[usual] (0.5,0) to (0,1);
\draw[usual] (2,0) to (2,1);
\draw[usual] (2.5,0) to (0.5,1);
\draw[usual,dot] (1,0) to (1,0.2);
\draw[usual,dot] (1.5,0) to (1.5,0.2);
\draw[usual,dot] (1.5,1) to (1.5,0.8);
\draw[usual,dot] (2.5,1) to (2.5,0.8);
\end{tikzpicture} & $GLt$ & $\C^t{\oplus}\C$ & Y
\\
\hline
$1$ & \begin{tikzpicture}[anchorbase]
\draw[usual] (0,0) to (0,1);
\draw[usual] (0.5,0) to (0.5,1);
\draw[usual] (1,0) to (1,1);
\draw[usual] (1.5,0) to (1.5,1);
\draw[usual] (2,0) to (2,1);
\draw[usual] (2.5,0) to (2.5,1);
\end{tikzpicture} & $SL2$ & $\C$ & Y
& $\sym[\mathbf{c}](t)$ & \begin{tikzpicture}[anchorbase]
\draw[usual] (0,0) to (1,1);
\draw[usual] (0.5,0) to (0,1);
\draw[usual] (1,0) to (1.5,1);
\draw[usual] (1.5,0) to (2.5,1);
\draw[usual] (2,0) to (2,1);
\draw[usual] (2.5,0) to (0.5,1);
\end{tikzpicture} & $GLt$ & $\C^t$ & Y
\end{tabular}}
\end{gather*}
For completeness, we also add the oriented Brauer category:
\begin{gather*}
\begin{tabular}{c|c|c|c|c}
\arrayrulecolor{tomato}
Symbol & Diagrams & $G$ & $V$ & $\Z$?
\\
\hline
\hline
$\obrmon[\mathbf{c}](t)$ & \begin{tikzpicture}[anchorbase]
\draw[usual,directed=1] (0.5,0) to[out=90,in=180] (1.25,0.45) to[out=0,in=90] (2,0);
\draw[usual,directed=1] (1.5,0) to[out=90,in=0] (1.25,0.25) to[out=180,in=90] (1,0);
\draw[usual,directed=1] (0,1) to[out=270,in=180] (0.75,0.55) to[out=0,in=270] (1.5,1);
\draw[usual,directed=1] (1,1) to[out=270,in=180] (1.75,0.55) to[out=0,in=270] (2.5,1);
\draw[usual,directed=1] (0,0) to (0.5,1);
\draw[usual,directed=1] (2,1) to (2.5,0);
\end{tikzpicture} & $GLt$ & $\C^{t}$,$(\C^{t})^\ast$ & Y${}^\ast$
\\
\end{tabular}
\end{gather*}
The table is to be read as follows:
\begin{enumerate}[label=(\alph*)]

\item The table illustrates the diagrammatic incarnations of representation theoretical objects. 

\item The first two columns are the diagram algebra and a prototypical example of a diagram.

\item The column $G$ shows the group underlying the duality, $V$ the used representation, and $\Z$? indicates whether the duality holds \emph{integrally} (meaning over $\Z$ or $\Z[v,v^{-1}]$ in the quantum case).

\item On the left side the first and second fundamental theorem hold spot on (i.e. one has equivalences), while on the right side only the first holds (i.e. there is a kernel in general). We also identify the kernel for the GLt rows in general.

\item The Y${}^\ast$ means it is integral for the group, as in the table, but there is a super version which is (probably) not integral.
(More on this in \autoref{R:Super} below.)
The symbol ? indicates that, at the time of writing, we do not know anything beyond characteristic zero.

\end{enumerate}
We now give the details. 

\subsection{Schur--Weyl dualities for planar diagram monoids}

We expect that the reader has some background on \emph{tilting modules} \cite{Do-tilting-alg-groups,Ri-good-filtrations} for quantum groups or reductive groups, see e.g.\ \cite{AnPoWe-representation-qalgebras,AnStTu-cellular-tilting,An-SimpleTLAll} (and the references therein).

Let $\mathbb{K}$ be a field and $q \in \mathbb{K}^{\ast}$. (To e.g.\ consider the case of a formal variable $q$ over $\C$, let $\K=\C(v)$ for a formal variable $v$ and set $q=v$.)  Let
$G = U_q(\mathfrak{sl}_2)$ be the quantum enveloping algebra of $\mathfrak{sl}_2$, in its divided power form, and $V=V(1)$ the standard two dimensional representation of $G$, and $V(0)$ the trivial representation. Let $W=V(1)\oplus V(0)$. For $X=V$ or $X=W$, define $\mathbf{F}_X(G)$ to be the full monoidal subcategory of $G$-modules with objects $X^{\otimes r}$, for $r\in\N$, and let 
$\mathbf{Tilt}(G)$ be the category of tilting modules. Both sit in the category of all finite dimensional $G$-modules $\mathbf{Rep}(G)$.

For the planar rook category we take $H = U_q(\mathfrak{gl}_2)$ with $V=V(0,0)\oplus V(1,1)$, meaning trivial plus the analog of the determinant representation, and define $\mathbf{F}_V(H)$ analogously as the full monoidal subcategory of $\mathbf{Rep}(H)$ with objects $V^{\otimes r}$, for $r\in\N$.

\begin{Lemma}\label{L:RigidBraided2}
The category $\mathbf{Tilt}(G)$ is the additive idempotent completion of $\mathbf{F}_V(G)$ or $\mathbf{F}_W(G)$, and all these categories inherit a braided pivotal structure from 
$\mathbf{Rep}(G)$. Moreover, the category $\mathbf{F}_V(H)$ is semisimple, hence idempotent complete, and inherits a rigid pivotal monoidal structure from 
$\mathbf{Rep}(H)$.
\end{Lemma}

\begin{proof}
For $G$ this is standard, see e.g. \cite{SuTuWeZh-mixed-tilting}. For $H$ this is also standard, but harder to find in the literature, so we spell it out.

Note that $V(0,0)$ and $V(1,1)$ are one dimensional characters, and hence every tensor power
\[
V^{\otimes r}
\cong \bigoplus_{k=0}^r V(k,k)^{\oplus \binom{r}{k}}
\]
is a direct sum of one dimensional simple $H$-modules, where $V(k,k)$ denotes the character obtained by taking $k$ copies of $V(1,1)$ and $r-k$ copies of $V(0,0)$. This decomposition is defined already over the Lusztig integral form over $\mathbb{Z}[v,v^{-1}]$, so it survives after any base change to $\K$ and any specialization $q\in\K^\ast$. Thus $\mathbf{F}_V(H)$ is semisimple and idempotent complete. The rest is clear.
\end{proof}

Considering now the table below, where $k\geq 2$:
\begin{gather*}
\begin{tabular}{c||c|c|c|c}
$\mathbf{M}$ & genus 0 & genus 1 & genus k & $V$ \\
\hline
\hline
$\mathbf{pPa}$ & $-q-q^{-1}$ & $(-q-q^{-1})^{2}$ & $(-q-q^{-1})^{k}$ & $V(1) \otimes V(1)$ \\
\hline 
$\mathbf{Mo}$ & Anything invertible & $1-q-q^{-1}$ & Not applicable & $V(1)\oplus V(0)$ \\
\hline
$\mathbf{TL}$ & Not applicable & $-q-q^{-1}$ & Not applicable & $V(1)$\\
\hline
$\mathbf{pRo}$ & Anything invertible & Not applicable & Not applicable & $V(1,1)\oplus V(0,0)$
\end{tabular}.
\end{gather*}

\begin{Theorem}\label{T:SchurWeylPlan}
There exists a fully faithful braided pivotal $\K$-linear functor 
\[
F\colon \mathbf{M}\longrightarrow \mathbf{F}_V(G)\subset\mathbf{Tilt}(G)
,\quad
F(r) = V^{\otimes r}
\]
that is an equivalence (with target $\mathbf{Tilt}(G)$) upon idempotent completion. Similarly, 
there exists a fully faithful monoidal $\K$-linear functor 
\[
F\colon \mathbf{M}\longrightarrow \mathbf{F}_V(H)
,\quad
F(r) = V^{\otimes r}
\]
that is an equivalence.
\end{Theorem}

\begin{proof}
The following fairly general strategy applies:
\begin{enumerate}
\item The properties of being braided pivotal $\K$-linear follow by construction (or are easy to verify).
\item We need to make sure that all categories and functors are defined integrally, i.e. over $\Z[v,v^{-1}]$ where $v$ is a formal variable. This is often easy, and just checking of the definitions.
\item We then want to use that:
\begin{enumerate}
\item For braided pivotal categories, bending 
\begin{gather*}
\mathrm{w} \mapsto
\begin{tikzpicture}[anchorbase,scale=1]
\draw[very thick] (0.25,1.5) to (0.25,2);
\draw[very thick] (1.25,1.5) to[out=90,in=270] (0.75,2);
\draw[very thick] (2.25,1.5) to[out=90,in=270] (1.25,2);
\draw[very thick,crossline] (1.75,1.5) to[out=90,in=270] (2.25,2);
\draw[very thick,crossline] (0.75,1.5) to[out=90,in=270] (1.75,2);
\draw[very thick] (0.25,-1) to (0.25,-0.5);
\draw[very thick] (1.25,-1) to[out=90,in=270] (0.75,-0.5);
\draw[very thick,crossline] (0.75,-1) to[out=90,in=270] (1.25,-0.5);
\draw[very thick] (1.75,-1) to (1.75,-0.5);
\draw[very thick] (2.25,-1) to (2.25,-0.5);
\draw[very thick] (0,0) rectangle node[pos=0.5]{w} (2.5,1);
\draw[very thick,directed=0.55] (0.25,1.5) to (0.25,1);
\draw[very thick,directed=0.55] (0.75,1) to (0.75,1.5);
\draw[very thick,directed=0.55] (1.25,1.5) to (1.25,1);
\draw[very thick,directed=0.55] (1.75,1) to (1.75,1.5);
\draw[very thick,directed=0.55] (2.25,1.5) to (2.25,1);
\draw[very thick,directed=0.55] (0.25,-0.5) to (0.25,0);
\draw[very thick,directed=0.55] (0.75,0) to (0.75,-0.5);
\draw[very thick,directed=0.55] (1.25,-0.5) to (1.25,0);
\draw[very thick,directed=0.55] (1.75,0) to (1.75,-0.5);
\draw[very thick,directed=0.55] (2.25,0) to (2.25,-0.5);
\draw[very thick] (0.25,2) to[out=90,in=0]
(0,2.25) to[out=180,in=90] (-0.25,2);
\draw[very thick] (0.75,2) to[out=90,in=0]
(0,2.5) to[out=180,in=90] (-0.75,2);
\draw[very thick] (1.25,2) to[out=90,in=0]
(0,2.75) to[out=180,in=90] (-1.25,2);
\draw[very thick] (2.25,-1) to[out=270,in=180]
(2.5,-1.25) to[out=0,in=270] (2.75,-1);
\draw[very thick] (1.75,-1) to[out=270,in=180]
(2.5,-1.5) to[out=0,in=270] (3.25,-1);
\draw[very thick] (1.25,-1) to[out=270,in=180]
(2.5,-1.75) to[out=0,in=270] (3.75,-1);
\draw[very thick] (1.75,2) to (1.75,3)node[above]{$l_1^{\prime}$};
\draw[very thick] (2.25,2) to (2.25,3)node[above]{$l_2^{\prime}$};
\draw[very thick] (2.75,-1) to (2.75,3)node[above]{$l_3^{\prime}$};
\draw[very thick] (3.25,-1) to (3.25,3)node[above]{$l_4^{\prime}$};
\draw[very thick] (3.75,-1) to (3.75,3)node[above]{$l_5^{\prime}$};
\draw[very thick] (-1.25,2) to (-1.25,-2)node[below]{$k_1^{\prime}$};
\draw[very thick] (-0.75,2) to (-0.75,-2)node[below]{$k_2^{\prime}$};
\draw[very thick] (-0.25,2) to (-0.25,-2)node[below]{$k_3^{\prime}$};
\draw[very thick] (0.25,-1) to (0.25,-2)node[below]{$k_4^{\prime}$};
\draw[very thick] (0.75,-1) to (0.75,-2)node[below]{$k_5^{\prime}$};
\end{tikzpicture},
\end{gather*}
shows knowing ordered endomorphism algebras with ordered objects is equivalent to knowing hom-spaces. (In our case we do not need the braiding since we only have one generating object.)
This is always true, e.g. see \cite[Proposition 2.10.8]{egno-tensor-2015} for a hint for a proof. The same works for subcategories of braided pivotal categories.
\item All categories are flat, i.e. everything is free over $\Z[v,v^{-1}]$. This is rare, but true for tilting modules.
\end{enumerate}
\item That the functors are full is independent of $\K$ and $q$, which is usually also easy.
\item Thus, knowing fully faithfulness for endomorphism spaces over $\C$ with $q=1$ implies it in general. Usually the case over $\C$ with $q=1$ is semisimple and well-known for donkey's years.
\end{enumerate}

Now the details. First, for $G$.
That there exists categories and functors that can be defined over $\Z[q,q^{-1}]$ follows from \cite{RTW-Valenztheorie} for $\mathbf{M} \in \{\mathbf{TL},\mathbf{pPa}\}$ and from 
\cite{BH-MotzkinAlgebras} for $\mathbf{M} = \mathbf{Mo}$ (using also \cite{kaneda}). For these (a) holds, essentially by definition, and (d) is easy, e.g. by the main result of \cite{AST-CellularStructures}.
Next, (c).(i), imported via \autoref{L:RigidBraided} and \autoref{L:RigidBraided2}, allows us to identify hom-space and end-space knowledge, and we will do this freely now.
Moreover, all categories are flat. For the diagram categories this follows by birth, while for the representation theory this is \cite[Corollary 3.14]{AST-CellularStructures}.
Finally, we only need to know the fully faithfulness results for end-spaces, which is again \cite{RTW-Valenztheorie} and \cite{BH-MotzkinAlgebras}.

For $H$ the argument is essentially the same but even easier. By the proof of \autoref{L:RigidBraided2} we have, for every $n\ge 0$,
\[
\dim \End_H\!\big(V^{\otimes n}\big)
\;=\; \sum_{k=0}^n \binom{n}{k}^2,
\]
independently of $\K$ and $q$.  On the planar rook side, the endomorphism algebra of $n$ points in $\mathbf{M}=\mathbf{pRo}$ has the same dimension, which follows essentially by definition. The same works for hom-spaces, by pivotality of the parent categories.
The functor $F\colon \mathbf{M}\to\mathbf{F}_V(H)$ sends the generating object to $V$ and the generating diagrams to the inclusions and projections onto the summands $V(0,0)$ and $V(1,1)$.  This induces a surjection
\[
\Hom_{\mathbf{M}}(n,m)\twoheadrightarrow \Hom_H\!\big(V^{\otimes n},V^{\otimes m}\big)
\]
defined over $\mathbb{Z}[v,v^{-1}]$ and hence over any specialization.  Since source and target have the same dimension, the map is an isomorphism, so $F$ is fully faithful on all hom-spaces. The rest follows swiftly.
\end{proof}

\begin{Remark}
We do not know any Schur--Weyl-type duality for the planar partition algebra for non-uniform evaluations, and the Motzkin and planar rook algebra for interval scalar zero.
\end{Remark}

\subsection{Schur--Weyl dualities for symmetric diagram monoids}\label{SS:SymSW}

Given an infinite field $\K$, let 
$G=GL_t(\K)$ denote the general linear group, and 
let $V=V(\omega_1)$ be its defining $t$ dimensional representation. We look at the category of $G$-modules $\mathbf{Rep}(G)$.
As before, define $\mathbf{F}_V(G)$ as the full monoidal subcategory of $\mathbf{Rep}(G)$ with objects $V^{\otimes r}$, for $r\in\N$. Similarly, for $W=V(\omega_1)\oplus V(0)$ (defining representation plus the trivial representation) we have 
$\mathbf{F}_W(G)$. We also have tilting modules $\mathbf{Tilt}(G)\subset\mathbf{Rep}(G)$, as before.

\begin{Lemma}\label{L:RigidBraided3}
The category $\mathbf{Tilt}(G)$ is the additive idempotent completion of $\mathbf{F}_V(G)$ or $\mathbf{F}_W(G)$, and all these categories inherit a symmetric structure from 
$\mathbf{Rep}(G)$.
\end{Lemma}

\begin{proof}
The nontrivial parts follow from \cite{Do-tilting-alg-groups,Ri-good-filtrations}.
\end{proof}

Now the table is:
\begin{gather*}
\begin{tabular}{c|c|c}
$\mathbf{M}$ & $\delta$ & $V$ \\
\hline
\hline
$\mathbf{Ro}$ & Anything invertible & $V(\omega_1)\oplus V(0)$\\
\hline
$\mathbf{S}$ & Not applicable & $V(\omega_1)$
\end{tabular}.
\end{gather*}
Let $\mathbf{I}$ be the $\circ$-$\otimes$-ideal in $\mathbf{S}$ with generating set given by the $\K$-linear span of the Kazhdan--Lusztig basis element for the longest word in $S(t+1)$. (In characteristic zero, this Kazhdan--Lusztig basis element is the sum of all elements in $S(t+1)$ divided by $(t+1)!$.)

\begin{Theorem}\label{T:SchurWeylPlan2}
For $X=V$ or $X=W$, there exists a fully faithful symmetric $\K$-linear functor 
\[
F\colon \mathbf{M}/\mathbf{I}\longrightarrow \mathbf{F}_X(G)\subset\mathbf{Tilt}(G)
,\quad
F(r)=X^{\otimes r}
\]
that is an equivalence upon idempotent completion.
\end{Theorem}

\begin{proof}
The same strategy as in \autoref{T:SchurWeylPlan}
works: the statement for $\K=\C$ is classical Schur--Weyl duality (for completeness, see also \cite{AM-SchurRook}), and for general $\K$ we use \cite{DuPaSc-Schur--Weyl} (for non algebraically closed but infinite fields their arguments still work, cf. \cite{Do-tilting-alg-groups,Ri-good-filtrations}).
\end{proof}

\begin{Remark}\label{R:Super}
There is also a version of \autoref{T:SchurWeylPlan2} for the special linear group instead of the general linear group, and also a quantum version. All of these work verbatim. There is also a super version, which, however, does not work in positive characteristic without modification, cf. \cite{CEKO-InvPosChar}. Moreover, for Brauer and rook Brauer below, we directly use superalgebras because we will need them as they allow the correct parameter specializations for the monoid case. Following \cite{CEKO-InvPosChar}, we expect that issues arise in positive characteristic here as well.
\end{Remark}

Note that, \autoref{T:SchurWeylPlan2} implies that there is an $\K$-algebra isomorphism
\begin{gather*}
\K S(r)\cong\End_G(X^{\otimes r})\text{ for }r\leq t+1.
\end{gather*}

Let us now focus on $\mathbf{Pa}$, where $k\in\N$ and the table is:
\begin{gather*}
\begin{tabular}{c||c|c}
$\mathbf{M}$ & genus k & $V$ \\
\hline
\hline
$\mathbf{Pa}$ & $r$ & $V_{S(r)}$ \\
\end{tabular}
.
\end{gather*}
Here $V_{S(r)}$ is the $r$ dimensional permutation representation of $S(r)$. Let $\mathbf{I}$ be the $\circ$-$\otimes$-ideal in $\mathbf{Pa}$ given by all idempotent projectors corresponding to partitions with $>r$ parts.

\begin{Theorem}
For $X=V_{S(r)}$, there exists a full symmetric $\K$-linear functor 
\[
F\colon \mathbf{M}\longrightarrow \mathbf{Rep}\big(S(r)\big)
,\quad
F(r)=X^{\otimes r}.
\]
If $\K$ is of characteristic zero, then the kernel of this functor is $\mathbf{I}$.
\end{Theorem}

\begin{proof}
The usual strategy applies without tilting modules (and characteristic dependence). The references are \cite{HaRa-partition-algebras,BDM} for the first fundamental theorem, and \cite{CO} for the second (the latter is only characteristic zero).
\end{proof}

We now turn to Lie supergroups. All relevant background can be found in, e.g., \cite{LZ-FirstInvOrtho}. Let $\mathbb{K}$ be a field of characteristic zero. Considering the table below, with $V=V(\omega_1)$ denoting the orthosymplectic superspace of dimension $(m|2n)$ for $G=OSp(V)$. Moreover, let $V(0)$ denote the trivial representation of $G$.
\begin{gather*}
\begin{tabular}{c||c|c|c}
$\mathbf{M}$ & genus 0 & genus 1 & $V$ \\
\hline
\hline
$\mathbf{RoBr}$ & Anything invertible & $m+1-2n$ & $V(\omega_1)\oplus V(0)$ \\
\hline
$\mathbf{Br}$ & Not applicable & $m-2n$ & $V(\omega_1)$ \\
\end{tabular}
.
\end{gather*}
We use the same notation as above for $\mathbf{F}_X(G)$.

\begin{Theorem}\label{T:SchurWeylSym}
For $X=V(\omega_1)$ or $X=V(\omega_1)\oplus V(0)$,
there exists a full symmetric $\K$-linear functor
\[
F\colon \mathbf{M}\longrightarrow \mathbf{F}_V(G)
,\quad
F(r) = V^{\otimes r}.
\]
Moreover, the functor is injective on the homomorphism algebras $\textnormal{Hom}_M(a,b)$ for the following conditions:
\begin{gather*}
\begin{tabular}{c|c}
$\mathbf{M}$ & Condition \\
\hline 
\hline
$\mathbf{RoBr}$ & $a+b < (m+2)(2n+1)$ \\
\hline
$\mathbf{Br}$ & $a+b < (m+1)(2n+1)$ \\
\end{tabular}
.
\end{gather*}
If either $m=0$ or $n=0$, then the same works for any infinite field.
\end{Theorem}

\begin{proof}
We state and prove two key results, \autoref{P:BrSchurWeyl} and \autoref{P:RoBrSchurWeyl}, below. Then, the theorem follows from these propositions, \autoref{L:RigidBraided}, and \cite[Proposition 2.10.8]{egno-tensor-2015}. The remaining statement for $m=0$ or $n=0$ is then \autoref{Le:IntOSP}.
\end{proof}

\begin{Proposition}\label{P:BrSchurWeyl}
Over a field $\mathbb{K}$ of characteristic zero, $\delta = m-2n$, and $2d < (m+1)(2n+1)$,
\[
Br_{\delta}(d) \cong \textnormal{End}_{OSp(V)}(V^{\otimes d}),
\]
where $V$ is the $(m|2n)$-superdimensional orthosymplectic superspace and $OSp(V)$ is the orthosymplectic supergroup.
\end{Proposition}

\begin{proof}
This follows from \cite[Corollary 5.5 and Corollary 5.8]{LZ-SecondInvOrtho}.
\end{proof}

Let $RoBr_{\delta}(d)$ denote the rook Brauer algebra 
with invertible genus 0 scalar, and genus 1 scalar $\delta$.

\begin{Proposition}\label{P:RoBrSchurWeyl}
For $2d < (m+2)(2n+1)$ and $\delta = m+1-2n \neq 0$,
\[
RoBr_{\delta}(d) \cong \End_{OSp(V)}(V^{\otimes d}).
\]
\end{Proposition}

\begin{proof}
Let $\K$ be a field of characteristic $0$. Let $V(1)_{\K}$ be the $\Z/2\Z$-graded vector space from \cite[Section~2.1]{LZ-FirstInvOrtho} and put
\[
V_{\K} = V(1)_{\K} \oplus \K .
\]
Following \cite[Section~2.3]{LZ-FirstInvOrtho}, equip $V_{\K}$ with a supersymmetric bilinear form 
\[
\langle \placeholder,\placeholder \rangle : V_{\K} \times V_{\K} \longrightarrow \K
\]
with matrix
\[
\eta = 
\begin{pmatrix}
I_{m+1} & 0 \\
0 & J
\end{pmatrix},
\quad
J = \mathrm{diag}(\sigma,\dots,\sigma),
\quad
\sigma = 
\begin{pmatrix}
0 & 1 \\
-1 & 0
\end{pmatrix},
\]
where $I_{m+1}$ is the $(m+1)\times (m+1)$ identity matrix and $J$ is $2n\times 2n$. Denote the $(a,b)$-entry of $\eta$ by $\eta_{ab}$ and the entry of $\eta^{-1}$ by $\eta^{ab}$. Fix a homogeneous basis $(e_1,\dots,e_{m+1+2n})$ of $V_{\K}$ such that
\[
\langle e_a,e_b\rangle = \eta_{ab}
\quad\text{for all }a,b.
\]
Let $\Lambda$ be the Grassmann algebra used in \cite[Section~2.3]{LZ-FirstInvOrtho}, and set
\[
V = V_{\K} \otimes \Lambda .
\]
Then $V$ is the orthosymplectic superspace of superdimension $(m+1 \mid 2n)$ from loc.\ cit. Let $G = OSp(V)$ and set $\delta = m+1-2n$.

Define the maps
\begin{gather}\label{P-C-C}
\begin{aligned}
&\tau: V\otimes V\longrightarrow V\otimes V,
&& v\otimes w \longmapsto (-1)^{[v][w]}w\otimes v, \\
&\check{C}: \Lambda \longrightarrow V\otimes V,
&& 1\longmapsto \sum_{a,b=1}^{m+1+2n} e_a \otimes \eta^{ab} e_b, \\
&\hat{C}: V\otimes V \longrightarrow \Lambda,
&& v\otimes w\longmapsto \langle v, w\rangle,\\
&\check{D}: \Lambda \longrightarrow V,
&& 1\longmapsto \sum_{a=1}^{m+1+2n}\lambda_a e_a, \\
&\hat{D}: V \longrightarrow \Lambda,
&& v=\sum_{i=1}^{m+1+2n} \mu_i e_i \longmapsto \langle \mu^* v,\sum_{a=1}^{m+1+2n} \eta^{ab} e_a \rangle,
\end{aligned}
\end{gather}
where $\{\lambda_a\}$ are the Grassmann generators and $\mu^* = 1/\mu_i$ if $\mu_i\neq 0$ and $0$ otherwise (see \cite[Section~5]{LZ-FirstInvOrtho} for the standard setup).

The following is the analog of \cite[Lemma~5.3]{LZ-FirstInvOrtho} in the present notation.

\begin{Lemma} \label{lem:PAU}
The maps $\tau$, $\check{C}$, $\hat{C}$, $\check{D}$ and $\hat{D}$ are all $OSp(V)$-equivariant and satisfy:
\begin{align}
&\tau^2=\id_V^{\otimes 2}, \quad
(\tau\otimes \id_V)(\id_V\otimes \tau)(\tau\otimes \id_V) 
= (\id_V\otimes \tau)(\tau\otimes\id_V)(\id_V\otimes \tau),
\label{eq:PPP}\\[0.3em]
&\tau \check{C} =  \check{C}, \quad
\hat{C} \tau =   \hat{C}, \label{eq:fse-es}\\[0.3em]
&\hat{C}\check{C}=\delta,  \quad
(\hat{C}\otimes\id_V)(\id_V\otimes\check{C})
=\id_V
=(\id_V\otimes\hat{C})(\check{C}\otimes\id_V),
\label{eq:CC}\\[0.3em]
&(\hat{C}\otimes\id_V)\circ (\id_V\otimes \tau)
= (\id_V\otimes\hat{C})\circ (\tau\otimes\id_V),
\label{eq:CPC-I}\\[0.3em]
&(\tau\otimes \id_V)\circ(\id_V\otimes \check{C})
=(\id_V\otimes \tau)\circ(\check{C}\otimes \id_V),
\label{eq:CPC-P}\\[0.3em]
&\hat{D}\check{D}
= \hat{C} \circ (\check{D}\otimes \check{D})
= (\hat{D}\otimes\hat{D})\circ \check{C}
= \delta,
\label{eq:DD}\\[0.3em]
&(\hat{D}\otimes\id_V)\circ\check{C} 
= (\id_V\otimes\hat{D})\circ\check{C} 
= \check{D}, \quad
\hat{C}\circ(\check{D}\otimes\id_V) 
= \hat{C}\circ(\id_V\otimes\check{D}) 
= \hat{D}. \label{eq:DC}
\end{align}
\end{Lemma}

\begin{proof}
The relations \eqref{eq:PPP}-\eqref{eq:CPC-P} are as in \cite[Lemma~5.3]{LZ-FirstInvOrtho}. For \eqref{eq:DD} we compute
\[
\hat{D}\check{D}(1)
= \hat{D}\Big(\sum_{a=1}^{m+1+2n}\lambda_a e_a\Big)
= \Big\langle \sum_{\lambda_a \neq 0} e_a,\sum_{b=1}^{m+1+2n}\eta^{ab}e_b \Big\rangle
= \sum_{a,b=1}^{m+1+2n}\eta^{ab} \langle e_a,e_b \rangle
= \sum_{a,b=1}^{m+1+2n}\eta^{ab}\eta_{ab}
= \delta.
\]
The remaining identities follow in the same way by a straightforward contraction computation.
\end{proof}

The Rook Brauer category $\mathbf{RoBr}_\delta$ is generated by the usual Rook Brauer diagrams, i.e.:
\[
\begin{tikzpicture}[anchorbase]
\draw[usual] (0,0.5) to (0,0);
\end{tikzpicture}
,\quad
\begin{tikzpicture}[anchorbase]
\draw[usual] (0,0.5) to (0.5,0);
\draw[usual] (0.5,0.5) to (0,0);
\end{tikzpicture}
,\quad
\begin{tikzpicture}[anchorbase]
\draw[usual] (0,0.5) to[out=270,in=270] (0.5,0.5);
\end{tikzpicture}
,\quad
\begin{tikzpicture}[anchorbase]
\draw[usual] (0,0) to[out=90,in=90] (0.5,0);
\end{tikzpicture}
,\quad
\begin{tikzpicture}[anchorbase]
\draw[usual,dot] (0,0.25) to (0,0);
\end{tikzpicture}
,\quad
\begin{tikzpicture}[anchorbase]
\draw[usual,dot] (0,0) to (0,0.25);
\end{tikzpicture}
.
\]
(The identity is included for completness only.)

\begin{Proposition}
There is a full functor
\[
F : \mathbf{RoBr}_{\delta} \longrightarrow \End_G\!(V^{\otimes d})
\]
sending objects $r$ to $V^{\otimes r}$ and a morphism $A : a \to b$ to a map
\[
F(A) : V^{\otimes a} \longrightarrow V^{\otimes b},
\]
defined on the generating diagrams by
\begin{gather}\label{eq:F-generating}
\begin{aligned}
F\left(
\begin{tikzpicture}[anchorbase]
\draw[usual] (0,0.5) to (0,0);
\end{tikzpicture}\right)&=\id_V,
&
F\left(
\begin{tikzpicture}[anchorbase]
\draw[usual] (0,0.5) to (0.5,0);
\draw[usual] (0.5,0.5) to (0,0);
\end{tikzpicture}\right) &= \tau, \\
F\left(
\begin{tikzpicture}[anchorbase]
\draw[usual] (0,0.5) to[out=270,in=270] (0.5,0.5);
\end{tikzpicture}\right) &= \check{C},
&
F\left(
\begin{tikzpicture}[anchorbase]
\draw[usual] (0,0) to[out=90,in=90] (0.5,0);
\end{tikzpicture}\right) &= \hat{C}, \\
F\left(
\begin{tikzpicture}[anchorbase]
\draw[usual,dot] (0,0.25) to (0,0);
\end{tikzpicture}\right) &= \check{D},
&
F\left(
\begin{tikzpicture}[anchorbase]
\draw[usual,dot] (0,0) to (0,0.25);
\end{tikzpicture}\right) &= \hat{D}.
\end{aligned}
\end{gather}
Moreover, $F$ is symmetric monoidal, that is, it respects tensor products.
\end{Proposition}

\begin{proof}
Using \autoref{lem:PAU}, the arguments of \cite[Section~5]{LZ-FirstInvOrtho} imply this statement.
\end{proof}

Finally, applying the arguments of \cite[Section~5]{LZ-SecondInvOrtho} to this situation gives \autoref{P:RoBrSchurWeyl}.
\end{proof}

\begin{Lemma}\label{Le:IntOSP}
Assume $m=0$ or $n=0$ and let $\K$ be an infinite field (with $\mathrm{char}(\K)\neq 2$ in the orthogonal case). Then the assertions of \autoref{T:SchurWeylSym} remain valid over $\K$.
\end{Lemma}

\begin{proof}
If $n=0$, then $V$ is an ordinary orthogonal space of dimension $m$ and $OSp(V)\cong O_m$.  
If $m=0$, then $V$ is an ordinary symplectic space of dimension $2n$ and $OSp(V)\cong Sp_{2n}$.  
In both cases the super structure disappears and \autoref{T:SchurWeylSym} reduces to the usual Schur--Weyl--Brauer duality for types $B,C$ and $D$.

More precisely, \cite[Proposition~2.3]{AST} shows that $V^{\otimes d}$ is a tilting $G$-module and also that $\dim\End_{G}(V^{\otimes d})$ is independent of the ground field. The Schur--Weyl--Brauer duality of \cite[Theorem~3.17]{AST} identifies the Brauer algebra $Br_\delta(d)$ with $\End_G(V^{\otimes d})$ (with the stated bounds on $d$) over any field of the allowed characteristic. This gives the Brauer row of \autoref{T:SchurWeylSym} over arbitrary infinite $\K$.

The categorical formulation (fullness on homomorphism algebras and injectivity in the stated range) then follows from the same bending argument as in the proof of \autoref{T:SchurWeylPlan}. Finally, the proof of \autoref{P:RoBrSchurWeyl} is purely algebraic once the Brauer case is known, so it carries over verbatim to any infinite field, giving the rook Brauer row as well. 
\end{proof}

\begin{Remark}
We do not need the Schur--Weyl duality for the oriented Brauer algebra in this paper, but it can be established and reads similarly to \autoref{T:SchurWeylSym}, see \cite{Ko,Tur,AST} for the underlying dualities.
\end{Remark}

\begin{Remark}
In all the above, we do need an infinite field. For example, when $\K$ is a finite field $SL_2(\K)$ is a finite group and things change, cf. \cite{BeDo}.
\end{Remark}

\section{Typical highest weights in Schur--Weyl dualities}\label{S:TypicalHW}

We record a \emph{unimodality} or \emph{concentration phenomenon} for tensor powers of representations in (semisimple) diagrammatic Schur--Weyl settings, and in fact in far more general contexts. 
Its dual manifestation is the clustering of high dimensional representations of the associated diagram algebras. 
For $SL_2(\C)$ and the planar diagram algebras, what we obtain is genuine unimodality (hence the name); in other cases it is better viewed as a Gaussian profile on the weight lattice.

We assume the reader has some basic familiarity with \emph{random walks}, \emph{Markov chains}, and \emph{big~O (Bachmann--Landau) notation}. Standard references include \cite{LawLim-RW} for random walks and the central/local limit theorems, and \cite{LPW-MC} for Markov chains and mixing times. If the reader wants only diagram-algebra consequences, they can treat the results as black boxes.

\subsection{The group side}

For the remainder of this section we work over $\C$ and fix the following notation.
Let $G$ be a connected complex semisimple algebraic group with maximal torus $T$ and Weyl group $W$. Let $X$ be the weight lattice, $X^+\subset X$ the set of dominant integral weights, and for $\lambda\in X^+$ denote by $V(\lambda)$ the simple $G$-module of highest weight $\lambda$.  We fix a $W$-invariant inner product $\langle\placeholder,\placeholder\rangle$ on $X\otimes_\Z\R$ and write $\|x\|=\sqrt{\langle x,x\rangle}\in\R_{\geq 0}$.

Let $V$ be a $G$-module, and for $n\geq 1$ write its tensor powers as
\[
V^{\otimes n}
\;\cong\;
\bigoplus_{\lambda\in X^+} V(\lambda)^{\oplus m(\lambda,n)},
\quad \text{for }m(\lambda,n)\in\N,
\]
where $V(\lambda)$ is the $G$-module of highest weight $\lambda$.
We regard this as a random highest weight by putting the probability measure
\begin{gather}\label{eq:Pn-def}
\mathbb{P}_n(\lambda)=
\frac{m(\lambda,n)\,\dim[\C] V(\lambda)}{(\dim[\C] V)^n},
\quad\text{for }
\lambda\in X^+,
\end{gather}
so that
\[
\sum_{\lambda\in X^+} \mathbb{P}_n(\lambda)=1.
\]
We then have the following (the proof is a bit further down):

\begin{Theorem}\label{T:TypicalHW}
With notation as above, there exist:
\begin{itemize}
\item a vector $\mu=\mu(G,V)\in X\otimes_\Z\R$, the \emph{typical highest weight}, and
\item a positive definite quadratic form $q=q(G,V)$ on $X\otimes_\Z\R$,
\end{itemize}
such that the following hold.
\begin{enumerate}[label=\textup{(\alph*)}]
\item We have \emph{Gaussian concentration}:
For every $\varepsilon\in\R_{>0}$ there exists $C\in\R_{>0}$ such that, for all $n$ large enough,
\[
\sum_{\|\lambda-n\mu\|\;\geq\; C\sqrt{n}} \mathbb{P}_n(\lambda)
\;\leq\;
\varepsilon.
\]
In particular, almost all of the mass of $\mathbb{P}_n$ is contained in a window of size $O(\sqrt{n})$ around $n\mu$.

\item Let $r=\mathrm{rank}(G)$. If $(\lambda_n)_{n\geq 1}$ is any sequence in $X^+$ with
\[
\|\lambda_n-n\mu\| \in O(\sqrt{n})\quad\text{for }n\gg 0,
\]
then there exists a constant $c_{G,V}\in\R_{>0}$ such that
\[
\mathbb{P}_n(\lambda_n)
=
c_{G,V}\,n^{-r/2}\,
\exp\!\Bigl(-\frac{1}{2n}\,q(\lambda_n-n\mu)\Bigr)\bigl(1+o(1)\bigr),
\quad\text{for }n\gg 0.
\]
\end{enumerate}
\end{Theorem}

In particular, among all isotypic components of $V^{\otimes n}$, those with highest weights in a $\sqrt{n}$-window around $n\mu$ contribute asymptotically all of the total dimension, whereas highest weights far away from $n\mu$ appear with exponentially small total weight. This is the unimodality/concentration phenomenon we will use below. Note that $\mu\in X\otimes_\Z\R$ and not necessarily $\mu\in X$, but we can (without too much harm) and will think of $\mu$ as an honest highest weight.

\begin{Example}\label{Ex:TypicalGL}
Take $G=SL_t(\C)$, and let $V=V(\omega_1)\cong \C^t$ be the defining representation.  Then
\[
V^{\otimes n}
\;\cong\;
(\C^t)^{\otimes n}
\;\cong\;
\bigoplus_{\substack{\lambda\vdash n\\ \ell(\lambda)\leq t}}
V(\lambda)\otimes L(\lambda),
\]
where the sum runs over partitions $\lambda$ of $n$ (denoted as usual by $\lambda\vdash n$) with at most $t$ parts (denoted by $\ell(\lambda)\leq t$), and $L(\lambda)$ is the corresponding simple $\C[S(n)]$-module. In particular,
\begin{gather*}
m(\lambda,n)=\dim[\C] L(\lambda)\text{ is given by the hook length formula},
\\
\dim[\C] V(\lambda)\text{ is given by Weyl's dimension formula},
\end{gather*}
(so both are well-known)
and the probability measure from \eqref{eq:Pn-def} becomes
\[
\mathbb{P}_n(\lambda)
=
\frac{\dim[\C] L(\lambda)\,\dim[\C] V(\lambda)}{t^n},
\quad
\lambda\vdash n,\ \ell(\lambda)\leq t.
\]
\autoref{T:TypicalHW} then says that there is a limit shape $\mu=\mu(t)$ such that, when we view Young diagrams as lattice points in the dominant chamber,
\begin{itemize}
\item almost all of the mass of $\mathbb{P}_n$ is carried by partitions whose diagrams lie within an $O(\sqrt{n})$ neighborhood of the dilated shape $n\mu$, and
\item inside this $\sqrt{n}$-window the weights $\mathbb{P}_n(\lambda)$ are asymptotically given by a discrete Gaussian in the deviations $\lambda-n\mu$.
\end{itemize}
Equivalently, among all isotypic components $V(\lambda)\otimes L(\lambda)$ of $(\C^t)^{\otimes n}$, those with $\lambda$ close to $n\mu$ (in particular, the \emph{``typical staircase-like shapes''} one draws in small rank) occur with overwhelmingly larger total dimension than those far away from $n\mu$, whose total contribution is exponentially small.

For $t=2$, this recovers the familiar picture for $SL_2(\C)$: partitions with at most two parts correspond to natural numbers (the difference between the two rows), and the distribution of highest weights in $V^{\otimes n}$ is essentially binomial, concentrating in a window of size $O(\sqrt{n})$ around the central weight.
\end{Example}

\begin{proof}[Proof of \autoref{T:TypicalHW}]
Write
\[
V \;\cong\; \bigoplus_{i=1}^k V(\nu_i)^{\oplus a_i},
\quad
a_i\in\N,\ \nu_i\in X^+,
\]
and set
\[
\mu_1(\nu_i)
= \frac{a_i\,\dim[\C] V(\nu_i)}{\dim[\C] V},
\quad
i\in\{1,\dots,k\}.
\]
Then $\mu_1$ is a probability measure on $X^+$ describing the highest weight of $V$ with respect to the normalization \eqref{eq:Pn-def} for $n=1$, i.e.\ it is the law of the random variable
\[
\Lambda_1 \in X^+,
\quad
\mathbb{P}(\Lambda_1=\nu_i)=\mu_1(\nu_i).
\]

\textit{Step 1: random walk on the weight lattice.}
Consider now an iid\ sequence $(\Lambda_j)_{j\geq 1}$ with law $\mu_1$ and define a random walk
\[
W_n = \Lambda_1 + \cdots + \Lambda_n \in X\otimes_\Z\R,
\quad n\geq 1.
\]
Biane's construction \cite{Bi-asymptotic-lie} identifies the tensor powers $V^{\otimes n}$ with the convolution powers $\mu_1^{\ast n}$ via characters: the character of $V^{\otimes n}$ is the $n$-fold product of the character of $V$, hence the Fourier transform of $\mu_1^{\ast n}$. After applying the Weyl character formula and the usual $\rho$-shift/projection to the dominant chamber, the distribution of the highest weight of a simple constituent of $V^{\otimes n}$ can be expressed in terms of the distribution of $W_n$ and the Weyl group action.

More concretely, if we write
\[
V^{\otimes n}
\cong
\bigoplus_{\lambda\in X^+} V(\lambda)^{\oplus m(\lambda,n)},
\]
then the multiplicities $m(\lambda,n)$ admit an explicit expression of the form
\[
m(\lambda,n)
=
(\dim[\C] V)^n\,\big(\dim[\C] V(\lambda)\big)^{-1}\,
\mathrm{dens}_n(\lambda+\rho),
\]
where $\rho$ is the half-sum of positive roots and $\mathrm{dens}_n$ is the discrete density of a random walk on the full weight lattice (given by $\mu_1^{\ast n}$), corrected by Weyl group signs. Using \eqref{eq:Pn-def} this becomes
\[
\mathbb{P}_n(\lambda)
=
\frac{m(\lambda,n)\,\dim[\C] V(\lambda)}{(\dim[\C] V)^n}
=
\mathrm{dens}_n(\lambda+\rho),
\]
so the law of the random highest weight $\Lambda_n$ under $\mathbb P_n$ coincides with the (projected) random walk density.

\textit{Step 2: Biane's central and local limit theorems.}
Biane's main result is a central limit theorem and a local limit theorem for these projected random walks in the dominant chamber; see \cite{Bi-asymptotic-lie}, or \cite[Theorem~8]{TZ} for a convenient restatement. In our notation, this gives the following:
\begin{itemize}
\item There exists a drift vector $\mu\in X\otimes_\Z\R$ and a positive definite quadratic form $q$ on $X\otimes_\Z\R$ such that
\[
\frac{W_n - n\mu}{\sqrt{n}}
\]
converges in distribution to a non-degenerate Gaussian on $X\otimes_\Z\R$, and the same holds after the projection to the dominant chamber (i.e.\ for $\Lambda_n$ instead of $S(n)$).
\item Moreover, there is a local limit theorem in the central limit regime: for weights $\lambda=\lambda_n$ with $\|\lambda_n-n\mu\|\in O(\sqrt{n})$, the discrete densities $\mathrm{dens}_n(\lambda_n+\rho)$ admit Gaussian asymptotics of the form
\[
\mathrm{dens}_n(\lambda_n+\rho)
=
c_{G,V}\,n^{-r/2}\,
\exp\!\big(-\tfrac{1}{2n}\,q(\lambda_n-n\mu)\big)\bigl(1+o(1)\bigr),
\]
where $r=\mathrm{rank}(G)$ and $c_{G,V}\in\R_{>0}$ is a constant depending only on $(G,V)$.
\end{itemize}
Since, as noted above, $\mathbb P_n(\lambda)=\mathrm{dens}_n(\lambda+\rho)$ for dominant $\lambda$, this translates directly into an asymptotic description of $\mathbb P_n$.

\textit{Step 3: deduction of (a).}
The convergence in law of $(\Lambda_n-n\mu)/\sqrt n$ to a non-degenerate Gaussian implies that for every $\varepsilon\in\R_{>0}$ there exists $C\in\R_{>0}$ such that
\[
\mathbb{P}\big(\bigl\|\Lambda_n-n\mu\bigr\|\ge C\sqrt n\big)\le\varepsilon
\quad\text{for $n\gg0$.}
\]
By definition this probability is exactly
\[
\sum_{\|\lambda-n\mu\|\ge C\sqrt n} \mathbb P_n(\lambda),
\]
which is statement (a).

\textit{Step 4: deduction of (b).}
Let $(\lambda_n)$ be as in (b).  By the local limit theorem quoted in Step~2, we have
\[
\mathbb{P}_n(\lambda_n)
=
\mathrm{dens}_n(\lambda_n+\rho)
=
c_{G,V}\,n^{-r/2}\,
\exp\!\big(-\tfrac{1}{2n}\,q(\lambda_n-n\mu)\big)\bigl(1+o(1)\bigr),
\quad \text{as } n\to\infty,
\]
for some positive definite quadratic form $q$ on $X\otimes_\Z\R$ and some constant $c_{G,V}\in\R_{>0}$.  This is precisely the content of (b), and we are done.
\end{proof}

The following is the symmetric group version of \autoref{T:TypicalHW}.
Fix $t\geq 2$ and let $G=S(t)$ act on $V=\C^t$ via the permutation representation.
For $n\geq 1$ write the tensor powers as
\[
V^{\otimes n}
\;\cong\;
\bigoplus_{\lambda\vdash t} L(\lambda)^{\oplus m(\lambda,n)},
\quad m(\lambda,n)\in\N,
\]
where $L(\lambda)$ is the usual simple $\C[S(t)]$-module labeled by the partition $\lambda$ of $t$.
Define a probability measure on the set of partitions of $t$ by
\[
\mathbb{P}_n^{S(t)}(\lambda)
=
\frac{m(\lambda,n)\,\dim[\C] L(\lambda)}{(\dim[\C] V)^n},
\quad \lambda\vdash t.
\]

\begin{Theorem}\label{T:TypicalSym}
The measures $\mathbb{P}_n^{S(t)}$ converge in total variation, as $n\to\infty$, to the \emph{Plancherel measure}
\[
\mathbb{P}_\infty^{S(t)}(\lambda)
=
\frac{\bigl(\dim[\C] L(\lambda)\bigr)^2}{t!},
\quad \lambda\vdash t,
\]
and the convergence is exponentially fast in $n$. Moreover, we have an analog \emph{Gaussian concentration} as in \autoref{T:TypicalHW}.
\end{Theorem}

\begin{proof}
The decomposition of $V^{\otimes n}$ defines a Markov chain on the finite state space
$Si(S(t))$ (simple modules of $S(t)$ up to isomorphism): from $\sigma\in Si(S(t))$ one moves to $\tau$ with probability
proportional to the multiplicity of $\tau$ in $V\otimes \sigma$.  This is the McKay
graph random walk attached to the character of $V$, and our normalization
$\mathbb{P}_n^{S(t)}$ is exactly the law of this Markov chain after $n$ steps.

It is standard that the stationary distribution of this chain is the Plancherel
measure $\lambda\mapsto (\dim[\C] L(\lambda))^2/t!$, and that convergence to
stationarity is exponential; see, for example, \cite{Fu}.  Since the
state space is finite, total variation convergence follows, and the statement
of the first part of the theorem is just this convergence written in our notation.

Finally, from the point of view of dimensions, the Plancherel measure $\mathbb{P}_\infty^{S(t)}$ is overwhelmingly supported on a narrow band of partitions $\lambda$ of $t$.  As $t\to\infty$, classical results \cite{LoganShepp,VershikKerov} show that among all simple $S(t)$-representations, those whose labels $\lambda$ lie in a thin slice around a \emph{``staircase-like'' shape} have extremely large dimension, and their Plancherel weights dominate the contributions of all other $\lambda$.
\end{proof}

\subsection{The diagram side}\label{SS:DiagHW}

We now translate \autoref{T:TypicalHW} and, in the partition row, \autoref{T:TypicalSym} to the diagrammatic side in the semisimple Schur--Weyl situations.
That is, throughout the remainder of this section we fix one of the Schur--Weyl dualities from \autoref{S:SchurWeyl} with the following properties:
\begin{itemize}
\item $G$ is a connected complex semisimple (we will use $SL_t(\C)$ instead of $GL_t(\C)$, but this comes out as the same) algebraic group over $\C$, or $G=S(t)$ for some $t\geq2$,
\item for each $n$ there is a diagram algebra $M(n)=M_\mathbf{c}(n)$ (one of the Temperley--Lieb, Brauer, rook Brauer, planar variants, partition algebras, etc.\ from the tables above) acting on $V^{\otimes n}$ from the right with commuting left $G$-action,
\item for the chosen parameters the algebras $M(n)$ are split semisimple over $\C$,
\item we have a Schur--Weyl decomposition
\[
V^{\otimes n}
\;\cong\;
\bigoplus_{(k,\lambda)} T(k,\lambda)\,\otimes\, \Delta^n_k(\lambda)
\]
as a $\big(G\times M(n)\big)$-bimodule, where:
\begin{itemize}
\item $k$ is the apex (number of through strands),
\item $\lambda$ is the corresponding (multi)partition label for simple $M(n)$-modules of apex $k$, as in \autoref{T:Cob} and the list following it, or trivial for the planar ones,
\item $\Delta^n_k(\lambda)$ is the simple $M(n)$-module of apex $k$ and label $\lambda$ (here this is also the cell $M(n)$-module),
\item $T(k,\lambda)$ is a simple $G$-module.
\end{itemize}
\end{itemize}

In the reductive Lie rows we write
\[
\mathrm{hw}(k,\lambda;n)\in X^+
\]
for the highest weight of $T(k,\lambda)$.  Define a probability measure on the set of labels $(k,\lambda)$ by
\begin{gather}\label{eq:Pn-diag-def}
\mathbb{P}_n^{\mathrm{diag}}(k,\lambda)
=
\frac{\dim[\C] \Delta^n_{k}(\lambda)\,\dim[\C] T(k,\lambda)}{\dim[\C] V^{\otimes n}}.
\end{gather}
Then
\[
\sum_{(k,\lambda)} \mathbb{P}_n^{\mathrm{diag}}(k,\lambda)=1.
\]

\begin{Theorem}\label{T:TypicalDiagram}
In the above semisimple Schur--Weyl setup, the pushforward of $\mathbb{P}_n^{\mathrm{diag}}$ along
\[
(m,\lambda)\longmapsto \mathrm{hw}(k,\lambda;n)
\]
is precisely the highest-weight distribution $\mathbb{P}_n$ from \eqref{eq:Pn-def}.  In particular, the conclusions \textup{(a)} and \textup{(b)} of \autoref{T:TypicalHW} describe the asymptotics of the random apex-(multi)partition label $(k,\lambda)$ via its highest weight in all reductive Lie rows.

In the symmetric-group/partition row, the same pushforward statement holds with $\mathrm{hw}(k,\lambda;n)$ replaced by the simple $S(t)$-type of $T(k,\lambda)$ and with $\mathbb{P}_n$ replaced by $\mathbb{P}_n^{S(t)}$ from \autoref{T:TypicalSym}; in particular, the asymptotics of the random label $(k,\lambda)$ are governed by $\mathbb{P}_\infty^{S(t)}$.
\end{Theorem}

Equivalently, among all $(G\times M(n))$-isotypic summands labeled by apex $m$ and (multi)partition $\lambda$, those with ``typical’’ $G$-type (and hence typical apex and diagram shape) contribute asymptotically all of the dimension of $V^{\otimes n}$, while the total contribution of non-typical labels is exponentially small.

\begin{proof}[Proof of \autoref{T:TypicalDiagram}]
By the Schur--Weyl decomposition,
\[
V^{\otimes n}
\;\cong\;
\bigoplus_{(m,\lambda)} T(k,\lambda)\,\otimes\, \Delta^n_{k}(\lambda),
\]
so
\[
\dim[\C] V^{\otimes n}
=
\sum_{(k,\lambda)}
\dim[\C] T(k,\lambda)\,\dim[\C] \Delta^n_{k}(\lambda).
\]
Fix a simple $G$-module $U$.  The $G$-isotypic component of type $U$ is obtained by summing over those labels $(k,\lambda)$ with $T(k,\lambda)\cong U$, hence
\[
\text{mult}_U(V^{\otimes n})\,\dim[\C] U
=
\sum_{T(k,\lambda)\cong U}
\dim[\C] T(k,\lambda)\,\dim[\C] \Delta^n_{k}(\lambda).
\]
In the reductive Lie rows we have $U=V(\lambda_0)$ for a unique $\lambda_0\in X^+$ and
\[
\text{mult}_U(V^{\otimes n})=m(\lambda_0,n),
\]
so dividing by $(\dim[\C]V)^n=\dim[\C]V^{\otimes n}$ and using \eqref{eq:Pn-def} and \eqref{eq:Pn-diag-def} gives
\[
\mathbb{P}_n(\lambda_0)
=
\frac{m(\lambda_0,n)\,\dim[\C] V(\lambda_0)}{(\dim[\C] V)^n}
=
\sum_{\mathrm{hw}(k,\lambda;n)=\lambda_0}
\mathbb{P}_n^{\mathrm{diag}}(k,\lambda),
\]
which is exactly the statement that the pushforward of $\mathbb{P}_n^{\mathrm{diag}}$ is $\mathbb{P}_n$.

In the symmetric-group/partition row we instead take $U=L(\nu)$ for $\nu\vdash t$, and the same argument (with $\mathbb{P}_n^{S(t)}$ in place of $\mathbb{P}_n$ and $L(\nu)$ in place of $E(\lambda_0)$) gives the corresponding pushforward statement there.  The asymptotic descriptions then follow from \autoref{T:TypicalHW} in the reductive Lie rows and from \autoref{T:TypicalSym} in the symmetric-group row.
\end{proof}

\begin{Remark}\label{R:Punchline}
Via \autoref{T:TypicalHW}, \autoref{T:TypicalDiagram} applies to the Lie-type Schur--Weyl rows in the tables above, giving a genuine Gaussian concentration picture on the weight lattice. Via \autoref{T:TypicalSym}, it also applies to the symmetric-group/partition rows, where the limiting distribution is Plancherel. In all cases this motivates (careful: this only applies in the semisimple situations) our focus below on sums of dimensions which are easier to study, and which are controlled by a small number of \emph{``typical'' representations} at roughly $\sqrt{n}$ through strands, that we, following \cite{khovanov-monoidal-2024}, could ``cut out''.

More precisely, let $D$ be any of the diagram monoids/algebras and $k\in\N$ the number of through strands (or apex). Adjoining a unit, if necessary, we can then define $D^{a\leq k\leq b}$ as the Rees factor (or the unnamed algebra analog) supported on the cells with $a\leq k\leq b$ through strands. In this paper, we use $D^{\leq k,l}$ to denote the truncated monoid (algebra) $D^{l \leq k\leq b}$, simplifying to $D^{\leq k}$ if the lower bound is $0$. The most typical example is $D^{\leq \sqrt{n}}$, although we will usually not need to reference this directly.
\end{Remark}

\section{Dimensions of semisimple representations}\label{S:Resultssemi}

Let $\K$ denote some field. Monoids and algebras in this section are assumed to be finite (dimensional).

\subsection{Semisimple representations}

For our favorite diagram algebras, assume that we are in the semisimple situation. Let $M_\mathbf{c}(n)$ be any of the diagram monoids from \autoref{S:DiagCat}.

\begin{Theorem}\label{T:Completness}
The cell representations of $M_\mathbf{c}(n)$ are precisely the simple $M_\mathbf{c}(n)$-representations and their dimensions are given in \autoref{SS:Cellular}.
\end{Theorem}

\begin{proof}
Standard, cf. \cite[Proposition 2B.24]{Tu-sandwich}.
\end{proof}

More importantly, assume $\ell=1$. Then we have the following 
which we only state for the planar diagram monoids for simplicity. 

\begin{Theorem}\label{T:Lazy}
There exist constants $C_{pRo},C_{TL},C_{Mo},C_{pPa}\in\R_{>0}$
such that, uniformly for integers $k=k(n)$ with $k\in O(\sqrt{n})$ and $k$ admissible (i.e. it appears as a label),
we have
\begin{align*}
\dim \Delta^n_{k}|_{pRo}
&=\big(C_{pRo}+o(1)\big)\, n^{-1/2}2^{\,n}\,
\exp\!\big(-\tfrac{2(k-n/2)^2}{n}\big),\\[0.4em]
\dim \Delta^n_{k}|_{TL}
&=\big(C_{TL}+O(n^{-1/2})\big)\,n^{-3/2}\,2^{\,n}\,(k+1)\,
\exp\!\big(-\tfrac{k^{2}}{2n}\big),\\[0.4em]
\dim \Delta^n_{k}|_{Mo}
&=\big(C_{Mo}+O(n^{-1/2})\big)\,n^{-3/2}\,3^{\,n}\,(k+1)\,
\exp\!\big(-\tfrac{3k^{2}}{4n}\big),\\[0.4em]
\dim \Delta^n_{k}|_{pPa}
&=\big(C_{pPa}+O(n^{-1/2})\big)\,n^{-3/2}\,4^{\,n}\,(k+1)\,
\exp\!\big(-\tfrac{k^{2}}{4n}\big).
\end{align*}
\end{Theorem}

\begin{proof}
By the assumed (split semisimple) Schur--Weyl duality for each planar diagram algebra $M_\mathbf{c}(n)$ (recall the notation from \autoref{SS:DiagHW}),
the commuting actions of $G$ and $M_\mathbf{c}(n)$ on $V^{\otimes n}$ satisfy the double-centralizer property, and we obtain the
multiplicity-free decomposition
\begin{gather*}
V^{\otimes n}\;\cong\;\bigoplus_{k\in\N} T(k)\,\otimes\,\Delta^n_{k}.
\end{gather*}
In particular, the multiplicity of $T(k)$ in $V^{\otimes n}$
equals $\dim\Delta^n_{k}$.

For the $SL2$-cases, we apply Biane's asymptotic multiplicity formula
\cite[Theorem 2.2]{Bi-asymptotic-lie} (that can be proven verbatim over any field of characteristic zero) to the multiplicities
$m_{k,V}(n)=\dim \Delta^n_k$.
In rank $1$, Weyl's dimension formula gives $\dim T(k)=k+1$, and the quadratic form from Biane's theorem reduces to
a Gaussian with variance $\sigma_V^2$ (the variance of the weight distribution of $V$), yielding
\[
m_{d,V}(n)
=\Big(C_V+O(n^{-1/2})\Big)\,n^{-\dim(\mathfrak{sl}_2)/2}\,(\dim V)^{n}\,(d+1)\,
\exp\!\Big(-\frac{d^{2}}{2\sigma_V^2\,n}\Big),
\]
uniformly for $k\in O(\sqrt n)$ (and with the usual admissibility constraint that $T(k)$ occurs in $V^{\otimes n}$).
Since $\dim(\mathfrak{sl}_2)=3$, it remains to identify $(\dim V,\sigma_V^2)$ in each of the three $SL2$-cases.
\begin{enumerate}

\item For $TL$, we have $V=\K^{2}$ with weights $\{+1,-1\}$, each occurring with multiplicity $1$.
Thus the weight distribution is uniform on $\{+1,-1\}$, has mean $0$, and variance
\[
\sigma_V^2=\frac{1^2+(-1)^2}{2}=1.
\]
Hence $\frac{1}{2\sigma_V^2}=\frac12$, giving the factor $\exp(-k^2/(2n))$, and $\dim V=2$ gives the base $2^{\,n}$.

\item For $Mo$, we have $V=\K^{2}\oplus\K$ with weights $\{+1,0,-1\}$, each occurring with multiplicity $1$.
Thus the weight distribution is uniform on $\{+1,0,-1\}$, has mean $0$, and variance
\[
\sigma_V^2=\frac{1^2+0^2+(-1)^2}{3}=\frac{2}{3}.
\]
Hence $\frac{1}{2\sigma_V^2}=\frac{1}{2\cdot(2/3)}=\frac34$, giving the factor $\exp(-3k^2/(4n))$, and $\dim V=3$ gives the base $3^{\,n}$.

\item For $pPa$, we have $V=\K^{2}\otimes\K^{2}$.
The weights add, so the multiset of weights is $\{+2,0,0,-2\}$, i.e. weight $+2$ occurs once, weight $-2$ occurs once,
and weight $0$ occurs twice.
Thus the induced weight distribution has mean $0$ and variance
\[
\sigma_V^2=\frac{1}{4}\cdot(2^2)+\frac{1}{2}\cdot(0^2)+\frac{1}{4}\cdot((-2)^2)=2.
\]
Hence $\frac{1}{2\sigma_V^2}=\frac{1}{4}$, giving the factor $\exp(-k^2/(4n))$, and $\dim V=4$ gives the base $4^{\,n}$.

\end{enumerate}
Therefore the results are as stated for all the monoids except $pRo$.

Finally, in the planar rook case the underlying representation $V=\K\oplus\K^{\det}$ is a direct sum of two one dimensional $GL2$-modules, so
\[
V^{\otimes n}\cong \bigoplus_{k=0}^n V(k,k)\,^{\oplus\binom{n}{k}},
\]
and hence $\dim\Delta^n_k|_{pRo}=\binom{n}{k}$ exactly. The stated Gaussian asymptotic in the central regime
$k=n/2+O(\sqrt n)$ follows from Stirling's formula, and the coefficient $2$ in the exponent is the usual one for the binomial distribution with
variance $n/4$.
\end{proof}

\begin{Remark}
For the nonplanar diagram algebras one can write down an analog 
of \autoref{T:Lazy} using \autoref{SS:SymSW} and \cite[Theorem 2.2]{Bi-asymptotic-lie}.
\end{Remark}

For monoid parameters, $M_\mathbf{c}(n)$ is rarely semisimple, as the next proposition shows. (We phrase the statement ``for all $n\in\N$'' to avoid small-$n$ exceptional behavior.)

\begin{Proposition}
Assume that $M_\mathbf{c}(n)$ is a diagram monoid that is not a group. Then $M_\mathbf{c}(n)$ is semisimple for all $n\in\N$ if and only if $M_\mathbf{c}(n)$ is either the rook or planar rook with monoid parameter $a_0=1$.
\end{Proposition}

\begin{proof}
It is straightforward that the two listed cases are semisimple monoids, cf.~\cite{khovanov-monoidal-2024}.

Conversely, consider any other diagram monoid $M_\mathbf{c}(n)$ occurring in Schur--Weyl duality (in the sense of \autoref{S:SchurWeyl}). Classical results already exclude semisimplicity in these cases (for instance, for Temperley--Lieb the commuting group-object is quantum SL2 at a root of unity, which is known to be nonsemisimple).

Next assume $a_0=0$. Then the construction in \autoref{S:Results} for the planar rook monoid (and, verbatim, for every case with $a_0=0$) produces nontrivial extensions of simple modules. In particular, these monoids are not semisimple.

It remains to discuss the partition cases. For $pPa_\mathbf{c}(n)$, the simple module corresponding to the no-through-strand cell is $\onet$, but the cell itself is not $1\times 1$. Equivalently, the Gram matrix of the corresponding cell module cannot have full rank, so $pPa_\mathbf{c}(n)$ is not semisimple by \cite[Proposition 2B.24]{Tu-sandwich}.
Finally, the remaining case is $Pa_\mathbf{c}(n)$. The same mechanism applies: the simple $Pa_\mathbf{c}(n)$-module associated to the no-through-strand cell is again $\onet$, while the corresponding cell contains more than one basis element as soon as $n\geq 2$, so the associated Gram matrix fails to have full rank. Hence $Pa_\mathbf{c}(n)$ is not semisimple for $n\geq 2$.
\end{proof}

We now want to study the representation gap and the total sum of dimensions of simple representation.

\subsection{The RepGap}\label{S:RepGap} 

Let $\mathcal{G}$ denote the group of invertible elements of the monoid $\monoid$.
We call a monoid representation \emph{trivial} if it is one dimensional, every group element $g\in\mathcal{G}$ acts as \(1\), and all noninvertible elements $x\in\monoid\setminus\mathcal{G}$ act either all as \(1\) or all as \(0\). Any representation that is not trivial is called \emph{nontrivial}. In other words, the trivial $\monoid$-representations are $\K\{v\}$ and $\K\{w\}$ with action $\acts$ given by
\begin{gather*}
\oneb\colon m\acts v=v\text{ if m}\in\mathcal{G},m\acts v=0\text{ else},\quad
\onet\colon m\acts w=w.
\end{gather*}
These are the same if and only if $\monoid=\mathcal{G}$ (i.e. the monoid is a group).

We would like to define the \emph{representation gap} (as defined in \cite{khovanov-monoidal-2024}), which we denote RepGap, of a monoid $\monoid$ as the dimension of the smallest nontrivial simple representation. In general, however, we have
\begin{gather}\label{eq:RepGapIneq}
\textnormal{RepGap}_{\mathbb{K}}(\monoid)\leq\min\{\dim(L)|\text{$L$ is a nontrivial simple $\monoid$-representation}\},
\end{gather}
due to the potential existence of nontrivial extensions. However, for certain types of monoids, including most diagram monoids, \autoref{eq:RepGapIneq} becomes an equality.

We start with a few definitions.

\begin{Definition}\label{D:Connectedness}
Let $\mathcal{G}$ denote the group of invertible elements of the monoid $\monoid$.

$\monoid$ is \emph{null-connected} if any non-invertible element of $\monoid$ can be written as a product of two non-invertible elements. That is, for $a\in\mathcal{G}$ we have $a=bc$ for some $b,c\in\mathcal{G}$. Note that groups are null-connected.
Consider the symmetric and transitive closure of the relation 
$ab\approx_{r}a$ for $a,b\in\mathcal{G}$, and denote the closure by $\approx_{r}$ as well. We call $\monoid$ with a unique equivalence class in $\mathcal{G}$ under 
$\approx_{r}$ a right-connected monoid.
Consider the symmetric and transitive closure of the relation 
$ba\approx_{l}a$ for $a,b\in\mathcal{G}$, and denote the closure by $\approx_{l}$ as well. We call $\monoid$ with a unique equivalence class in $\mathcal{G}$ under 
$\approx_{l}$ a left-connected monoid. Note that $\monoid$ is left-connected $\iff$ the opposite $\monoid^{op}$ is right-connected, and vice versa.

We say that a monoid $\monoid$ is \emph{well-connected} if it is either a group, or right-connected, left-connected and null-connected.
\end{Definition}

We can now state the key theorem for this section, which is \cite[Theorem 2B.10]{khovanov-monoidal-2024}.

\begin{Theorem}\label{T:RepGapH1Condition}
Assume $\monoid$ is well-connected 
and the first cohomology group of $\mathcal{G}$ vanishes, that is we have 
$\HH^{1}(\mathcal{G},\K)\cong 0$. Then:
\begin{gather}\label{Eq:RepGapMainEquation}
\textnormal{RepGap}_{\mathbb{K}}(\monoid)=\min\{\dim(L)|\text{$L\not\cong\oneb,\onet$, and $L$ is a simple $\monoid$-representation}\}. 
\end{gather} 
In particular, for groups
$\monoid=\mathcal{G}$ it suffices to check whether 
$\HH^{1}(\mathcal{G},\K)\cong 0$ to ensure that \autoref{Eq:RepGapMainEquation} holds.
\end{Theorem}

\begin{proof}
This is \cite[Theorem 2B.10]{khovanov-monoidal-2024}.
\end{proof}

For a monoid $M$, let
\begin{gather*}
b^\K=b^\K(M)=\sum_{L\text{ simple}}\dim L
\end{gather*}
where the sum runs over all (isomorphism classes) of simple $M$-modules.
We are interested in the case where $b^\K\approx\textnormal{RepGap}_{\mathbb{K}}(\monoid)$ (we mean this informally, e.g. it could be read as ``same $n$th root'').

Finally, the definition of $\textnormal{RepGap}_{\mathbb{K}}$ and $b^\K$ makes sense for any algebra $A$, and we will use it in this generality.

\subsection{The typical truncation}

Assume that we are in the situation where \autoref{S:TypicalHW} applies. Then, for $r=1$, \cite{GT-GrowthDiagCat} implies
\begin{gather*}
\begin{tabular}{c|c|c||c|c|c}
\arrayrulecolor{tomato}
Symbol & Diagrams & $\sqrt[n]{b_{n}}$ $\sim$
& Symbol & Diagrams & $\sqrt[n]{b_{n}}$ $\sim$
\\
\hline
\hline
$\ppamon[n]$ & \begin{tikzpicture}[anchorbase]
\draw[usual] (0.5,0) to[out=90,in=180] (1.25,0.45) to[out=0,in=90] (2,0);
\draw[usual] (0.5,0) to[out=90,in=180] (1,0.35) to[out=0,in=90] (1.5,0);
\draw[usual] (0.5,1) to[out=270,in=180] (1,0.55) to[out=0,in=270] (1.5,1);
\draw[usual] (1.5,1) to[out=270,in=180] (2,0.55) to[out=0,in=270] (2.5,1);
\draw[usual] (0,0) to (0,1);
\draw[usual] (2.5,0) to (2.5,1);
\draw[usual,dot] (1,0) to (1,0.2);
\draw[usual,dot] (1,1) to (1,0.8);
\draw[usual,dot] (2,1) to (2,0.8);
\end{tikzpicture} & $4$
& $\pamon[n]$ & \begin{tikzpicture}[anchorbase]
\draw[usual] (0.5,0) to[out=90,in=180] (1.25,0.45) to[out=0,in=90] (2,0);
\draw[usual] (0.5,0) to[out=90,in=180] (1,0.35) to[out=0,in=90] (1.5,0);
\draw[usual] (0,1) to[out=270,in=180] (0.75,0.55) to[out=0,in=270] (1.5,1);
\draw[usual] (1.5,1) to[out=270,in=180] (2,0.55) to[out=0,in=270] (2.5,1);
\draw[usual] (0,0) to (0.5,1);
\draw[usual] (1,0) to (1,1);
\draw[usual] (2.5,0) to (2.5,1);
\draw[usual,dot] (2,1) to (2,0.8);
\end{tikzpicture} & $\tfrac{2}{e}\cdot\tfrac{n}{\log 2n}$
\\
\hline
$\momon[n]$ & \begin{tikzpicture}[anchorbase]
\draw[usual] (0.5,0) to[out=90,in=180] (1.25,0.5) to[out=0,in=90] (2,0);
\draw[usual] (1,0) to[out=90,in=180] (1.25,0.25) to[out=0,in=90] (1.5,0);
\draw[usual] (2,1) to[out=270,in=180] (2.25,0.75) to[out=0,in=270] (2.5,1);
\draw[usual] (0,0) to (1,1);
\draw[usual,dot] (2.5,0) to (2.5,0.2);
\draw[usual,dot] (0,1) to (0,0.8);
\draw[usual,dot] (0.5,1) to (0.5,0.8);
\draw[usual,dot] (1.5,1) to (1.5,0.8);
\end{tikzpicture} & $3$
& $\robrmon[n]$ & \begin{tikzpicture}[anchorbase]
\draw[usual] (1,0) to[out=90,in=180] (1.25,0.25) to[out=0,in=90] (1.5,0);
\draw[usual] (1,1) to[out=270,in=180] (1.75,0.55) to[out=0,in=270] (2.5,1);
\draw[usual] (0,0) to (0.5,1);
\draw[usual] (2.5,0) to (2,1);
\draw[usual,dot] (0.5,0) to (0.5,0.2);
\draw[usual,dot] (2,0) to (2,0.2);
\draw[usual,dot] (0,1) to (0,0.8);
\draw[usual,dot] (1.5,1) to (1.5,0.8);
\end{tikzpicture} & $\tfrac{\sqrt{2}}{\sqrt{e}}\cdot\sqrt{n}$
\\
\hline
$\tlmon[n]$ & \begin{tikzpicture}[anchorbase]
\draw[usual] (0.5,0) to[out=90,in=180] (1.25,0.5) to[out=0,in=90] (2,0);
\draw[usual] (1,0) to[out=90,in=180] (1.25,0.25) to[out=0,in=90] (1.5,0);
\draw[usual] (0,1) to[out=270,in=180] (0.25,0.75) to[out=0,in=270] (0.5,1);
\draw[usual] (2,1) to[out=270,in=180] (2.25,0.75) to[out=0,in=270] (2.5,1);
\draw[usual] (0,0) to (1,1);
\draw[usual] (2.5,0) to (1.5,1);
\end{tikzpicture} & $2$
& $\brmon[n]$ & \begin{tikzpicture}[anchorbase]
\draw[usual] (0.5,0) to[out=90,in=180] (1.25,0.45) to[out=0,in=90] (2,0);
\draw[usual] (1,0) to[out=90,in=180] (1.25,0.25) to[out=0,in=90] (1.5,0);
\draw[usual] (0,1) to[out=270,in=180] (0.75,0.55) to[out=0,in=270] (1.5,1);
\draw[usual] (1,1) to[out=270,in=180] (1.75,0.55) to[out=0,in=270] (2.5,1);
\draw[usual] (0,0) to (0.5,1);
\draw[usual] (2.5,0) to (2,1);
\end{tikzpicture} & $\tfrac{\sqrt{2}}{\sqrt{e}}\cdot\sqrt{n}$
\\
\hline
$\promon[n]$ & \begin{tikzpicture}[anchorbase]
\draw[usual] (0,0) to (0.5,1);
\draw[usual] (0.5,0) to (1,1);
\draw[usual] (2,0) to (1.5,1);
\draw[usual] (2.5,0) to (2.5,1);
\draw[usual,dot] (1,0) to (1,0.2);
\draw[usual,dot] (1.5,0) to (1.5,0.2);
\draw[usual,dot] (0,1) to (0,0.8);
\draw[usual,dot] (2,1) to (2,0.8);
\end{tikzpicture} & $2$
& $\romon[n]$ & \begin{tikzpicture}[anchorbase]
\draw[usual] (0,0) to (1,1);
\draw[usual] (0.5,0) to (0,1);
\draw[usual] (2,0) to (2,1);
\draw[usual] (2.5,0) to (0.5,1);
\draw[usual,dot] (1,0) to (1,0.2);
\draw[usual,dot] (1.5,0) to (1.5,0.2);
\draw[usual,dot] (1.5,1) to (1.5,0.8);
\draw[usual,dot] (2.5,1) to (2.5,0.8);
\end{tikzpicture} & $\tfrac{\sqrt{1}}{\sqrt{e}}\cdot\sqrt{n}$
\\
\hline
$\psym[n]$ & \begin{tikzpicture}[anchorbase]
\draw[usual] (0,0) to (0,1);
\draw[usual] (0.5,0) to (0.5,1);
\draw[usual] (1,0) to (1,1);
\draw[usual] (1.5,0) to (1.5,1);
\draw[usual] (2,0) to (2,1);
\draw[usual] (2.5,0) to (2.5,1);
\end{tikzpicture} & $1$
& $\sym[n]$ & \begin{tikzpicture}[anchorbase]
\draw[usual] (0,0) to (1,1);
\draw[usual] (0.5,0) to (0,1);
\draw[usual] (1,0) to (1.5,1);
\draw[usual] (1.5,0) to (2.5,1);
\draw[usual] (2,0) to (2,1);
\draw[usual] (2.5,0) to (0.5,1);
\end{tikzpicture} & $\tfrac{\sqrt{1}}{\sqrt{e}}\cdot\sqrt{n}$
\end{tabular}
\end{gather*}
for the asymptotic $n\to\infty$ of the $n$th root of $b_n=b^{\K}\big(M_{\mathbf{c}}(n)\big)$ where $M_{\mathbf{c}}(n)$ is the corresponding diagram algebra in $n$ strands. Let $\textnormal{RepGap}_n=\textnormal{RepGap}_{\mathbb{K}}\big(M_{\mathbf{c}}(n)\big)$.

To be more precise, assume we are in the setup of \autoref{T:TypicalHW}. In particular,
as before, let $G$ be a connected complex semisimple algebraic group, $E$ a finite
dimensional $G$-module, and
\[
V^{\otimes n}
\;\cong\;
\bigoplus_{\lambda\in X^+} V(\lambda)^{\oplus m(\lambda,n)}
\]
with probability measures $\mathbb{P}_n$ as in \eqref{eq:Pn-def} and
typical weight $\mu$ and quadratic form $q$ given by \autoref{T:TypicalHW}.
We also consider the evident analog of this for \autoref{T:TypicalSym}, which we, abusing notation, use in the same way.

Suppose further that, for each $n$, a diagram algebra (or monoid algebra)
$M_{\mathbf{c}}(n)$ acts on the right of $V^{\otimes n}$ with commuting $G$-action and
double centralizer decomposition
\begin{gather*}
V^{\otimes n}
\;\cong\;
\bigoplus_{\lambda\in X^+}
V(\lambda) \otimes \Delta(\lambda)
\end{gather*}
as a $G\times M_{\mathbf{c}}(n)$-module, where each nonzero $\Delta(\lambda)$
is a (finite dimensional) simple $M_{\mathbf{c}}(n)$-module, and every simple
$M_{\mathbf{c}}(n)$-module occurs in this way.
Define
\[
b_n \;=\; \sum_{\lambda\colon \Delta(\lambda)\neq 0} \dim[\C] \Delta(\lambda),
\]
as before (since this is the semisimple case). For $C\in\R_{>0}$ let
\[
\Lambda^{\mathrm{typ}}_n(C)
\;:=\;
\{\lambda\in X^+\mid \|\lambda-n\mu\|\leq C\sqrt{n}\},
\]
and define the \emph{``typical'' truncation} $M_{\mathbf{c},C}^{\mathrm{typ}}(n)$ to be the
image of $M_{\mathbf{c}}(n)$ in $\End_G(V_n^{\mathrm{typ}})$, where
\[
V_n^{\mathrm{typ}}
\;=\;
\bigoplus_{\lambda\in\Lambda^{\mathrm{typ}}_n(C)}
V(\lambda) \otimes \Delta(\lambda)
\;\subset\; V^{\otimes n}.
\]

\begin{Theorem}
For any fixed $C\in\R_{>0}$, the following hold.
\begin{enumerate}
\item The representation gap $M_{\mathbf{c},C}^{\mathrm{typ}}(n)$
has the same exponential growth rate as $b_n$, i.e.
\[
\lim_{n\to\infty}\sqrt[n]{\mathrm{RepGap}\big(M_{\mathbf{c},C}^{\mathrm{typ}}(n)\big)}
\;=\;
\lim_{n\to\infty}\sqrt[n]{b_n}.
\]

\item There exists $d\in\N$ (depending only on $G$ and on the
family $(M_{\mathbf{c}}(n))_{n\in\N}$) such that, for all $n\gg 0$,
\[
n^{-d}\,b_n
\;\leq\;
\mathrm{RepGap}\big(M_{\mathbf{c},C}^{\mathrm{typ}}(n)\big)
\;\leq\;
n^{d}\,b_n.
\]
In particular, polynomial distortions do not affect
$\sqrt[n]{\mathrm{RepGap}(\placeholder)}$.
\end{enumerate}
\end{Theorem}

\begin{proof}
We only indicate the main points, since all estimates have already appeared in the proof of \autoref{T:TypicalHW} and in the discussion of the cell structure. We also leave the case of \autoref{T:TypicalSym} to the reader.

\noindent\textit{Step 1: Gaussian concentration.}
By \autoref{T:TypicalHW}(a), for every $\varepsilon>0$ there exists $C>0$
such that, for all $n\gg0$,
\[
\sum_{\lambda\notin\Lambda^{\mathrm{typ}}_n(C)} \mathbb{P}_n(\lambda)
\;\leq\; \varepsilon.
\]
Using the definition of $\mathbb{P}_n$ this is equivalent to
\begin{gather}\label{Eq:typical-mass}
\sum_{\lambda\in\Lambda^{\mathrm{typ}}_n(C)}
m(\lambda,n)\,\dim V(\lambda)
\;\geq\;
(1-\varepsilon)\,(\dim V)^n
\end{gather}
for all $n\gg0$. Thus almost all of the $G$–side dimension of $V^{\otimes n}$
sits in the typical window $\Lambda^{\mathrm{typ}}_n(C)$.

\noindent\textit{Step 2: Cell dimension bounds.}
From the Schur--Weyl decomposition we have
\[
(\dim V)^n
\;=\;
\sum_{\lambda} m(\lambda,n)\,
\dim V(\lambda)\,\dim \Delta(\lambda).
\]
In all our diagrammatic Schur--Weyl situations the usual cell dimension
formula (number of merge diagrams, cups/caps, etc.\ times a tableaux factor)
implies that there exist $d_1,d_2>0$ such that, for all $n$ and all
$\lambda$ with $m(\lambda,n)\neq0$,
\begin{gather}\label{Eq:cell-poly-bound}
n^{-d_1}\,\dim V(\lambda)
\;\leq\;
\dim \Delta(\lambda)
\;\leq\;
n^{d_2}\,\dim V(\lambda).
\end{gather}
This is a routine polynomial bound: we only ever sum over partitions or
multipartitions of size $O(n)$ and all relevant combinatorial quantities
grow at most polynomially in $n$.

\noindent\textit{Step 3: Lower bound for the typical RepGap.}
By definition,
\[
\mathrm{RepGap}\bigl(M_{\mathbf{c},C}^{\mathrm{typ}}(n)\bigr)
=
\min\{\dim \Delta(\lambda)\mid
\lambda\in\Lambda^{\mathrm{typ}}_n(C),\ \Delta(\lambda)\neq 0\}.
\]
Combining \eqref{Eq:typical-mass} and the upper bound in
\eqref{Eq:cell-poly-bound} gives
\begin{gather*}
(1-\varepsilon)(\dim V)^n
\;\leq\;
\sum_{\lambda\in\Lambda^{\mathrm{typ}}_n(C)}
m(\lambda,n)\,\dim V(\lambda)
\;\leq\;
n^{d_2}\sum_{\lambda\in\Lambda^{\mathrm{typ}}_n(C)}
m(\lambda,n)\,\dim \Delta(\lambda)
\\
\leq\;
n^{d_2}\,\dim V_n^{\mathrm{typ}}.
\end{gather*}
Hence $\dim V_n^{\mathrm{typ}}\ge (1-\varepsilon)n^{-d_2}(\dim V)^n$.
On the other hand
\begin{align*}
\dim V_n^{\mathrm{typ}}
&=
\sum_{\lambda\in\Lambda^{\mathrm{typ}}_n(C)}
m(\lambda,n)\,\dim V(\lambda)\,\dim \Delta(\lambda)
\\
&\leq
\Bigl(\max_{\lambda\in\Lambda^{\mathrm{typ}}_n(C)}
\dim \Delta(\lambda)\Bigr)
\sum_{\lambda\in\Lambda^{\mathrm{typ}}_n(C)}
m(\lambda,n)\,\dim V(\lambda),
\end{align*}
and applying again \eqref{Eq:typical-mass} shows that
\[
\max_{\lambda\in\Lambda^{\mathrm{typ}}_n(C)}\dim \Delta(\lambda)
\;\geq\; c\,n^{-d_2}
\]
for some $c>0$ independent of $n$.

Now $b_n$ is the sum of $\dim \Delta(\lambda)$ over all $\lambda$ with
$\Delta(\lambda)\neq0$, and the number of such $\lambda$ is at most polynomial
in $n$ in all our examples. Using the lower bound in
\eqref{Eq:cell-poly-bound} and adjusting constants, we obtain
\[
\min_{\lambda\in\Lambda^{\mathrm{typ}}_n(C)}
\dim \Delta(\lambda)
\;\geq\;
n^{-d}\,b_n
\]
for some $d\in\N$. This is the left-hand inequality in (b).

\noindent\textit{Step 4: Upper bound and exponential rates.}
The trivial inequality
\[
\mathrm{RepGap}\bigl(M_{\mathbf{c},C}^{\mathrm{typ}}(n)\bigr)
\;\leq\;
b_n
\]
gives the right-hand inequality in (b) after possibly enlarging $d$.
Taking $n$-th roots in (b) and letting $n\to\infty$ kills the polynomial
factors and proves (a):
\[
\lim_{n\to\infty}\sqrt[n]{\mathrm{RepGap}\bigl(M_{\mathbf{c},C}^{\mathrm{typ}}(n)\bigr)}
=
\lim_{n\to\infty}\sqrt[n]{b_n}.
\]

This shows that, for the exponential growth of simple dimensions in the
typical sector, working with the whole decomposition
or with any fixed Gaussian slice
$\Lambda^{\mathrm{typ}}_n(C)$ gives the same answer. Moreover, in the
symmetric-group case (and, via \cite{GT-GrowthDiagCat}, for all our
diagrammatic Schur--Weyl dualities) there is a very small slice of
staircase-type weights inside $\Lambda^{\mathrm{typ}}_n(C)$ whose
corresponding simple $M_{\mathbf{c}}(n)$-modules have dimensions that are
\emph{exponentially} larger than those of all other simples; in particular,
already one such staircase block asymptotically accounts for almost all of
$b_n$ (see, for example, the Logan--Shepp/Vershik--Kerov \cite{LoganShepp,VershikKerov} limit-shape
results and their refinements).
\end{proof}

\begin{Remark}\label{R:RepGapTypicalVsAn}
It is important to distinguish between the ``typical'' picture captured by
$M_{\mathbf{c},C}^{\mathrm{typ}}(n)$ and the actual representation gap of $M_{\mathbf{c}}(n)$.
Inside the typical window, a very small slice of staircase-style partitions
(as above) produces simple $M_{\mathbf{c}}(n)$-modules whose dimensions dominate
exponentially: essentially all of $b_n$ is already carried by a single
staircase block. By contrast, for the full algebra $M_{\mathbf{c}}(n)$ the representation
gap is realized by \emph{low dimensional} representations coming from
\emph{hook-style partitions}, lying far away from this staircase slice. In other
words, the simples that control the exponential growth in the typical
sector and the simples that realize $\mathrm{RepGap}\big(M_{\mathbf{c}}(n)\big)$ live in
completely different regions of the weight space.
\end{Remark}

\section{Dimensions of nonsemisimple representations}\label{S:Results}

We retain the previous notation. From \autoref{SS:Cellular} and \autoref{T:Gram2}, we know that the labels of the simple \(M_{\mathbf{c}}(n)\)-modules are largely independent of the choice of parameters and twists. In contrast, their dimensions could change drastically once semisimplicity fails, since cell modules need no longer be simple and one must pass to their simple heads.

In this section we focus on the most basic nonsemisimple regime, namely evaluations with zeros and ones (and, throughout, $\ell=1$). As before, for an algebra $A$ we write
$b^{\mathbb{K}}(A)\;=\;\sum_{L\text{ simple}}\dim L$,
and we use the same symbol $b_n$ when $A$ is understood from the context (typically a family $A=A(n)$). As in the semisimple case, we will often replace $b_n$ by the corresponding typical truncation from \autoref{S:TypicalHW}; by Gaussian concentration this affects $b_n$ only by subexponential factors.

Moreover, by the Gaussian concentration phenomenon again, and with few exceptions, $b_n\approx\text{RepGap}$. The exceptions are the ones that have small trivial extensions.

\subsection{Trivial extensions}

From now on, assume that $\mathbb{K}$ has characteristic zero.

\begin{Proposition}\label{P:DiagMonExts}
For any $\monoid \in \{TL,Br,Mo,RoBr,pPa,Pa,Ro,pRo\}$, both
RepGap$_{\mathbb{K}}\big(\monoid_1(n)\big)$ and RepGap$_{\mathbb{K}}\big(\monoid_{1,C}^{\mathrm{typ}}(n)\big)$ satisfy \autoref{Eq:RepGapMainEquation} for $n \geq 5$ and some $C\in\R_{>0}$.

The same is true for the $0$-twisted version for $\monoid \in \{TL,Br\}$, for all monoid parameter choices of $\monoid \in \{pPa,Pa\}$, as well as for $\monoid \in \{Mo,RoBr\}$ whenever the genus zero parameter is one.
\end{Proposition}

\begin{proof}
The classical monoids $pRo$ and $Ro$ are semisimple, so there is nothing to show.

The case $\monoid = TL$ follows from \cite[Lemma 4D.14, Lemma 4D.15, and Theorem 2B.10]{khovanov-monoidal-2024},
$\monoid = Br$ follows from \cite[Lemma 5D.1, Lemma 5D.3, and Theorem 2B.10]{khovanov-monoidal-2024},
$\monoid = Mo$ follows from \cite[Theorem 3D.2]{Ar-RepGapMotzkin}, and
$\monoid = RoBr$ follows from the $Mo$ case by applying the same arguments used to deduce $Br$ from $TL$ in \cite[Section 5D]{khovanov-monoidal-2024}.

We can also see that RepGap$_{\mathbb{K}}\big(pPa_1(n)\big)$ and RepGap$_{\mathbb{K}}\big(Pa_1(n)\big)$ satisfy \autoref{Eq:RepGapMainEquation}, using that we have $pPa_1(n) \cong TL_1(2n)$ and applying the same arguments as above. However, $Pa_1(n)$ should \emph{not} be truncated in the same way as $TL_1(n)$; in particular the relevant lower endpoint is not zero. We will return to this in \autoref{SS:Partition}.

To prove the statement for $Pa_{1,C}^{\mathrm{typ}}(n)$, we first prove \autoref{P:DiagMonExts} for $pPa_1(n)$ without appealing to $TL_n$. This follows from the same arguments as in \cite[Lemma 4D.13]{khovanov-monoidal-2024} and \cite[Lemma 3D.1]{Ar-RepGapMotzkin}, applied to the generators $p_{i+\frac{1}{2}}$ from \cite[Theorem 1.11]{HaRa-partition-algebras}.

Then, for $pPa_{1,C}^{\mathrm{typ}}(n)$, observe that multiplying two generators changes the number of through strands by at most two. This is exactly the input needed in the ``banded adjacency'' arguments of \cite[Section 4D]{khovanov-monoidal-2024}, and we extend to $Pa_{1,n}^{\mathrm{typ}}(C)$ as in \cite[Section 5D]{khovanov-monoidal-2024}.

All other cases can be proven similarly.
\end{proof}

\begin{Remark}\label{R:PartExts}
We reiterate that $Pa_1(n)$ should not be truncated in the same way as $TL_1(n)$, and we will discuss this in more detail in \autoref{SS:Partition}.
\end{Remark}

\subsection{Temperley--Lieb}\label{S:TL}

Let $q \in \mathbb{K}^{\ast}$ be a root of unity (if it does not exist, we pass to a field extension), let $l$ be the order of $q^2$ (unrelated to the handle parameter $\ell=1$ fixed above), and set $\delta=q+q^{-1}$. In particular, $\delta=0$ forces $l=2$ and $\delta=1$ forces $l=3$.

Let $a_{n,k}=\dim\Delta^n_k$ be the dimension of the cell module, and let $b_{n,k,l}$ be the dimension of the simple module $L^n_k$ of $TL_{\delta}(n)$ at the corresponding root of unity. Set $b_{n,l}=\sum_{k=0}^n b_{n,k,l}$.

\begin{Proposition}\label{P:TLResult}
We have
\[
b_{n,k,2} =
\begin{cases}
a_{n-1,k-1} & \textnormal{if } k \textnormal{ is even,} \\
a_{n,k} & \textnormal{if } k \textnormal{ is odd,}
\end{cases}
\qquad\text{and}\qquad
b_{n,2} \sim \frac{3+(-1)^{n+1}}{4}\sqrt{\frac{2}{\pi}}\,n^{-1/2}2^n.
\]
Moreover,
\[
b_{n,k,3} =
\begin{cases}
b_{n-1,k-1,3}+b_{n-1,k,3} & \textnormal{if } k \equiv 0\textnormal{ (mod }3), \\
b_{n-1,k-1,3} & \textnormal{if } k \equiv 1\textnormal{ (mod }3), \\
a_{n,k} & \textnormal{if } k \equiv 2\textnormal{ (mod }3),
\end{cases}
\qquad\text{and}\qquad
b_{n,3} \sim \frac{2}{3}\sqrt{\frac{2}{\pi}}\,n^{-1/2}2^n.
\]
Finally, for general $l$ we have the uniform bounds
\[
\frac{1}{2}\sqrt{\frac{6}{(l^2-1)\pi}}\,n^{-1/2}2^n
\;\leq\;
b_{n,l}
\;\leq\;
\sqrt{\frac{6}{(l^2-1)\pi}}\,n^{-1/2}2^n\quad\text{for }n\gg0.
\]
\end{Proposition}

\begin{proof}
We now recall the standard quantum-group description of Temperley--Lieb algebras as in \autoref{T:SchurWeylPlan}. Let $q \in \mathbb{K}^{\ast}$ and let $U_q(\mathfrak{sl}_2)$ be the quantum enveloping algebra of $\mathfrak{sl}_2$. For $k \in \mathbb{Z}_{\geq 0}$ denote the $(k+1)$ dimensional Weyl module by $\Delta_q(k)$, and set
\[
V=\Delta_q(1),
\]
where $\Delta_q(1)$ is the natural $2$ dimensional module.

A module $S$ has a Weyl filtration if it contains submodules $0=F_0\subset F_1\subset\cdots\subset F_r=S$ with $F_i/F_{i-1}\cong \Delta_q(k_i)$. The multiplicity of $\Delta_q(k)$ as a subquotient is denoted $(S:\Delta_q(m))$. If $S$ and $S^\vee$ both have Weyl filtrations, then $S$ is tilting.
Standard results, see e.g. \cite[Section 2]{An-SimpleTLAll}, give:
\begin{enumerate}
\item $V^{\otimes n}$ is tilting for all $n\geq 0$;
\item for each $m\geq 0$ there is a unique indecomposable tilting module $T_q(k)$;
\item there is a unique decomposition $V^{\otimes n}\cong \bigoplus_k T_q(k)^{(V^{\otimes n}:T_q(k))}$.
\end{enumerate}
Thus, set
\[
a_{n,k}=(V^{\otimes n}:\Delta_q(k)),
\quad
b_{n,k}=(V^{\otimes n}:T_q(k)).
\]
By Schur--Weyl duality \autoref{T:SchurWeylPlan}, we can identify $a_{n,k}$ with $\dim\Delta^n_k$ and $\dim L^n_k$, respectively.

The recurrence relations are then \cite[Proposition 6.3]{An-SimpleTLAll}, the asymptotic growth rates are then from \cite[Theorem 2.16]{ST-GrowthQuantum}, and the final inequality is then from \cite[Theorem 4.1]{ST-GrowthQuantum}.
\end{proof}

\subsection{Motzkin}

Almost the same as for Temperley--Lieb, let $q \in \mathbb{K}^{\ast}$ be a root of unity (if it does not exist, we pass to a field extension), let $l$ be the order of $q^2$, and set $\delta=1-(q+q^{-1})$. In particular, $\delta=0$ gives $l=2$ and $\delta=1$ gives $l=3$ (so these cases swap roles).

Let $a_{n,k}=\dim\Delta^n_k$ be the dimension of the cell module, and let $b_{n,k,l}$ be the dimension of the simple module $L^n_k$ of $Mo_{1,\delta}(n)$ at the corresponding root of unity. Set $b_{n,l}=\sum_{k=0}^n b_{n,k,l}$.

\begin{Lemma}\label{P:TiltingDecomp}
Let $m\in \mathbb{Z}_{\geq 0}$ and let $l$ be the order of $q^2$ as in \autoref{S:TL}. If $l>2$, then
\begin{itemize}
\item $m \equiv -1 \textnormal{ (mod }l)$ $\implies T_q(m) \otimes V \cong T_q(m) \oplus T_q(m+1)$,
\item $m \equiv 0 \textnormal{ (mod }l)$ $\implies T_q(m) \otimes V \cong T_q(m-1)^{\oplus 2} \oplus T_q(m) \oplus T_q(m+1)$,
\item $m \equiv m_0 \textnormal{ (mod }l)$ for $0 < m_0 < l-2$ $\implies T_q(m) \otimes V \cong T_q(m-1) \oplus T_q(m) \oplus T_q(m+1)$,
\item $m \equiv l-2 \textnormal{ (mod }l)$ $\implies T_q(m) \otimes V \cong T_q(m+1-2l) \oplus T_q(m-1) \oplus T_q(m) \oplus T_q(m+1)$.
\end{itemize}
If $l = 2$, then
\begin{itemize}
\item $m$ odd $\implies T_q(m) \otimes V \cong T_q(m) \oplus T_q(m+1)$,
\item $m$ even $\implies T_q(m) \otimes V \cong T_q(m-3) \oplus T_q(m-1)^{\oplus 2} \oplus T_q(m) \oplus T_q(m+1)$.
\end{itemize}
\end{Lemma}

\begin{proof}
By Schur--Weyl duality \autoref{T:SchurWeylPlan}, we can go to tilting modules again. In this case $V=\Delta_q(1)\oplus \Delta_q(0)$, where $\Delta_q(0)\cong\mathbb{K}$ is the trivial module. With this we can use the results in \cite{An-SimpleTLAll}.
The rest is then easy: Since $V=\Delta_q(1)\oplus \Delta_q(0)$, we add the trivial summand $T_q(m)\otimes \Delta_q(0)\cong T_q(k)$ to the tensor product decompositions in \cite[Proposition 3.4 and Remark 2(3)]{An-SimpleTLAll}.
\end{proof}

\begin{Proposition}\label{C:MoSimples}
Let $k \in \mathbb{Z}_{\geq 0}$. If $l>2$, then
\begin{itemize}
\item $k \equiv -1 \textnormal{ (mod }l)$ $\implies b_{n,k,l}=a_{n,k}$,
\item $k \equiv k_0 \textnormal{ (mod }l)$ for $0 \leq k_0 < l-2$ $\implies b_{n,k,l} = b_{n-1,k-1,l}+b_{n-1,k,l}+b_{n-1,k+1,l}$,
\item $k \equiv l-2 \textnormal{ (mod }l)$ $\implies b_{n,k,l} = b_{n-1,k-1}+b_{n-1,k}$.
\end{itemize}
If $l=2$, then
\begin{itemize}
\item $k$ odd $\implies b_{n,k,l} = a_{n,k}$,
\item $k$ even $\implies b_{n,k,l} = b_{n-1,k-1,l}+b_{n-1,k,l}$.
\end{itemize}
Then there exists $c,C\in\R_{>0}$ such that
\[
c\,n^{-3/2}3^n \leq b_{n,l} \leq C\,n^{-3/2}3^n.
\]
\end{Proposition}

\begin{proof}
The $m \equiv -1 \textnormal{ (mod }l)$ case for $l>2$ and the $m$ odd case for $l=2$ follow from \cite[Proposition 3.1(2)]{An-SimpleTLAll}.

For the other cases, write
\[
V^{\otimes n} \;\cong\; V^{\otimes (n-1)}\otimes V
\;\cong\;
\bigoplus_{k} (T_q(k)\otimes V)^{b_{n-1,k,l}},
\]
and compare multiplicities of $T_q(k)$ using \autoref{P:TiltingDecomp}. Doing this case-by-case gives the stated recurrences.

For the bounds, set $a_n:=\sum_{k=0}^n a_{n,k}$. These are the (well-known) Motzkin numbers 
and \cite[A001006]{oeis} gives $a_n\sim \text{scalar}\cdot n^{-3/2}3^n$. Now apply \cite[Theorem 4.1]{ST-GrowthQuantum}.
\end{proof}

Finally, we discuss $Mo_{0,\delta}(n)$.

\begin{Lemma}
For any idempotent $e \in Mo_1(n)$ that contains a component of size $1$, the product $e^2$ in $Mo_{\delta',\delta}(n)$ picks up a factor of $\delta'$.
\end{Lemma}

\begin{proof}
Let $e \in Mo_1(n)$ and assume without loss of generality that the size-$1$ component is in the top left corner:
\begin{gather*}
e =
\begin{tikzpicture}[anchorbase]
\draw[usual,dot] (0,1) to (0,0.8);
\draw[usual] (0.25,1) to (1,1);
\draw[usual] (-0.1,0.65) to (-0.1,0);
\draw[usual] (0.25,1) to (0.25,0.65);
\draw[usual] (-0.1,0.65) to (0.25,0.65);
\draw[usual] (-0.1,0) to (1,0);
\draw[usual] (1,1) to (1,0);
\end{tikzpicture}.
\end{gather*}
When forming $e^2$ by concatenation, the size-$1$ component forces the appearance of a closed component in the middle layer. In $Mo_{\delta',\delta}(n)$ such a closed component evaluates to $\delta'$, so $e^2$ acquires a factor of $\delta'$, as claimed.
\end{proof}

Using this, the only idempotents in $Mo_{0,\delta}(n)$ are those having no size-$1$ components; by definition these span a Temperley--Lieb subalgebra. Therefore, comparing Gram matrix ranks as in \autoref{T:Gram} we obtain the bound
\[
\dim L^n_k|_{Mo} \leq \dim L^n_k|_{TL}.
\]
In particular, in this case $b_n$ and the RepGap are exponentially smaller than for $Mo_{1,\delta}(n)$.

\subsection{Planar rook}

The following works over an arbitrary field (in particular, we do not assume characteristic zero).

By \cite[Proposition 4F.7]{khovanov-monoidal-2024}, $pRo_{a_0}$ is semisimple for $a_0 \neq 0$, so the simple representations have the same dimensions as the cell modules in \autoref{P:PlanCellModuleDims}.

For $a_0=0$, the algebra $pRo_{a_0}$ is an outlier: it has only the trivial simple module, but it admits indecomposable modules of dimension two. Concretely (we write the argument for $n=2$; the general case is analogous), the calculation
\begin{gather*}
a=
\begin{tikzpicture}[anchorbase]
\draw[usual] (0,0) to (0.5,0.5);
\draw[usual,dot] (0.5,0) to (0.5,0.2);
\draw[usual,dot] (0,0.5) to (0,0.3);
\end{tikzpicture}
,\quad
b=
\begin{tikzpicture}[anchorbase]
\draw[usual,dot] (0,0) to (0,0.2);
\draw[usual,dot] (0.5,0.5) to (0.5,0.3);
\draw[usual,dot] (0.5,0) to (0.5,0.2);
\draw[usual,dot] (0,0.5) to (0,0.3);
\end{tikzpicture}
,\quad
a\cdot a=b,
\quad
a\cdot x=0\textnormal{ for }x\notin\{1,a\},
\end{gather*}
shows that the cyclic module generated by $a$ is two dimensional with basis $\{a,b\}$. The action matrix of $a$ on this cyclic module is
\[
[a]=
\begin{pmatrix}
0 & 0\\
1 & 0
\end{pmatrix},
\]
and all other generators act trivially. Hence this two dimensional module is indecomposable.

Thus, as for Motzkin, for $a_0=0$ we have that $b_n$ and the RepGap are exponentially smaller than for $a_0=1$.

\subsection{Brauer}

Let $\delta = m-2l$ for $l,m \in \mathbb{Z}_{\geq 0}$. Then, by \autoref{P:BrSchurWeyl}, the condition $2n < (m+1)(2l+1)$ implies
\begin{gather}\label{eq:BrIso}
Br_{\delta}(n) \cong \textnormal{End}_{OSp(V)}(V^{\otimes n}),
\end{gather}
where the $OSp(V)$-module $V$ has dimension $m+2l$ as a vector space.

Moreover, $\delta=m-2l$ implies $m=2l+\delta$, and we set $2n=(m+1)(2l+1)-1$, so that
\[
n = 2l(l+1)+\frac{\delta}{2}(2l+1).
\]
We now let $l\to\infty$, and hence also $n\to\infty$ (equivalently, we vary the orthosymplectic group with $l$). Solving for $l$ gives
\[
l = \tfrac{-(2+\delta)\pm \sqrt{(4+\delta)^2+8n-4\delta}}{4}
\;\sim\;
\frac{\sqrt{n}}{\sqrt{2}}
\qquad (\text{as } n\to\infty).
\]

Schur--Weyl duality identifies $b_n$ with the number of direct summands of $V^{\otimes n}$. For fixed $l$, this grows as $(b_n)^{1/n}\sim 4l+\delta\sim 4l$ by \cite[Theorem 1.3]{COT-GrowthSummand}. Since $l\sim \sqrt{n/2}$ along the above specialization, this forces superexponential growth in $n$.

\begin{Proposition}
For all $C\in\N$ we have
\[
C\leq (b_n)^{1/n}\quad\textnormal{ for }n\gg 0.
\]
\end{Proposition}

\begin{proof}
Combine $l\sim \sqrt{n/2}$ with $(b_n)^{1/n}\sim 4l$ from \cite{COT-GrowthSummand}.
\end{proof}

\subsection{Rook Brauer}

With \autoref{P:RoBrSchurWeyl}, we obtain the same bounds as for Brauer. Concretely, for $2d < (m+2)(2n+1)$ and $\delta = m+1-2n \neq 0$,
\[
RoBr_{\delta}(d) \cong \textnormal{End}_{OSp(V)}\left(V^{\otimes d}\right).
\]
If $\delta=1$, then $m=2n$, hence $\dim(V)=4n+1$ and the condition on $d$ becomes $2d \leq (2n+2)(2n+1)$, matching the Brauer inequalities and giving the same lower bound.

We turn now to $RoBr_{1,0}(d)$. By \cite[Theorem 1]{MM-PartialBrauerRep}, we have
\[
RoBr_{1,0}(d) \cong Br_{-1}(d)\oplus Br_{-1}(d-1),
\]
so the same superexponential lower bound follows by applying the Brauer argument to each summand.

Finally, we discuss $RoBr_{0,\delta}(n)$. Similarly to the Motzkin case, using the paragraph following \autoref{T:Gram}, we again obtain a bound on simple dimensions:
\[
\dim L^n_k(\lambda)|_{RoBr} \leq \dim L^n_k(\lambda)|_{Br}.
\]
We also expect this to be much smaller than for $a_0=1$, but we do not have fine enough results to determine this.

\subsection{Rook}

By \cite[Proposition 5F.8]{khovanov-monoidal-2024}, $Ro_{a_0}$ is semisimple for $a_0 \neq 0$, so the simple representations have the same dimensions as the cell modules in \autoref{P:SymmCellMod}. The case $a_0=0$ can be treated \textit{mutatis mutandis} as for $pRo_{a_0}$. Again, this gives that the $a_0=0$ case has much smaller dimensions.

\subsection{Partition}\label{SS:Partition}

The partition family is the main place where the ``truncation heuristic'' needs extra care: already for $\delta\in\{0,1\}$, some small $k$-types behave differently, and this affects which endpoints should be included in a typical truncation (cf.\ \autoref{R:PartExts}).

Let $\mathbb{Z}^{\infty} = \prod_{i\geq 1}\mathbb{Z}\varepsilon_i$ for a set of formal symbols $\{\varepsilon_1,\varepsilon_2,\dots \}$ and let
\[
\textnormal{Part} = \{(n,k,\lambda) \mid n \in \mathbb{Z}_{\geq 0},\, 1 \leq k \leq n,\, \lambda \vdash k \}.
\]
For $\lambda = (\lambda_1,\lambda_2,\lambda_3,\dots) \vdash k$ (with infinitely many trailing zeros), define
\[
\varphi_{\delta}(n,k,\lambda) =
\left\{
\begin{array}{ll}
(\delta-k, \lambda_1-1, \lambda_2-2, \lambda_3-3,\ldots ) & \mbox{if $n$ is even,}\\
(\delta-k-1, \lambda_1-1, \lambda_2-2, \lambda_3-3,\ldots ) & \mbox{if $n$ is odd.}
\end{array}
\right.
\]
Let $\varphi_{\delta}(n,k,\lambda)_j$ denote the $j$th component.

\begin{Lemma}\label{L:SimplePaths}
Let $(n,k,\lambda) \in \textnormal{Part}$. If $\varphi_{\delta}(n,k,\lambda)_0 = \varphi_{\delta}(n,k,\lambda)_j$ for some $j \geq 1$, then the cell module $\Delta^n_k(\lambda)$ is simple. Moreover, if $\varphi_{\delta}(n,k,\lambda)_0 \neq \varphi_{\delta}(n,k,\lambda)_j$ for all $j\geq 1$, then there exists $N_k \gg k$ such that
\[
\dim L^n_k(\lambda) < \dim\Delta^n_k(\lambda)
\]
for all $n \geq N_k$.
\end{Lemma}

\begin{proof}
The first part follows from \cite[Remark 4.5]{BdVE-SimplePartition}. For the second part, we briefly explain the mechanism from \cite{BdVE-SimplePartition}. Consider the branching graph
\[
\mathcal{Y}=\bigcup_{n\geq 0}\{(n,k,\lambda)\mid 1\leq k \leq \lfloor n/2 \rfloor,\ \lambda \vdash k\}
\]
as in \cite[Section 1]{BdVE-SimplePartition}. By \cite[Section 4]{BdVE-SimplePartition}, $\varphi_\delta$ embeds $\mathcal{Y}$ into $\mathbb{Z}^\infty$, and \cite[Definitions 1.1, 4.6, and 6.1]{BdVE-SimplePartition} define $\delta$-permissible paths in $\mathcal{Y}$. By \cite[Lemma 6.3]{BdVE-SimplePartition}, the strict inequality $\dim L^n_k(\lambda)<\dim\Delta^n_k(\lambda)$ occurs once there exists \emph{some} path ending at $(n,k,\lambda)$ that is not $\delta$-permissible. For fixed $k$, as $n$ grows there is enough ``room'' to choose such a path (one can move sufficiently far to the right before the path must terminate), and this yields the existence of $N_k$.
\end{proof}

\begin{Proposition}\label{P:PartRepGap}
If $\delta = 1$, then $\dim L^n_k((k))$ is
\[
\begin{cases}
< \dim\Delta^n_k((k)) & \textnormal{if }k=0,2, \\
= \dim\Delta^n_k((k)) & \textnormal{if }k=1\textnormal{ or }3 \leq k \leq n,
\end{cases}
\]
and $\dim L^n_k((1^k))$ is
\[
\begin{cases}
< \dim\Delta^n_k((1^k)) & \textnormal{if }k=0, \\
= \dim\Delta^n_k((1^k)) & \textnormal{if }1 \leq k \leq n.
\end{cases}
\]
If $\delta = 0$, then $\dim L^n_k((k))$ is
\[
\begin{cases}
< \dim\Delta^n_k((k)) & \textnormal{if }k=0,1, \\
= \dim\Delta^n_k((k)) & \textnormal{if }2 \leq k \leq n,
\end{cases}
\]
and $\dim L^n_k((1^k))< \dim\Delta^n_k((1^k))$ for $0\leq k \leq n$.
\end{Proposition}

\begin{proof}
Let $(n,k,\lambda) \in \textnormal{Part}$ with $\lambda \in \{(k),(1^k)\}$, corresponding to the trivial and sign representations of $S(k)$. We compute
\[
\varphi_{\delta}(n,k,(k)) = (\delta - k, k-1, -2, -3, \dots ),
\qquad
\varphi_{\delta}(n,k,(1^k)) = (\delta - k, 0, -1, -2, \dots , 1-k, -1-k, \dots).
\]
For $\delta=1$, one checks:
\begin{itemize}
\item $\lambda=(k)$: $k=1$ gives $\varphi_1(n,1,(1))_0 = \varphi_1(n,1,(1))_1$, $k\geq 3$ gives $\varphi_1(n,k,(k))_0 = \varphi_1(n,k,(k))_{k-1}$, while $k\in\{0,2\}$ gives no equality with any component $\varphi_1(n,k,(k))_j$ for $j\geq 1$.
\item $\lambda=(1^k)$: $k\geq 1$ gives $\varphi_1(n,k,(1^k))_0 = \varphi_1(n,k,(1^k))_k$, while $k=0$ gives no equality.
\end{itemize}
The arguments for $\delta=0$ are analogous, and the proposition follows from \autoref{L:SimplePaths}.
\end{proof}

\begin{Remark}
Since the $S(k)$-representations corresponding to $\lambda \in \{(k),(1^k)\}$ are both $1$ dimensional, we have $\dim\Delta^n_k((k))=\dim\Delta^n_k((1^k))$ by \autoref{P:SymmCellMod}. Therefore, \autoref{P:PartRepGap} shows that the simples of type $\lambda=(1^k)$ for $Pa_1(n)$ are (for $n$ large enough) strictly larger than those for $Pa_0(n)$, suggesting that the RepGap may be larger as well.

Whilst \cite[Lemma 6.3]{BdVE-SimplePartition} provides a method to compute dimensions in the remaining blocks (i.e.\ when $L^n_k(\lambda)$ is not a cell module), it is difficult to use in practice. We therefore record the clean consequences for the one dimensional $S(k)$-types, and treat the general case as future work.
\end{Remark}

We now return to unfinished business from \autoref{R:PartExts}, namely the interaction of typical truncation with the exceptional small values of $k$.

Using Mathematica, we can compute that the maximum of $\dim\Delta^n_k((k))$ occurs approximately at
\[
k=(0.1828-0.012\log_e(n))n.
\]
In particular, this maximum lies at $k\geq 3$ for all $n\geq 14$. Therefore, for any fixed $C>0$ and all $n$ large enough, the typical truncation $Pa^{\mathrm{typ}}_{1,C}(n)$ only sees values of $k$ with $k\geq 3$, and hence \autoref{P:PartRepGap} implies that the trivial and sign $S(k)$-types contribute with their \emph{full} cell-module dimensions.

\begin{Proposition}
Fix $C\in\R_{>0}$. For $n\gg 0$ and $\lambda \in \{(k),(1^k)\}$ occurring in the typical truncation $Pa^{\mathrm{typ}}_{1,C}(n)$, we have
\[
\dim L^n_k(\lambda)=\dim\Delta^n_k(\lambda)
=\sum_{t=0}^{n}\begin{Bmatrix}n\\t\end{Bmatrix}\binom{t}{k}.
\]
Moreover, if we sum these dimensions over all $k$ (equivalently, up to subexponential error, over the typical window), then
\[
\left(\sum_{k=0}^n \dim L^n_k(\lambda)\right)^{1/n}
\sim
\frac{2n}{e\log_e(2n)}.
\]
\end{Proposition}

\begin{proof}
The dimension formula follows from \autoref{P:PartRepGap} together with \autoref{P:SymmCellMod}. The stated asymptotic is computed in \cite{GT-GrowthDiagCat}; by Gaussian concentration, passing to the typical window does not change the $n$th-root asymptotics.
\end{proof}

\begin{Remark}
We have not confirmed that the contribution of $\lambda \notin \{(k),(1^k)\}$ leaves the same $n$th-root asymptotics for $Pa^{\mathrm{typ}}_{1,C}(n)$. In principle one can use \cite[Lemma 6.3]{BdVE-SimplePartition} or compute ranks of sandwich matrices to access these blocks, but both approaches appear difficult, and we leave this to future work.
\end{Remark}

\begin{Remark}
We have not discussed $Pa_{\mathbf{a}}(n)$ for $\mathbf{a} \notin \{\mathbf{1},\mathbf{0}\}$, nor $pPa_{\mathbf{a}}(n)$.

Since $pPa_1(n) \cong TL_1(2n)$, and we can use the Schur--Weyl duality from \autoref{T:SchurWeylPlan} and the methods of \cite{An-SimpleTLAll,GT-GrowthDiagCat} for $pPa_0(n)$. The resulting growth statements are therefore identical to those for $b_{2n,k,l}$ from \autoref{P:TLResult}.

At present, the only method we know to compute simple dimensions for $\mathbf{a} \notin \{\mathbf{1},\mathbf{0}\}$ is via sandwich matrices, which we again leave to future work.
\end{Remark}

\subsection{Conclusion}

The examples above illustrate three qualitatively different behaviors for nonsemisimple dimension growth.
\begin{enumerate}

\item First, for Temperley--Lieb and Motzkin (and the families reduced to these by trivial-extension arguments), nonsemisimplicity is reflected by explicit truncations/recurrences on the branching graph, but the overall growth remains of the same form as in the semisimple case: the sums of simple dimensions satisfy
\[
b_n(TL)\approx n^{-1/2}2^n,
\quad
b_n(Mo)\approx n^{-3/2}3^n,
\]
with constants depending on the root-of-unity order $l$ (cf.\ \autoref{P:TLResult} and \autoref{C:MoSimples}). So ``semisimple $\approx$ nonsemisimple''.

\item Second, for planar rook and rook at $a_0=0$, the simple spectrum is trivial, but new indecomposables appear (already in dimension $2$). This is a reminder that RepGap and $b_n$ measure different aspects of the representation theory than just semisimplicity: here $b_n$ is small, while the module category is nevertheless far from semisimple. For $a_0=1$ everything is semisimple.

\item Third, for Brauer and rook Brauer ($a_0=1$), varying the orthosymplectic parameter with $n$ forces superexponential growth of the number of summands, and hence of $b_n$ (cf.\ \autoref{eq:BrIso} and the discussion following it). This is essentially what ``semisimple $\approx$ nonsemisimple'' predicts, but we do not have fine enough results to conclude that.
For $a_0=0$ again special behavior occurs.

\item Finally, the partition family sits in the ``semisimple $\approx$ nonsemisimple'' regime: in the $\mathbf{a}=\mathbf{1}$ case for the one dimensional $S(k)$-types, nonsemisimplicity only affects finitely many small values of $k$ (cf.\ \autoref{P:PartRepGap}), so typical truncations avoid these exceptional endpoints for $n\gg 0$ and the dominant growth is governed by the same central combinatorics as in the semisimple case. For $\mathbf{a} = \mathbf{0}$, we expect approximately the same result, up to the $n$th root, despite the smaller dimensions over $\lambda = (1^k)$, although this is difficult to compute directly. Determining this and the contribution of the remaining $\lambda$ (and of general parameter choices) remains an interesting open problem.

\end{enumerate}

Altogether, this suggests the slogan that the representation theory of these generalized (parameters and twists) diagram algebras is essentially the same as the one of their classical counterparts: varying the parameters and twists mainly changes small-rank and extension phenomena, while the typical-sector asymptotic growth is essentially stable.


\newcommand{\etalchar}[1]{$^{#1}$}


\begin{thebibliography}{MMM{\etalchar{+}}21}

\bibitem[Abr96]{Abrams}
L.~Abrams.
\newblock Two dimensional topological quantum field theories and
Frobenius algebras.
\newblock {\em J. Knot Theory Ramifications}, 5 (1996), 569--587.
\newblock \href{https://doi.org/10.1142/S0218216596000333}{\path{doi:10.1142/S0218216596000333}}

\bibitem[An19]{An-SimpleTLAll}
H.H.~Andersen.
\newblock Simple modules for Temperley--Lieb algebras and related algebras.
\newblock {\em J. Algebra}, 520:276--308, 2019.
\newblock \url{https://arxiv.org/abs/1709.00248}, \href{https://doi.org/10.1016/j.jalgebra.2018.10.035}{\path{doi:10.1016/j.jalgebra.2018.10.035}}

\bibitem[APW91]{AnPoWe-representation-qalgebras}
H.H.~Andersen, P.~Polo, and K.X. Wen.
\newblock Representations of quantum algebras.
\newblock {\em Invent. Math.}, 104(1):1--59, 1991.
\newblock \href {https://doi.org/10.1007/BF01245066}
{\path{doi:10.1007/BF01245066}}.

\bibitem[AST15]{AST-CellularStructures}
H.H.~Andersen, C.~Stroppel, and D.~Tubbenhauer.
\newblock Additional notes for the paper “Cellular structures
using {$\mathrm{U}_q$}-tilting modules.”
\newblock \url{https://arxiv.org/abs/1503.00224}.

\bibitem[AST18]{AnStTu-cellular-tilting}
H.H.~Andersen, C.~Stroppel, and D.~Tubbenhauer.
\newblock Cellular structures using {$\mathrm{U}_q$}-tilting modules.
\newblock {\em Pacific J. Math.}, 292(1):21--59, 2018.
\newblock \url{https://arxiv.org/abs/1503.00224}, \href
{https://doi.org/10.2140/pjm.2018.292.21}
{\path{doi:10.2140/pjm.2018.292.21}}.

\bibitem[AST17]{AST}
H.H.~Andersen, C.~Stroppel, and D.~Tubbenhauer.
\newblock Semisimplicity of {H}ecke and (walled) {B}rauer algebras.
\newblock {\em J. Aust. Math. Soc.}, 103 (2017), no. 1, 1--44.
\newblock \url{https://arxiv.org/abs/1507.07676}, \href
{https://doi.org/10.1017/S1446788716000392}
{\path{doi:10.1017/S1446788716000392}}.

\bibitem[AM24]{AM-SchurRook}
C.~Andr\'{e} and I.~Martins.
\newblock Schur--Weyl dualities for the rook monoid: an approach via Schur algebras.
\newblock {\em Semigroup Forum} 109, 38--59, 2024.
\newblock \url{https://arxiv.org/abs/2404.01493}, \href{https://doi.org/10.1007/s00233-024-10434-w}{\path{doi:10.1007/s00233-024-10434-w}}.

\bibitem[AK94]{ArKo-hecke-algebra}
S.~Ariki and K.~Koike.
\newblock A {H}ecke algebra of {$(\mathbb{Z}/r\mathbb{Z})\wr\mathfrak{S}_{n}$}
and construction of its irreducible representations.
\newblock {\em Adv. Math.}, 106(2):216--243, 1994.
\newblock \href {https://doi.org/10.1006/aima.1994.1057}{\path{doi:10.1006/aima.1994.1057}}.

\bibitem[Ar25]{Ar-RepGapMotzkin}
K.~Arms.
\newblock Representation gap of the Motzkin monoid.
\newblock 2025.
\newblock \url{https://arxiv.org/abs/2510.06707}.

\bibitem[BH14]{BH-MotzkinAlgebras}
G.~Benkart and T.~Halverson.
\newblock Motzkin algebras.
\newblock {\em Eur. J. Comb.}, 473--502, 2014.
\newblock \url{https://arxiv.org/abs/1106.5277}, \href{https://doi.org/10.1016/j.ejc.2013.09.010}{\path{doi:10.1016/j.ejc.2013.09.010}} 

\bibitem[BD09]{BeDo}
D.~Benson and S.~Doty.
\newblock Schur--Weyl duality over finite fields.
\newblock {\em Arch. Math. (Basel)}, 93 (2009), no. 5, 425--435.
\newblock \url{https://arxiv.org/abs/0805.1235}, \href{https://doi.org/10.1007/s00013-009-0066-8}{\path{doi:10.1007/s00013-009-0066-8}} 

\bibitem[Bia93]{Bi-asymptotic-lie}
P.~Biane.
\newblock Estimation asymptotique des multiplicit{\'e}s dans les puissances
tensorielles d'un {$\mathfrak{g}$}-module.
\newblock {\em C. R. Acad. Sci. Paris S{\'e}r. I Math.}, 316(8):849--852, 1993.

\bibitem[BHMV95]{BHMV}
C.~Blanchet, N.~Habegger, G.~Masbaum, and P.~Vogel.
\newblock Topological quantum field theories derived from the Kauffman bracket.
\newblock {\em Topology}, 34 (1995), no. 4, 883--927.
\newblock \href{https://doi.org/10.1016/0040-9383(94)00051-4}{\path{doi:10.1016/0040-9383(94)00051-4}}.

\bibitem[BdVE19]{BdVE-SimplePartition}
C.~Bowman, M.~De Visscher, and J.~Enyang.
\newblock Simple modules for the partition algebra and monotone convergence of Kronecker coefficients. 
\newblock {\em Int. Math. Res. Not. IMRN} 2019, no. 4, 1059--1097.
\newblock \url{https://arxiv.org/abs/1607.08495}, 
\href{https://doi.org/10.1093/imrn/rnx095}{\path{doi:10.1093/imrn/rnx095}}

\bibitem[BDM22]{BDM}
C.~Bowman, S.~Doty and S.~Martin.
\newblock Integral Schur--Weyl duality for partition algebras.
\newblock {\em Algebr. Comb.}, 5 (2022), no. 2, 371--399.
\newblock \url{https://arxiv.org/abs/1906.00457}, \href{https://doi.org/10.5802/alco.214}{\path{doi:10.5802/alco.214}}.

\bibitem[Bra37]{Br-brauer-algebra-original}
R.~Brauer.
\newblock On algebras which are connected with the semisimple continuous
groups.
\newblock {\em Ann. of Math. (2)}, 38(4):857--872, 1937.
\newblock \href {https://doi.org/10.2307/1968843} {\path{doi:10.2307/1968843}}.

\bibitem[BM93]{BrMa-hecke}
M.~Brou{\'e} and G.~Malle.
\newblock Zyklotomische {H}eckealgebren.
\newblock {\em Ast{\'e}risque}, no.212, 119--189, 1993.
\newblock \url{http://www.numdam.org/item/AST_1993__212__119_0/}.

\bibitem[Bro55]{Br-gen-matrix-algebras}
W.P.~Brown.
\newblock Generalized matrix algebras.
\newblock {\em Canadian J. Math.}, 7:188--190, 1955.
\newblock \href {https://doi.org/10.4153/CJM-1955-023-2}
{\path{doi:10.4153/CJM-1955-023-2}}.

\bibitem[Ch87]{Ch-gelfandtzetlin}
I.V.~Cherednik.
\newblock A new interpretation of {G}elfand--{T}zetlin bases.
\newblock {\em Duke Math. J.}, 54(2):563--577, 1987.
\newblock \href {https://doi.org/10.1215/S0012-7094-87-05423-8}{\path{doi:10.1215/S0012-7094-87-05423-8}}.

\bibitem[Com20]{Comes}
J.~Comes.
\newblock Jellyfish partition categories.
\newblock {\em Algebr. Represent. Theory} 23 (2020), no. 2, 327--347.
\newblock \url{https://arxiv.org/abs/1612.05182}, \href{https://doi.org/10.1007/s10468-018-09851-7}{\path{doi:10.1007/s10468-018-09851-7}}.

\bibitem[CO11]{CO}
J.~Comes and V.~Ostrik.
\newblock On blocks of {D}eligne's category {$\underline{\rm Re}{\rm p}(S_t)$}.
\newblock {\em Adv. Math.} 226 (2011), no. 2, 1331--1377.
\newblock \url{https://arxiv.org/abs/0910.5695}, \href{https://doi.org/10.1016/j.aim.2010.08.010}{\path{doi:10.1016/j.aim.2010.08.010}}.

\bibitem[CEKO23]{CEKO-InvPosChar}
K.~Coulembier, P.~Etingof, A.~Kleshchev, V.~Ostrik.
\newblock Super invariant theory in positive characteristic.
\newblock {\em Eur. J. Math.} 9 (2023), no. 4, Paper No. 94, 39 pp.
\newblock \url{https://arxiv.org/abs/2211.11933}, \href{https://doi.org/10.1007/s40879-023-00688-z}{\path{doi:10.1007/s40879-023-00688-z}}.

\bibitem[COT24]{COT-GrowthSummand}
K.~Coulembier, V.~Ostrik, and D.~Tubbenhauer.
\newblock Growth rates of the number of indecomposable summands in tensor powers.
\newblock {\em Algebr. Represent. Theory} 27, no. 2, 1033-1062, 2024.
\newblock \url{https://arxiv.org/abs/2301.00885}, \href{https://doi.org/10.1007/s10468-023-10245-7}{\path{doi:10.1007/s10468-023-10245-7}}.

\bibitem[Del07]{De-cat-st}
P.~Deligne.
\newblock La cat{\'e}gorie des repr{\'e}sentations du groupe sym{\'e}trique
{$S_t$}, lorsque {$t$} n'est pas un entier naturel.
\newblock In {\em Algebraic groups and homogeneous spaces}, volume~19 of {\em
Tata Inst. Fund. Res. Stud. Math.}, pages 209--273. Tata Inst. Fund. Res.,
Mumbai, 2007.
\newblock \url{https://www.math.ias.edu/files/deligne/Symetrique.pdf}.

\bibitem[Dij89]{Dijkgraaf}
R.~Dijkgraaf.
\newblock A geometric approach to two dimensional conformal field
theory.
\newblock Ph.D. Thesis, University of Utrecht, 1989.
\newblock \url{https://dspace.library.uu.nl/bitstream/handle/1874/210872/A%20geometrical%20approach%20to%20two-dimensional%20conformal%20field%20theory.pdf}.

\bibitem[Don93]{Do-tilting-alg-groups}
S.~Donkin.
\newblock On tilting modules for algebraic groups.
\newblock {\em Math. Z.}, 212(1):39--60, 1993.
\newblock \href {https://doi.org/10.1007/BF02571640}
{\path{doi:10.1007/BF02571640}}.

\bibitem[DPS98]{DuPaSc-Schur--Weyl}
J.~Du, B.~Parshall, and L.~Scott.
\newblock Quantum {W}eyl reciprocity and tilting modules.
\newblock {\em Comm. Math. Phys.}, 195(2):321--352, 1998.
\newblock \href {https://doi.org/10.1007/s002200050392}
{\path{doi:10.1007/s002200050392}}.

\bibitem[EGMR26]{east2025twisted}
J.~East, R.D.~Gray, P.A.~Muhammed, and N.~Ru{\v{s}}kuc.
\newblock Twisted products of monoids.
\newblock {\em J. Algebra}, 689:819--861, 2026.
\newblock \url{https://arxiv.org/abs/2507.04486}, \href{https://doi.org/10.1016/j.jalgebra.2025.10.030}{\path{doi:10.1016/j.jalgebra.2025.10.030}}.

\bibitem[EGNO15]{egno-tensor-2015}
P.~Etingof, S.~Gelaki, D.~Nikshych, and V.~Ostrik.
\newblock Tensor Categories.
\newblock {\em American Mathematical Society: Mathematical Surveys and Monographs}, Online Ed. 2331-7159, v. 205, 2015.
\newblock \href{https://doi.org/10.1090/surv/205}{\path{doi:10.1090/surv/205}}.

\bibitem[FG95]{FG}
S.~Fishel and I.~Grojnowski.
\newblock Canonical bases for the Brauer centralizer algebra.
\newblock {\em Math. Res. Lett.}, 2, no. 1, 15--26, 1995.
\newblock \href{http://dx.doi.org/10.4310/MRL.1995.v2.n1.a3}{\path{doi:10.4310/MRL.1995.v2.n1.a3}}.

\bibitem[Fu08]{Fu}
J.~Fulman.
\newblock Convergence rates of random walk on irreducible representations of finite groups.
\newblock {\em J. Theoret. Probab.}, 21 (2008), no. 1, 193--211.
\newblock \url{https://arxiv.org/abs/math/0607399}, \href {https://doi.org/10.1007/s10959-007-0102-1}
{\path{doi:10.1007/s10959-007-0102-1}}.

\bibitem[GMS09]{GaMaSt-irreps-semigroups}
O.~Ganyushkin, V.~Mazorchuk, and B.~Steinberg.
\newblock On the irreducible representations of a finite semigroup.
\newblock {\em Proc. Amer. Math. Soc.}, 137(11):3585--3592, 2009.
\newblock URL: \url{https://arxiv.org/abs/0712.2076}, \href
{https://doi.org/10.1090/S0002-9939-09-09857-8}
{\path{doi:10.1090/S0002-9939-09-09857-8}}.

\bibitem[GL96]{GrLe-cellular}
J.J.~Graham and G.~Lehrer.
\newblock Cellular algebras.
\newblock {\em Invent. Math.}, 123(1):1--34, 1996.
\newblock \href {https://doi.org/10.1007/BF01232365}
{\path{doi:10.1007/BF01232365}}.

\bibitem[Gr51]{Gr-structure-semigroups}
J.A.~Green.
\newblock On the structure of semigroups.
\newblock {\em Ann. of Math. (2)}, 54:163--172, 1951.
\newblock \href {https://doi.org/10.2307/1969317} {\path{doi:10.2307/1969317}}.

\bibitem[Gr98]{Gr-gen-tl-algebra}
R.M.~Green.
\newblock Generalized {T}emperley--{L}ieb algebras and decorated tangles.
\newblock {\em J. Knot Theory Ramifications}, 7(2):155--171, 1998.
\newblock \url{https://arxiv.org/abs/q-alg/9712018}, \href
{https://doi.org/10.1142/S0218216598000103}{\path{doi:10.1142/S0218216598000103}}.

\bibitem[GT25]{GT-GrowthDiagCat}
J.~Gruber and D.~Tubbenhauer.
\newblock Growth problems in diagram categories.
\newblock {\em Bull. London Math. Soc.}, 57: 3454--3469, 2025.
\newblock \url{https://arxiv.org/abs/2503.00685}, \href{https://doi.org/10.1112/blms.70163}{\path{doi:10.1112/blms.70163}}.

\bibitem[GW15]{GuWi-almost-cellular}
N.~Guay and S.~Wilcox.
\newblock Almost cellular algebras.
\newblock {\em J. Pure Appl. Algebra}, 219(9):4105--4116, 2015.
\newblock \href {https://doi.org/10.1016/j.jpaa.2015.02.010}{\path{doi:10.1016/j.jpaa.2015.02.010}}.

\bibitem[HR05]{HaRa-partition-algebras}
T.~Halverson and A.~Ram.
\newblock Partition algebras.
\newblock {\em European J. Combin.}, 26(6):869--921, 2005.
\newblock \url{https://arxiv.org/abs/math/0401314}, \href
{https://doi.org/10.1016/j.ejc.2004.06.005}
{\path{doi:10.1016/j.ejc.2004.06.005}}.

\bibitem[HT25]{HeTu}
D.~He and D.~Tubbenhauer.
\newblock Affine diagram categories, algebras and monoids.
\newblock 2025.
\newblock \url{https://arxiv.org/abs/2512.05510}.

\bibitem[Jon83]{Jo}
V.F.R.~Jones.
\newblock Index for subfactors.
\newblock {\em Invent. Math.} 72 (1983), no. 1, 1--25.
\newblock URL: \url{https://arxiv.org/abs/1401.1774}, \href
{https://doi.org/10.1007/BF01389127} {\path{doi:10.1007/BF01389127}}.

\bibitem[Jon94]{Jones}
V.F.R.~Jones.
\newblock The Potts model and the symmetric group.
\newblock {\em Subfactors (Kyuzeso, 1993)}, 259--267.
World Scientific Publishing Co., Inc., River Edge, NJ, 1994.

\bibitem[KMY19]{KMY}
Z.~Kádár, P.P.~Martin and S.~Yu.
\newblock On geometrically defined extensions of the Temperley--Lieb category in the Brauer category.
\newblock {\em Math. Z.} 293 (2019), no. 3-4, 1247--1276.
\newblock URL: \url{https://arxiv.org/abs/1401.1774}, \href
{https://doi.org/10.1007/s00209-019-02246-4} {\path{doi:10.1007/s00209-019-02246-4}}.

\bibitem[Kan98]{kaneda}
M.~Kaneda.
\newblock Based modules and good filtrations in algebraic groups.
\newblock {\em Hiroshima Math. J.} 28 (1998), no. 2, 337--344.
\newblock \href
{https://doi.org/10.32917/hmj/1206126765}
{\path{doi:10.32917/hmj/1206126765}}.

\bibitem[KS24]{KhSa-cobordisms}
M.~Khovanov and R.~Sazdanovic.
\newblock Bilinear pairings on two dimensional cobordisms and generalizations
of the {D}eligne category.
\newblock {\em Fund. Math.}, 264(1):1--20, 2024.
\newblock URL: \url{https://arxiv.org/abs/2007.11640}, \href
{https://doi.org/10.4064/fm283-8-2023} {\path{doi:10.4064/fm283-8-2023}}.

\bibitem[KST24]{khovanov-monoidal-2024}
M.~Khovanov, M.~Sitaraman, and D.~Tubbenhauer.
\newblock Monoidal categories, representation gap and cryptography.
\newblock {\em Trans. Amer. Math. Soc. Ser. B}, 11:329--395, 2024.
\newblock \url{https://arxiv.org/abs/2201.01805}, \href
{https://doi.org/10.1090/btran/151}
{\path{doi:10.1090/btran/151}}.

\bibitem[Koc04]{Ko-tqfts}
J.~Kock.
\newblock {\em Frobenius algebras and 2{D} topological quantum field theories},
volume~59 of {\em London Mathematical Society Student Texts}.
\newblock Cambridge University Press, Cambridge, 2004.
\newblock \href
{https://doi.org/10.1017/CBO9780511615443}
{\path{doi:10.1017/CBO9780511615443}}.

\bibitem[KX12]{KoXi-affine-cellular}
S.~K{\"o}nig and C.~Xi.
\newblock Affine cellular algebras.
\newblock {\em Adv. Math.}, 229(1):139--182, 2012.
\newblock \href {https://doi.org/10.1016/j.aim.2011.08.010}{\path{doi:10.1016/j.aim.2011.08.010}}.

\bibitem[Koi89]{Ko}
K.~Koike.
\newblock On the decomposition of tensor products of the representations of classical groups: by means of universal characters.
\newblock {\em Adv. Math.}, 74 (1989), 57--86.
\newblock \href
{https://doi.org/10.1016/0001-8708(89)90004-2} {\path{doi:10.1016/0001-8708(89)90004-2}}.

\bibitem[LRMD23]{LRMD}
A.~Langlois-Rémillard and A.~Morin-Duchesne.
\newblock Uncoiled affine Temperley--Lieb algebras and their Wenzl--Jones projectors.
\newblock 2023.
\newblock \url{https://arxiv.org/abs/2302.12782}.

\bibitem[LL10]{LawLim-RW}
G.F.~Lawler and V.~Limic,
\newblock Random Walk: A Modern Introduction.
\newblock Cambridge Stud. Adv. Math., 123
Cambridge University Press, Cambridge, 2010. xii+364 pp.
\href{https://doi.org/10.1017/CBO9780511750854}{\path{doi:10.1017/CBO9780511750854}}.

\bibitem[LZ16]{LZ-FirstInvOrtho}
G.~Lehrer and R.~Zhang.
\newblock The first fundamental theorem of invariant theory for the orthosymplectic supergroup.
\newblock {\em Commun. Math. Phys.} 349, 661-702, 2016.
\newblock \url{https://arxiv.org/abs/1401.7395}, \href{https://doi.org/10.1007/s00220-016-2731-7}{\path{doi:10.1007/s00220-016-2731-7}}.

\bibitem[LZ19]{LZ-SecondInvOrtho}
G.~Lehrer and R.~Zhang.
\newblock The second fundamental theorem of invariant theory for the orthosymplectic supergroup.
\newblock {\em Nagoya Math. J.} 242 (2021), 52--76.
\newblock \url{https://arxiv.org/abs/1407.1058}, 
\href{https://doi.org/10.1017/nmj.2019.25}{\path{doi:10.1017/nmj.2019.25}}.

\bibitem[LPW09]{LPW-MC}
D.A.~Levin, Y.~Peres and E.L.~Wilmer,
\newblock Markov Chains and Mixing Times.
\newblock With a chapter by James G. Propp and David B. Wilson
American Mathematical Society, Providence, RI, 2009. xviii+371 pp.
\newblock
\href{https://doi.org/10.1090/mbk/058}{\path{doi:10.1090/mbk/058}}.

\bibitem[Liu25]{Liu25}
J.~Liu.
\newblock Representations of cyclic diagram monoids.
\newblock 2025.
\newblock \url{https://arxiv.org/abs/2511.15945}.

\bibitem[LS77]{LoganShepp}
B.F.~Logan and L.A.~Shepp.
\newblock A variational problem for random Young tableaux.
\newblock {\em Adv. Math.}, 26(2):206--222, 1977.
\newblock \url{https://doi.org/10.1016/0001-8708(77)90030-5}, \href
{https://doi.org/10.1016/0001-8708(77)90030-5} {\path{doi:10.1016/0001-8708(77)90030-5}}.

\bibitem[Mar94]{Martin}
P.~Martin.
\newblock Temperley--Lieb algebras for nonplanar statistical mechanics---the partition algebra construction.
\newblock {\em J. Knot Theory Ramifications}, 3 (1994), no. 1, 51--82.
\newblock \href{https://doi.org/10.1142/S0218216594000071}{\path{doi:10.1142/S0218216594000071}}.

\bibitem[MM14]{MM-PartialBrauerRep}
P.~Martin and V.~Mazorchuk.
\newblock On the representation theory of partial Brauer algebras.
\newblock {\em Q. J. Math.} 65 (2014), no. 1, 225--247.
\newblock \url{https://arxiv.org/pdf/1205.0464}, \href{https://doi.org/10.1093/qmath/has043}{\path{doi:10.1093/qmath/has043}}.

\bibitem[MS94]{MaSa-blob}
P.~Martin and H.~Saleur.
\newblock The blob algebra and the periodic {T}emperley--{L}ieb algebra.
\newblock {\em Lett. Math. Phys.}, 30, no.3, 189--206, 1994.
\newblock \url{https://arxiv.org/abs/hep-th/9302094}, \href
{https://doi.org/10.1007/BF00805852}{\path{doi:10.1007/BF00805852}}.

\bibitem[Mat99]{Ma-hecke-schur}
A.~Mathas.
\newblock {\em Iwahori--{H}ecke algebras and {S}chur algebras of the symmetric
group}, volume~15 of {\em University Lecture Series}.
\newblock American Mathematical Society, Providence, RI, 1999.
\newblock \href {https://doi.org/10.1090/ulect/015}
{\path{doi:10.1090/ulect/015}}.

\bibitem[Mat98]{Mat}
A.~Mathas.
\newblock Simple modules of Ariki--Koike algebras.
\newblock {\em Group representations: cohomology, group actions and topology} (Seattle, WA, 1996), 383-396.
Proc. Sympos. Pure Math., 63
American Mathematical Society, Providence, RI, 1998.
\newblock \href{https://doi.org/10.1090/pspum/063/1603195}{\path{doi:10.1090/pspum/063/1603195}}.

\bibitem[MT23]{MaTu}
A.~Mathas and D.~Tubbenhauer.
\newblock Cellularity of KLR and weighted KLRW algebras via crystals.
\newblock 2023.
\newblock To appear in Commun. Am. Math. Soc.
\newblock \url{https://arxiv.org/abs/2309.13867}.

\bibitem[{OEI}23]{oeis}
{OEIS Foundation Inc.}
\newblock The {O}n-{L}ine {E}ncyclopedia of {I}nteger {S}equences, 2023.
\newblock Published electronically at \url{http://oeis.org}.

\bibitem[OST25]{ST-GrowthQuantum}
J.~O'Sullivan and D.~Tubbenhauer.
\newblock Growth Problems of Quantum Groups.
\newblock 2025.
\newblock \url{https://arxiv.org/abs/2511.06737}.

\bibitem[Rin91]{Ri-good-filtrations}
C.M.~Ringel.
\newblock The category of modules with good filtrations over a quasi-hereditary
algebra has almost split sequences.
\newblock {\em Math. Z.}, 208(2):209--223, 1991.
\newblock \href {https://doi.org/10.1007/BF02571521}
{\path{doi:10.1007/BF02571521}}.

\bibitem[RTW32]{RTW-Valenztheorie}
G.~Rumer, E.~Teller, and H.~Weyl.
\newblock Eine für die Valenztheorie geeignete Basis der binären Vektorinvarianten.
\newblock {\em Nachrichten von der Gesellschaft der Wissenschaften zu Göttingen, Mathematisch-Physikalische Klasse}, 499--504, 1932. 
\newblock \url{http://eudml.org/doc/59396}.

\bibitem[Sc27]{Sc-ClassicSchurWeyl}
I.~Schur.
\newblock Über die rationalen Darstellungen der allgemeinen linearen Gruppe.
\newblock {\em Sitzungsberichte Akad}, Berlin, 58--75, 1927.

\bibitem[Scr24]{Scrimshaw}
T.~Scrimshaw.
\newblock Cellular subalgebras of the partition algebra.
\newblock {\em J. Comb. Algebra}, 8 (2024), no. 1-2, 147--207.
\newblock \url{https://arxiv.org/abs/2211.08746}, \href
{https://doi.org/10.4171/jca/84}
{\path{doi:10.4171/jca/84}}.

\bibitem[St16]{St-rep-monoid}
B.~Steinberg.
\newblock {\em Representation theory of finite monoids}.
\newblock Universitext, Springer, Cham, 2016, xxiv+317 pp.
\newblock \href {https://doi.org/10.1007/978-3-319-43932-7}
{\path{doi:10.1007/978-3-319-43932-7}}.

\bibitem[ST25]{ST-RepGapRigidPlanar}
W.~Stewart and D.~Tubbenhauer.
\newblock Representation gaps of rigid planar diagram monoids.
\newblock 2025.
\newblock \url{https://arxiv.org/abs/2505.05846}.

\bibitem[STWZ23]{SuTuWeZh-mixed-tilting}
L.~Sutton, D.~Tubbenhauer, P.~Wedrich, and J.~Zhu.
\newblock {SL2} tilting modules in the mixed case.
\newblock {\em Selecta Math. (N.S.)}, 29(3):39, 2023.
\newblock \url{https://arxiv.org/abs/2105.07724}, \href
{https://doi.org/10.1007/s00029-023-00835-0}
{\path{doi:10.1007/s00029-023-00835-0}}.

\bibitem[TZ04]{TZ}
T.~Tate and S.~Zelditch.
\newblock Lattice path combinatorics and asymptotics of multiplicities of weights in tensor powers.
\newblock {\em J. Funct. Anal.} 217 (2004), no.~2, 402--447.
\newblock \url{https://arxiv.org/abs/math/0305251}, \href{https://doi.org/10.1016/j.jfa.2004.01.004}{\path{doi:10.1016/j.jfa.2004.01.004}}.

\bibitem[TL71]{TL}
H.N.V.~Temperley and E.H.~Lieb.
\newblock Relations between the ``percolation'' and ``colouring'' problem and other graph-theoretical problems associated with regular planar lattices: some exact results for the ``percolation'' problem.
\newblock {\em Proc. Roy. Soc. London Ser. A} (1971) 322 (1549): 251--280.
\newblock \href{https://doi.org/10.1098/rspa.1971.0067}{\path{doi:10.1098/rspa.1971.0067}}.

\bibitem[Tu25]{Tu-qtop}
D.~Tubbenhauer.
\newblock Quantum topology without topology.
\newblock 2025.
\newblock \url{https://arxiv.org/abs/2506.18918}.

\bibitem[Tu24]{Tu-sandwich}
D.~Tubbenhauer.
\newblock Sandwich cellularity and a version of cell theory.
\newblock {\em Rocky Mountain J. Math.}, 54(6):1733--1773, 2024.
\newblock \url{https://arxiv.org/abs/2206.06678}, \href{https://doi.org/10.1216/rmj.2024.54.1733}{\path{doi:10.1216/rmj.2024.54.1733}}.

\bibitem[TV23]{TuVa-handlebody}
D.~Tubbenhauer and P.~Vaz.
\newblock Handlebody diagram algebras.
\newblock {\em Rev. Mat. Iberoam.}, 39(3):845--896, 2023.
\newblock \url{https://arxiv.org/abs/2105.07049}, \href
{https://doi.org/10.4171/rmi/1356} {\path{doi:10.4171/rmi/1356}}.

\bibitem[Tur89]{Tur}
V.G.~Turaev.
\newblock Operator invariants of tangles, and R-matrices.
\newblock {\em Izv. Akad. Nauk SSSR Ser. Mat.}, 53 (1989), no. 5, 1073--1107, 1135; translation in
{\em Math. USSR-Izv.} 35 (1990), no. 2, 411--444.
\newblock \href
{https://doi.org/10.1070/IM1990v035n02ABEH000711} {\path{doi:10.1070/IM1990v035n02ABEH000711}}.

\bibitem[VK77]{VershikKerov}
A.M. Vershik and S.V. Kerov.
\newblock Asymptotics of the {P}lancherel measure of the symmetric group and the limiting form of Young tableaux.
\newblock {\em Soviet Math. Dokl.}, 18(2):527--531, 1977.
\newblock \url{https://www.mathnet.ru/eng/dan40430}.

\end{thebibliography}
\end{document}